\documentstyle[11pt]{article}

\newtheorem{defn}{{\bf Definition}}[section]
\newtheorem{eg}[defn]{{\bf Example}}
\newtheorem{lemma}[defn]{{\bf Lemma}}
\newtheorem{prop}[defn]{{\bf Proposition}}

\newtheorem{cor}[defn]{{\bf Corollary}}
\newtheorem{remark}[defn]{{\bf Remark}}
\newtheorem{conj}[defn]{{\bf Conjecture}}
\newtheorem{qn}[defn]{{\bf Question}}

\font\bbb=msbm9 scaled\magstep1

\newcommand{\CC}{\mbox{\bbb C}}
\newcommand{\FF}{\mbox{\bbb F}}
\newcommand{\HH}{\mbox{\bbb H}}

\newcommand{\QQ}{\mbox{\bbb Q}}
\newcommand{\RR}{\mbox{\bbb R}}
\newcommand{\ZZ}{\mbox{\bbb Z}}

\font\bb=msbm7 scaled\magstep1

\newcommand{\F}{\mbox{\bb F}}

\newcommand{\intC}{C^{^{{\circ}}}}
\newcommand{\intP}{P^{^{{\circ}}}}

\newcommand{\bdC}{C^{^{{\bullet}}}}
\newcommand{\bdP}{P^{^{{\bullet}}}}

\newcommand{\intcA}{{\rm conv}(A)^{{}^{\circ}}}

\newcommand{\intcoal}{\langle \alpha\rangle^{{}^{\circ}}}
\newcommand{\intcobe}{\langle \beta\rangle^{{}^{\circ}}}
\newcommand{\intcoga}{\langle \gamma\rangle^{{}^{\circ}}}
\newcommand{\intcosi}{\langle \sigma\rangle^{{}^{\circ}}}
\newcommand{\intcomu}{\langle \mu\rangle^{{}^{\circ}}}

\newcommand{\TPSS}{S^{\hspace{.2mm}2}\! \times \hspace{-3.3mm}_{-} \,
S^{\hspace{.1mm}1}}
\newcommand{\TPSSF}{S^{\hspace{.2mm}3}\! \times \hspace{-3.8mm}_{-} \,
S^{\hspace{.1mm}1}}
\newcommand{\TPSSD}{S^{\hspace{.2mm}d-1}\! \times \hspace{-3.8mm}_{-} \, S^{\hspace{.1mm}1}}
\newcommand{\TPSSDS}{S^{\hspace{.2mm}d-1}\! \times \hspace{-3.3mm}_{-} \, S^{\hspace{.1mm}1}}

\newcommand{\hoop}{-\hspace{-2mm}-\hspace{-1.5mm}\circ}
\newcommand{\hoops}{-\hspace{-2mm}-\hspace{-1.2mm}\circ}

\begin{document}

\title{\bf On stellated spheres, shellable balls, lower bounds and
a combinatorial criterion for tightness }
\author{{\bf Bhaskar Bagchi}$^{\rm \,a}$, {\bf Basudeb
Datta}$^{\rm \,b,}$\footnote{Corresponding author. \newline
\mbox{} \hspace{3mm} {\em E-mail addresses:} bbagchi@isibang.ac.in
(B. Bagchi), dattab@math.iisc.ernet.in (B. Datta).} $^{, 1}$
}

\date{}

\maketitle

\vspace{-4mm}

\noindent {\small $^{\rm a}$Theoretical Statistics and Mathematics
Unit, Indian Statistical Institute,  Bangalore 560\,059, India.

\smallskip

\noindent $^{\rm b}$Department of Mathematics, Indian Institute of
Science, Bangalore 560\,012,  India.
}

\footnotetext[1]{Research partially supported by grant from UGC
Centre for Advanced Study. }

\begin{center}

\date{January 28, 2012}

\end{center}

\hrule

\medskip

 \noindent {\bf Abstract}

\medskip

{\small We introduce the $k$-stellated (combinatorial) spheres and
compare and contrast them with $k$-stacked (triangulated) spheres.
It is shown that for $d \geq 2k$, any $k$-stellated sphere of
dimension $d$ bounds a unique and canonically defined $k$-stacked
ball. In parallel, any $k$-stacked polytopal sphere of dimension
$d\geq 2k$ bounds a unique and canonically defined $k$-stacked
(polytopal) ball, which answers a question of McMullen. We
consider the class ${\cal W}_k(d)$ of combinatorial $d$-manifolds
with $k$-stellated links.  For $d\geq 2k+2$, any member of ${\cal
W}_k(d)$ bounds a unique and canonically defined ``$k$-stacked"
$(d+1)$-manifold.

We introduce the mu-vector of simplicial complexes, and show that
the mu-vector of any 2-neighbourly simplicial complex dominates
its vector of Betti numbers componentwise, and the two vectors are
equal precisely when the complex is tight. When $d\geq 2k$, we
are able to estimate/compute certain alternating sums of the
mu-numbers of any 2-neighbourly member of ${\cal W}_k(d)$.
This leads to a lower bound theorem for such triangulated manifolds.
As an application, it is shown that any $(k+1)$-neighbourly
member of ${\cal W}_k(d)$ is tight, subject only to an extra
condition on the $k^{\rm th}$ Betti number in case $d=2k+1$. This
result more or less settles a recent conjecture of Effenberger,
and it also provides a uniform and conceptual tightness proof for
all the known tight triangulated manifolds, with only two
exceptions. It is shown that any polytopal upper bound sphere of
odd dimension $2k+1$ belongs to the class ${\cal W}_k(2k+1)$, thus
generalizing a theorem (the $k=1$ case) due to Perles. This shows
that the case $d=2k+1$ is indeed exceptional for the tightness
theorem.

We also prove a lower bound theorem for triangulated manifolds in
which the members of ${\cal W}_1(d)$ provide the equality case. This
generalises a result (the $d=4$ case) due to Walkup and K\"{u}hnel.
As a consequence, it is shown that every tight member of ${\cal W}_1(d)$
is strongly minimal, thus providing substantial evidence in favour of
a conjecture of K\"{u}hnel and Lutz.}

\bigskip

{\small


\medskip

 }
\bigskip

\hrule

\newpage

\tableofcontents

\section{Summary of results}

But for some exceptions in Section 3, all simplicial complexes
considered here are finite and abstract. By a triangulated
sphere/ball/manifold, we mean an abstract simplicial complex whose
geometric carrier is a sphere/ball/manifold. We identify two
complexes if they are isomorphic.

In this paper, we introduce the class $\Sigma_k(d)$, $0\leq k \leq
d +1$, of {$k$-stellated} triangulated $d$-spheres and compare it
with the class ${\cal S}_k(d)$, $0\leq k \leq d$, of {$k$-stacked}
triangulated $d$-spheres. We have the filtration
$$
\Sigma_0(d) \subseteq \Sigma_1(d) \subseteq \cdots \subseteq
\Sigma_{d}(d) \subseteq \Sigma_{d +1}(d)
$$
of the class of all combinatorial $d$-spheres, and the comparable
filtration
$$
{\cal S}_0(d) \subseteq {\cal S}_1(d) \subseteq \cdots \subseteq
{\cal S}_{d}(d)
$$
of the class of all triangulated $d$-spheres. The {\em standard
$d$-sphere} $S^{\,d}_{d + 2}$ is the unique $(d + 2)$-vertex
triangulation of the $d$-sphere. It may be described as the
boundary complex of the $(d + 1)$-dimensional geometric simplex.
The standard sphere $S^{\,d}_{d + 2}$ is the unique member of
$\Sigma_0(d) = {\cal S}_0(d)$. We also have the equality
$\Sigma_1(d) = {\cal S}_1(d)$. In the existing literature, the
members of ${\cal S}_1(d)$ are known as the $d$-dimensional {\em
stacked spheres}. For $d\geq 2k - 1$, we have the inclusion
$\Sigma_k(d) \subseteq {\cal S}_k(d)$. However, for each $k \geq
2$, there are $k$-stacked spheres which are not $k$-stellated.

In parallel with these classes of triangulated spheres, we also
consider the classes $\widehat{\Sigma}_k(d)$ and $\widehat{\cal
S}_k(d)$ of {$k$-shelled} $d$-balls and {$k$-stacked} $d$-balls,
respectively. We have the filtration
$$
\widehat{\Sigma}_0(d) \subseteq \widehat{\Sigma}_1(d) \subseteq
\cdots \subseteq \widehat{\Sigma}_{d}(d)
$$
of the class of all shellable $d$-balls, and the comparable
filtration
$$
\widehat{\cal S}_0(d) \subseteq \widehat{\cal S}_1(d) \subseteq
\cdots \subseteq \widehat{\cal S}_{d}(d)
$$
of the class of all triangulated $d$-balls. The {\em standard
$d$-ball} $B^{\,d}_{d + 1}$ is the unique $(d + 1)$-vertex
triangulation of the $d$-dimensional ball. It may be described as
the face complex of the $d$-dimensional geometric simplex. The
standard ball $B^{\,d}_{d + 1}$ is the unique member of
$\widehat{\Sigma}_0(d)= \widehat{\cal S}_0(d)$. We also have the
equality $\widehat{\Sigma}_1(d) = \widehat{\cal S}_1(d)$ and for
all $d \geq k$ we have the inclusion $\widehat{\Sigma}_k(d)
\subseteq \widehat{\cal S}_k(d)$. However, for each $k \geq 2$,
there are $k$-stacked balls which are not $k$-shelled.

While a {\em $k$-stellated $d$-sphere} is defined as a
triangulated $d$-sphere which may be obtained from $S^{\,d}_{d +
2}$ by a finite sequence of bistellar moves of index $< k$, a {\em
$k$-shelled $d$-ball} is a triangulated $d$-ball obtained from
$B^{\,d}_{d+1}$ by a finite sequence of shelling moves of index $<
k$. A {\em $k$-stacked $d$-ball} is a triangulated $d$-ball all
whose faces of codimension $k+1$ (i.e., dimension $d-k-1$) are in
its boundary. A {\em $k$-stacked $d$-sphere} is a triangulated
$d$-sphere which may be represented as the boundary of a
$k$-stacked $(d+1)$-ball. The boundary of any $k$-shelled
$(d+1)$-ball is a $k$-stellated $d$-sphere. Conversely, when
$d\geq 2k-1$, any $k$-stellated $d$-sphere may be represented as
the boundary of a $k$-shelled $(d+1)$-ball. A triangulated ball is
$k$-shelled if and only if it is $k$-stacked and shellable. Each
$k$-stacked (respectively $k$-shelled) ball is the antistar of a
vertex in a $k$-stacked (respectively $k$-stellated) sphere. We
prove that, when $d\geq 2k$, for any $k$-stellated $d$-sphere $S$,
there is a unique $k$-stacked $(d+1)$-ball $\overline{S}$ whose
boundary is $S$. The ball $\overline{S}$ has a natural and
intrinsic description in terms of the combinatorics of $S$. We
show that this result is also valid if $d\geq 2k$ and $S$ is a
polytopal $k$-stacked $d$-sphere, thus answering a query implicit
in \cite{mc2} raised by McMullen in the context of equality in
GLBT (generalized lower bound theorem) for polytopal spheres. For
general $k$-stacked spheres, we can only prove that, when $d\geq
2k+1$, a $k$-stacked $d$-sphere $S$ bounds a unique $k$-stacked
ball $\overline{S}$. However, there seems to be no combinatorial
description of $\overline{S}$ in this generality.

The entire $g$-vector (equivalently, $f$-vector) of a
$k$-stellated $d$-sphere is determined by the $k$ numbers $g_1,
\dots, g_k$. This is actually true, more generally, of $k$-stacked
$d$-spheres. However, for a $k$-stellated sphere of dimension $d
\geq 2k-1$, these $k$ components of its $g$-vector have an
interesting geometric interpretation. For $1 \leq i\leq k$, $g_i$
is the number of bistellar moves of index $i-1$ in any sequence of
bistellar moves (of index $\leq k-1$) used to obtain the given
$d$-sphere from the standard $d$-sphere.

Next, we introduce the class ${\cal W}_k(d)$, $0\leq k \leq d$, of
(combinatorial) $d$-manifolds whose vertex-links are $k$-stellated
spheres. This class may be compared with the generalized Walkup
classes ${\cal K}_k(d)$ of triangulated $d$-manifolds all whose
vertex-links are $k$-stacked spheres. We have the inclusions
$\Sigma_k(d) \subseteq {\cal W}_k(d)$, ${\cal S}_k(d) \subseteq
{\cal K}_k(d)$ and, for $d \geq 2k$, ${\cal W}_k(d) \subseteq
{\cal K}_k(d)$. In consequence, all the proper face-links (of
dimension $\geq k-1$) of members of ${\cal W}_k(d)$ are
$k$-stellated, and all the proper face-links (of dimension $\geq
k$) of members of ${\cal K}_k(d)$ are $k$-stacked. We have the
filtration
$$
{\cal W}_0(d) \subseteq {\cal W}_1(d) \subseteq \cdots \subseteq
{\cal W}_{d}(d)
$$
of the class of all closed combinatorial $d$-manifolds, and
the corresponding filtration
$$
{\cal K}_0(d) \subseteq {\cal K}_1(d) \subseteq \cdots \subseteq
{\cal K}_{d-1}(d)
$$
of the class of all triangulated closed $d$-manifolds all
whose vertex-links are triangulated $(d-1)$-spheres. Again, the
standard sphere $S^{\,d}_{d + 2}$ is the only member of ${\cal
W}_0(d) = {\cal K}_0(d)$, and we have ${\cal W}_1(d) = {\cal
K}_1(d)$. For $d\geq 2k+2$, any member $M$ of ${\cal W}_k(d)$ is
the boundary of a canonically defined $(d+1)$-manifold
$\overline{M}$ such that all the faces of codimension $k+1$ in
$\overline{M}$ belong to the boundary $M$.

The entire $g$-vector (equivalently,
face-vector) of any member of ${\cal W}_k(d)$ is determined by the
$(k+1)$ numbers $g_1, \dots, g_{k+1}$. In particular, when $d\geq
2k$ is even, the Euler characteristic $\chi$ of any member of
${\cal W}_k(d)$ is determined by the $(k+1)^{\rm st}$ component
$g_{k+1}$ of its $g$-vector by the formula
$$
(-1)^k{d+2 \choose k+1}(\chi - 2) =  2\,g_{k+1}.
$$

If $M \in {\cal W}_k(d)$ is 2-neighbourly, with
Betti numbers $\beta_i$ (with respect to any field) and
$g$-components $g_i$, then we show that
\begin{enumerate}
\item[{\rm (a)}] $g_{l+1} \geq {\displaystyle {d+2
\choose l+1} \sum_{i=1}^{l}} (-1)^{l-i}\beta_i$ \, for \, $1 \leq
l < k$, \, when $d\geq 2k$, \vspace{-2mm}
\item[{\rm (b)}] $g_{k+1} \geq
{\displaystyle {d+2 \choose k+1} \sum_{i=1}^{k}}
(-1)^{k-i}\beta_i$, \, when $d = 2k+1$,
\item[{\rm (c)}] $g_{k+1} =
{\displaystyle {d+2 \choose k+1} \sum_{i=1}^{k}}
(-1)^{k-i}\beta_i$, \, when $d\geq 2k+2$, and
\item[{\rm (d)}] $\beta_i = 0$ \, for \, $k+1 \leq i \leq d-k-1$,
\, when $d \geq 2k+2$.
\end{enumerate}
Since the components of the face vector are non-negative
linear combinations of the $g$-numbers, this result may be
interpreted as a lower bound theorem for 2-neighbourly
members of ${\cal W}_k(d)$.

We also prove a lower bound theorem for general triangulated
closed manifolds. When $d\geq 4$, among all triangulated
closed $d$-manifolds with given first Betti number and given
number of vertices, members of ${\cal W}_1(d)$ (when they exist)
minimize the face vector componentwise. The case $d=4$ of this
result is due to Walkup and K\"{u}hnel.

Any member of $\Sigma_k(d)$ (or even of ${\cal S}_k(d)$) is at
most $k$-neighbourly, unless it is the standard sphere. In
consequence, any member of ${\cal W}_k(d)$ (or ${\cal K}_k(d)$),
other than the standard sphere, is at most $(k+1)$-neighbourly.
This leads us to consider the class ${\cal W}_k^{\ast}(d)$ (and
${\cal K}_k^{\ast}(d)$) consisting of all $(k+1)$-neighbourly
members of ${\cal W}_k(d)$ (respectively  of ${\cal K}_k(d)$). We
have ${\cal W}_k^{\ast}(d) \subseteq {\cal K}_k^{\ast}(d)$ for all
$k$ and $d$. These classes have no member other than the standard
spheres unless $d \geq 2k$. Generalizing a result of Perles, we
prove that any upper bound polytopal sphere (for instance, cyclic
sphere) of odd dimension $2k+1$ belongs to the class ${\cal
W}^{\ast}_k(2k+1)$.

Moreover, when $d \geq 2k+2$ and $k \geq 2$, any member of
${\cal W}^{\ast}_k(d)$ has the same integral homology as the
connected sum of $\beta$ copies of $S^{\,k} \times S^{\,d-k}$,
where the non-negative integer $\beta$ is given by the formula
$$
{{m+k-d-2 \choose k+1}} = {{d+2 \choose k+1}}\beta,
$$
with $m = f_0$, the number of vertices. This result may be
compared with Kalai's theorem\,: for $d \geq 4$, any member of
${\cal W}_1(d)$ triangulates the connected sum of finitely many
copies of $S^1 \times S^{\,d-1}$ or $\TPSSD$.

Recall that a connected simplicial complex $X$ is said to be {\em
tight} with respect to a field $\FF$ if the inclusion map from any
induced (full) subcomplex of $X$ into $X$ is injective at the
level of $\FF$-homology. In case of a closed manifold $X$, this
has the following geometric interpretation\,: $X$ is $\FF$-tight
if the standard geometric realization of $X$ in $\RR^{\,n-1}$
($n=$ number of vertices of $X$) is ``as convex as possible"
subject to the constraint imposed by its homology with
$\FF$-coefficients. Our interest in the notion of tightness stems
from the following conjecture of K\"uhnel and Lutz (which seems to
be borne out by all the known examples)\,: any tight triangulated
manifold has the componentwise minimum face vector among all
triangulations of the same manifold!
As a consequence of our lower bound theorem for general triangulated
manifolds, we show that this conjecture is valid for all tight
members of ${\cal W}_1(d)$. Since any connected
induced subcomplex of a tight simplicial complex is obviously
tight, tightness imposes an extremely powerful constraint on the
possible combinatorics of a simplicial complex. For instance, any
$\FF$-tight simplicial complex is necessarily 2-neighbourly and
any $\FF$-tight triangulated closed manifold is
$\FF$-orientable. Thus, it is not surprising that, apart from
three infinite families (including the trivial family of standard
spheres), only seventeen sporadic examples of tight triangulated
manifolds (of dimension $>2$) are known so far.

All this makes it very important to obtain usable combinatorial
criteria for tightness of triangulated manifolds. In this paper,
we introduce the mu-vector of a simplicial complex (with respect
to a given field) and compare it with its beta-vector (i.e., the
vector of Betti numbers over the same field). It is shown that, in
general, for 2-neighbourly simplicial complexes, the alternating
sums of the components of the mu-vector dominate the corresponding
sums for the beta-vector; the two vectors coincide precisely   in
the case of tight complexes. This paraphrases the ``combinatorial
strong Morse inequality". We succeed in explicitly computing or
estimating certain functionals of the mu-vectors of 2-neighbourly
members of ${\cal W}_k(d)$, $d\geq 2k$,  entirely in terms of
their $g$-vectors. The lower bound theorem for ${\cal W}_k(d)$ stated
above (as well as the determination of the homology type of
members of ${\cal W}_k^{\ast}(d)$) is a consequence of this theory.
It also leads to the following
combinatorial criterion for tightness\,: for $d \neq 2k+1$, any
member $M$ of ${\cal W}^{\ast}_k(d)$ is tight. This is with
respect to any field in case $k \geq 2$, and with respect to a
field $\FF$ for which $M$ is $\FF$-orientable in case $k = 1$.
(The example of upper bound polytopal spheres of odd dimension
shows that the case $d=2k+1$ is a genuine exception to this
result. We also find a characterization of the tight members of
${\cal W}^{\ast}_k(2k+1)$, thus covering this exceptional case.)
The $k=1$ case of this theorem is a recent result due to
Effenberger. Effenberger also conjectured the tightness of all
members of the supposedly larger class ${\cal K}^{\ast}_k(d)$.
Thus this paper, which is largely motivated by Effenberger's work,
partially settles his conjecture. It may also be pointed
out that we do not know of a single member of ${\cal
K}^{\ast}_k(d)$ which is not in ${\cal W}^{\ast}_k(d)$.

In the final section of this paper, we present various examples,
counter examples, questions and conjectures related to the above
results. For instance, we show that for each $k \geq 2$, there are
$k$-stacked triangulated $d$-spheres which are not even
$(d+1)$-stellated (i.e., not combinatorial spheres) and
$k$-stacked combinatorial $d$-spheres which are not $d$-stellated.
Recently, Klee and Novik found an extremely beautiful construction
of a $(2d+4)$-vertex triangulation $M$ of $S^{\,k} \times
S^{\,d-k}$ for all pairs $0\leq k\leq d$. We show that, for $d\geq
2k$, these triangulations are in ${\cal W}_k(d)$. Klee and Novik
obtained their triangulation $M$ as the boundary complex of a
triangulated $(d+ 1)$-manifold $\overline{M}$. For
$d \geq 2k+2$, this is an instance of our canonical construction
$M \mapsto \overline{M}$. As an application, we show that, for
$d\neq 2k$, the full automorphism group of the Klee-Novik
triangulation is a group of order $4d+8$, already found by these
authors. This makes it interesting to determine the full
automorphism group of the Klee-Novik manifolds for $d=2k$.

We show that the tightness of most of the known tight manifolds
follows from our result. This provides a unified and conceptual
proof of tightness of these manifolds, where the previous proofs
were mostly by computer-aided case by case analysis.

In view of our (rather isolated) result on polytopal spheres, it
seems natural to conjecture that for polytopal spheres of
dimension $d\geq 2k$, the notions ``$k$-stellated" and
``$k$-stacked" coincide. We also pose a general lower bound
conjecture to which members of ${\cal W}_k(d)$ (or ${\cal
K}_k(d)$) should provide the cases of equality. This is related to
a recent work of Novik and Swartz, who proved a previous
conjecture of Kalai.

A preliminary version of this paper was posted in the arXiv
\cite{bd13_v1} 
and presented in the workshop `Topological and Geometric
Combinatorics', February 6 - 12, 2011 at MFO, Oberwolfach, Germany.

\section{Bistellar moves and shelling moves}

A $d$-dimensional simplicial complex is called {\em pure} if all
its maximal faces (called {\em facets}) are $d$-dimensional. A
$d$-dimensional pure simplicial complex is said to be a {\em weak
pseudomanifold} if each of its $(d - 1)$-faces is in at most two
facets. For a $d$-dimensional weak pseudomanifold $X$, the {\em
boundary} $\partial X$ of $X$ is the pure subcomplex of $X$ whose
facets are those $(d-1)$-dimensional faces of $X$ which are
contained in unique facets of $X$. The {\em
dual graph} $\Lambda(X)$ of a weak pseudomanifold $X$ is the graph
whose vertices are the facets of $X$, where two facets are
adjacent in $\Lambda(X)$ if they intersect in a face of
codimension one. A {\em pseudomanifold} is a weak pseudomanifold
with a connected dual graph. All connected triangulated manifolds
are automatically pseudomanifolds.

For any two simplicial complexes $X$ and $Y$, their {\em join} $X
\ast Y$ is the simplicial complex whose faces are the disjoint
unions of the faces of $X$ with the faces of $Y$. (Here we adopt
the convention that the empty set is a face of every simplicial
complex.)

For a finite set $\alpha$, let $\overline{\alpha}$ (respectively
$\partial \alpha$) denote the simplicial complex whose faces are
all the subsets (respectively, all proper subsets) of $\alpha$.
Thus, if $\#(\alpha)= n \geq 2$, $\overline{\alpha}$ is a copy of
the standard triangulation $B^{\,n - 1}_n$ of the
$(n-1)$-dimensional ball, and $\partial \alpha$ is a copy of the
standard triangulation $S^{\,n-2}_n$ of the $(n-2)$-dimensional
sphere. Thus, for any two disjoint finite sets $\alpha$ and
$\beta$, $\overline{\alpha}\ast\partial\beta$ and $\partial\alpha
\ast \overline{\beta}$ are two triangulations of a ball; they have
identical boundaries, namely $(\partial {\alpha}) \ast (\partial
{\beta})$.

A subcomplex $Y$ of a simplicial complex $X$ is said to be an {\em
induced} (or {\em full}\,) subcomplex if every face of $X$
contained in the vertex-set of $Y$ is a face of $Y$. If $X$ is a
$d$-dimensional simplicial complex with an induced subcomplex
$\overline{\alpha} \ast \partial \beta$ ($\alpha \neq \emptyset$,
$\beta \neq \emptyset$)
of dimension $d$ (thus, $\dim(\alpha) + \dim(\beta) = d$), then $Y
:= (X\setminus (\overline{\alpha} \ast \partial \beta)) \cup
(\partial\alpha \ast \overline{\beta})$ is clearly another
triangulation of the same topological space $|X|$. In this case,
$Y$ is said to be obtained from $X$ by the {\em bistellar move}
$\alpha \mapsto \beta$. If $\dim(\beta) = i$ ($0\leq i \leq d$),
we say that $\alpha \mapsto \beta$ is a {\em bistellar move of
index $i$} (or an {\em $i$-move}, in short). Clearly, if $Y$ is
obtained from $X$ by an $i$-move $\alpha \mapsto \beta$ then $X$
is obtained from $Y$ by the (reverse) $(d- i)$-move $\beta \mapsto
\alpha$. Notice that, in case $i=0$, i.e., when $\beta$ is a
single vertex, we have $\partial \beta = \{\emptyset\}$ and hence
$\overline{\alpha} \ast \partial \beta = \overline{\alpha}$.
Therefore, our requirement that $\overline{\alpha} \ast \partial
\beta$ is the induced subcomplex of $X$ on $\alpha \sqcup \beta$
means that $\beta$ is a new vertex, not in $X$. Thus, a $0$-move
creates a new vertex, and correspondingly a $d$-move deletes an
old vertex. For $0 < i < d$, any $i$-move preserves the
vertex-set; these are sometimes called the {\em proper bistellar
moves}. For a thorough treatment of bistellar moves, see \cite{bl},
for instance.

A triangulation $X$ of a manifold is called a {\em combinatorial
manifold} if its geometric carrier $|X|$ is a piecewise linear
(pl) manifold with the pl structure induced from $X$. A
combinatorial triangulation of a sphere/ball  is called a {\em
combinatorial sphere/ball} if it induces the standard pl structure
(namely, that of the standard sphere/ball) on its geometric
carrier. Equivalently (cf. \cite{li, p}), a simplicial complex is
a combinatorial sphere (or ball) if it is obtained from a standard
sphere (respectively, a standard ball) by a finite sequence of
bistellar moves. In general, a triangulated manifold is a
combinatorial manifold if and only if the link of each of its
vertices is a combinatorial sphere or combinatorial ball.
(Recall that the {\em link} of
a vertex $x$ in a complex $X$, denoted by ${\rm lk}_X(x)$, is
the subcomplex $\{\alpha \in X \, : \, x\not\in \alpha, \alpha
\sqcup\{x\} \in X\}$. Also, the {\em star} of $x$ in $X$, denoted
by ${\rm st}_X(x)$, is the cone $x \ast {\rm lk}_X(x)$.
The {\em antistar} of $x$ in $X$, denoted by ${\rm ast}_X(x)$,
is the subcomplex $\{\alpha\in X \, : \, x\not\in\alpha\}$.) This
leads us to introduce\,:

\begin{defn}\hspace{-1.8mm}{\bf .} \label{stellated-sphere}
{\rm For $0\leq k \leq d+1$, a $d$-dimensional simplicial complex
$X$ is said to be {\em $k$-stellated} if $X$ may be obtained from
$S^{\,d}_{d+2}$ by a finite sequence of bistellar moves, each of
index $< k$. By convention, $S^{\,d}_{d+2}$ is the only
$0$-stellated simplicial complex of dimension $d$.

Clearly, for $0\leq k \leq l \leq d+1$, $k$-stellated implies
$l$-stellated. All $k$-stellated simplicial complexes are
combinatorial spheres. We let $\Sigma_k(d)$ denote the class of
all $k$-stellated $d$-spheres. By Pachner's theorem (\cite{p}),
$\Sigma_{d+1}(d)$ consists of all combinatorial $d$-spheres.
 }
\end{defn}

By definition, $X\in \Sigma_k(d)$ if and only if there is a
sequence $X_0, X_1, \dots, X_n$ of $d$-dimensional simplicial
complexes such that $X_0 = S^{\,d}_{d + 2}$, $X_n = X$ and, for $0
\leq j < n$, $X_{j + 1}$ is obtained from $X_j$ by a single
bistellar move of index $\leq k-1$. The smallest such integer $n$
is said to be the {\em length} of $X\in \Sigma_k(d)$ and is
denoted by $l(X)$. For $X, Y \in \Sigma_k(d)$, we say that $Y$ is
shorter than $X$ if $l(Y) < l(X)$. Thus, $S^{\,d}_{d + 2}$ is the
unique shortest member of $\Sigma_k(d)$ (of length 0), and every
other member of $\Sigma_k(d)$ can be obtained from a shorter
member by a single bistellar move of index $< k$. Thus, induction
on the length is a natural method for proving results about the
class $\Sigma_k(d)$.

Let $X$, $Y$ be two pure simplicial complexes of dimension $d$. We
say that $X$ is obtained from $Y$ by the {\em shelling move}
$\alpha \leadsto \beta$ if $\alpha$ and $\beta \neq \emptyset$ are
disjoint faces of $X$ such that (i) $Y \subseteq X$, and $\alpha
\sqcup \beta$ is the only facet of $X$ which is not a facet of
$Y$, and (ii) the induced subcomplex of $Y$ on the vertex set of
$\alpha \sqcup \beta$ is $\overline{\alpha}\ast \partial \beta$.
If $\dim(\beta) = i$, we say that the shelling move $\alpha
\leadsto \beta$ is of index $i$. (Clearly, $\dim(\alpha) +
\dim(\beta) = d - 1$, so that $0 \leq i \leq d$).

We say that a $d$-dimensional simplicial complex $X$ {\em
shellable} if $X$ is obtained from the standard $d$-ball
$B^{\,d}_{d + 1}$ by a finite sequence of shelling moves. Clearly,
each shelling move increases the number of facets by one, so that
- when $X$ is shellable, the number of shelling moves needed to
obtain $X$ from $B^{\,d}_{d + 1}$ is one less than the number of
facets of $X$.

Let $X$ and $Y$ be $d$-dimensional pseudomanifolds. If $X$ is
obtained from $Y$ by the shelling move $\alpha \leadsto \beta$
then $X = Y \cup \overline{\alpha \sqcup \beta}$, $Y \cap
\overline{\alpha \sqcup \beta} = \overline{\alpha}\ast \partial
\beta$. (Since $X$ is a pseudomanifold, it follows that
$\overline{\alpha}\ast \partial \beta \subseteq \partial Y$.) If
the move is of index $<d$, then $\overline{\alpha} \ast \beta$ is
a combinatorial $(d-1)$-ball; if it is of index $d$ (so that
$\alpha = \emptyset$), $\overline{\alpha} \ast \partial \beta$
($=\partial \beta$) is a combinatorial $(d-1)$-sphere. Therefore,
if $Y$ is a combinatorial $d$-ball, then $X$ is also a
combinatorial $d$-ball in case the shelling move is of index $<d$,
and $X$ is a combinatorial $d$-sphere if the shelling move is of
index $d$. (Also note that $Y$ can't be a combinatorial sphere
since a $d$-dimensional pseudomanifold without boundary can't be
properly contained in a $d$-pseudomanifold with or without
boundary.) From these observations, it is immediate by an
induction on the number of facets that a shellable pseudomanifold
of dimension $d$ is either a combinatorial ball or a combinatorial
sphere. (This result appears to be due to Danaraj and Klee
\cite{dk}.) Also if $X$ is a shellable $d$-pseudomanifold, then
among the shelling moves used to obtain $X$ from $B^{\,d}_{d
+1}$, only the last move can be of index $d$; this happens if and
only if $X$ is a $d$-sphere. These considerations lead us to
introduce\,:

\begin{defn}\hspace{-1.8mm}{\bf .} \label{k-shelled-ball}
{\rm For $0\leq k \leq d$, a $d$-dimensional pseudomanifold is
said to be {\em $k$-shelled} if it may be obtained from the
standard $d$-ball $B^{\,d}_{d+1}$ by a finite sequence of shelling
moves, each of index $< k$. By convention, $B^{\,d}_{d+1}$ is the
only $0$-shelled pseudomanifold of dimension $d$. }
\end{defn}

Clearly, all $k$-shelled pseudomanifolds are combinatorial balls.
Also, for $0\leq k \leq l\leq d$, $k$-shelled implies $l$-shelled.
By $\widehat{\Sigma}_k(d)$, $0\leq k\leq d$, we denote the class
of all $k$-shelled $d$-balls. Thus $\widehat{\Sigma}_d(d)$
consists of all the shellable $d$-balls. Note that, while all
shellable balls are combinatorial balls, the converse is false.

Unlike the case of bistellar moves, the reverse of a shelling move
is not a shelling move. Nonetheless, the two notions are closely
related, as the following lemma shows.

\begin{lemma}\hspace{-1.8mm}{\bf .} \label{L1}
If a triangulated $(d+1)$-ball $X$ is obtained from a triangulated
$(d+1)$-ball $Y$ by a shelling move $\alpha \leadsto \beta$ of
index $i \leq d$ then the triangulated $d$-sphere $\partial X$ is
obtained from the triangulated $d$-sphere $\partial Y$ by the
bistellar move $\alpha \mapsto \beta$ of index $i$.
\end{lemma}

\noindent {\bf Proof.} Let $\sigma = \alpha \sqcup \beta$. Thus,
$\sigma$ is the only facet of $X$ which is not in $Y$. Since $Y
\subseteq X$ are $(d+1)$-dimensional pseudomanifolds, it follows
that (i) a boundary $d$-face of $Y$ is not a boundary $d$-face of
$X$ if and only if (it is a face of $Y$ and) it is contained in
$\sigma$, i.e., if and only if it is a facet of
$\overline{\alpha}\ast \partial \beta$, and (ii) a boundary $d$
face of $X$ is not a face of $Y$ if and only if it is a facet of
$\overline{\beta}\ast \partial \alpha$. Since $\partial X$ and
$\partial Y$ are pure simplicial complexes of dimension $d$, the
result follows. \hfill $\Box$

\medskip

As an immediate consequence of this lemma, we have\,:

\begin{cor}\hspace{-1.8mm}{\bf .} \label{C1}
If $B$ is a $k$-shelled $(d+1)$-ball then $\partial B$ is a
$k$-stellated $d$-sphere.
\end{cor}

For a simplicial complex $X$, say of dimension $d$, and a
non-negative integer $m \leq d$, the {\em $m$-skeleton}  of $X$,
denoted by ${\rm skel}_m(X)$, is the subcomplex of $X$ consisting
of all its faces of dimension $\leq m$. We recall\,:

\begin{defn}\hspace{-1.8mm}{\bf .} \label{k-stacked-ball}
{\rm For $0\leq k \leq d+1$, a triangulated $(d+1)$-dimensional
ball $B$ is said to be {\em $k$-stacked} if all the faces of $B$
of codimension (at least) $k+1$ lie in its boundary; i.e., if
${\rm skel}_{d - k}(B) = {\rm skel}_{d - k}(\partial B)$. A
triangulated $d$-sphere $S$ is said to be {\em $k$-stacked} if
there is a $k$-stacked $(d+ 1)$-ball $B$ such that $\partial B =
S$.
We let ${\cal S}_k(d)$ and  $\widehat{\cal S}_k(d)$ denote the
class of all $k$-stacked $d$-spheres and of all $k$-stacked
$d$-balls respectively. }
\end{defn}

Clearly, we have ${\cal S}_0(d) \subseteq
{\cal S}_1(d) \subseteq \cdots \subseteq {\cal S}_d(d)$ and
$\widehat{\cal S}_0(d) \subseteq \widehat{\cal S}_1(d) \subseteq
\cdots \subseteq \widehat{\cal S}_d(d)$. Trivially, the standard
$d$-ball is the only member of $\widehat{\cal S}_0(d)$, and hence
the standard $d$-sphere is the only member of ${\cal S}_0(d)$. Our
first proposition shows that ${\cal S}_d(d)$ consists of all the
triangulated $d$-spheres. Notice that, trivially, $\widehat{\cal
S}_d(d)$ consists of all triangulated $d$-balls.

\begin{prop}\hspace{-1.8mm}{\bf .} \label{P1}
Every triangulated $d$-sphere is $d$-stacked.
\end{prop}

\noindent {\bf Proof.} Let $S$ be a triangulated $d$-sphere. Fix a
vertex $x$ of $S$. Let $A_x$ be the antistar of $x$ in $S$.
Set $B_x =
\overline{\{x\}} \ast A_x$. It is shown in Lemma 9.1 of \cite{bd9}
that $B_x$ is a triangulated $(d+1)$-ball. Clearly, $B_x$ has the
same vertex set as $S = \partial B_x$. Therefore, $B_x$ is a
$d$-stacked $(d+1)$-ball and (hence) $S$ is a $d$-stacked
$d$-sphere. \hfill $\Box$

\begin{prop}\hspace{-1.8mm}{\bf .} \label{P2}
Let $B$ be a triangulated $(d+1)$-ball. Then $B$ is $k$-shelled if
and only if $B$ is shellable and $k$-stacked.
\end{prop}

\noindent {\bf Proof.} Suppose $B$ is $k$-shelled. Then, of
course, $B$ is shellable. We prove that $B$ is $k$-stacked by
induction on the number of facets of $B$. If $B$ has only one
facet then $B = B^{\,d + 1}_{d + 2}$, the standard ball, and the
result is trivial. Otherwise, $B$ is obtained from a $k$-shelled
ball $B^{\hspace{.2mm} \prime}$ (with one less facet) by a single
shelling move $\alpha \leadsto \beta$ of index $\leq k - 1$. By
induction hypothesis, ${\rm skel}_{d - k}(B^{\hspace{.2mm}
\prime}) = {\rm skel}_{d - k}(\partial B^{\hspace{.2mm} \prime})$,
and by Lemma \ref{L1}, $\partial B$ is obtained from $\partial
B^{\hspace{.2mm} \prime}$ by the bistellar move $\alpha \mapsto
\beta$ of index $\leq k-1$.

Let $\gamma$ be a face of $B$ of dimension $\leq d-k$. Since
$\dim(\alpha) \geq d-k+1$, $\gamma \not \supseteq \alpha$. If
$\gamma$ is a face of $B^{\hspace{.2mm} \prime}$ then (as
$B^{\hspace{.2mm} \prime}$ is $k$-stacked), $\gamma \in \partial
B^{\hspace{.2mm} \prime}$. Since $\gamma \not\supseteq \alpha$,
and $\partial B$ is obtained from $\partial B^{\hspace{.2mm}
\prime}$ by the bistellar move $\alpha \mapsto \beta$, it follows
that $\gamma \in \partial B$. If, on the other hand, $\gamma$ is
not a face of $B^{\hspace{.2mm} \prime}$ then $\beta \subseteq
\gamma \subseteq \alpha\sqcup\beta$ and hence we have $\gamma \in
\overline{\beta}\ast \partial \alpha \subseteq \partial B$. Thus
$\gamma \in \partial B$ in either case. So, $B$ is $k$-stacked.
This proves the ``only if\," part.

The ``if part\," is also proved by induction on the number of facets
of $B$. Suppose $B$ is a $k$-stacked shellable $(d+1)$-ball. If $B
= B^{\,d+1}_{d+2}$, then $B$ is vacuously $k$-shelled. Else, $B$
is obtained from a shellable $(d+1)$-ball $B^{\hspace{.2mm}
\prime}$ (with one less facet) by a single shelling move $\alpha
\leadsto \beta$. By Lemma \ref{L1}, $\partial B$ is obtained from
$\partial B^{\hspace{.2mm} \prime}$ by the bistellar move $\alpha
\mapsto \beta$. Hence $\alpha \not\in \partial B$ but $\alpha \in
B$. Since $B$ is $k$-stacked, it follows that $\dim(\alpha) \geq
d-k+1$, and hence $\dim(\beta) \leq k-1$. Thus, the shelling move
$\alpha \leadsto \beta$ is of index $\leq k-1$. Let $\gamma \in
B^{\hspace{.2mm} \prime}$, $\dim(\gamma) \leq d-k$. Since
$B^{\hspace{.2mm} \prime}\subseteq B$, it follows that $\gamma \in
B$. Since $\dim(\gamma) \leq d-k$ and $B$ is $k$-stacked, it
follows that $\gamma \in \partial B$. As $\beta \not\in
B^{\hspace{.2mm} \prime}$ and $\gamma \in B^{\hspace{.2mm}
\prime}$, we also have $\gamma \not\supseteq \beta$. Thus $\gamma
\not\supseteq \beta$, $\gamma \in \partial B$ and $\partial B$ is
obtained from $\partial B^{\hspace{.2mm} \prime}$ by the bistellar
move $\alpha \mapsto \beta$. Hence $\gamma \in \partial
B^{\hspace{.2mm} \prime}$. This shows that $B^{\hspace{.2mm}
\prime}$ is $k$-stacked. As $B^{\hspace{.2mm} \prime}$ is
$k$-stacked and shellable, the induction hypothesis implies that
$B^{\hspace{.2mm} \prime}$ is $k$-shelled. Since $B$ is obtained
from $B^{\hspace{.2mm} \prime}$ by a shelling move of index $\leq
k-1$, it follows that $B$ is also $k$-shelled. This completes the
induction. \hfill $\Box$

\bigskip

Thus we have $\widehat{\Sigma}_k(d) \subseteq \widehat{\cal
S}_k(d)$. Our next result gives a one-sided relationship between
$k$-stacked spheres and $k$-stacked balls on one hand, and between
$k$-stellated spheres and $k$-shelled balls on the other hand.

\begin{prop}\hspace{-1.8mm}{\bf .} \label{P3}
Let $B$ be a triangulated ball. \vspace{-2mm}
\begin{enumerate}
\item[$(a)$] If $B$ is $k$-stacked then there is a $k$-stacked
sphere $S$ such that $B$ is the antistar of a vertex in $S$.
\vspace{-2mm}
 \item[$(b)$] If $B$ is $k$-shelled then there is a $k$-stellated
sphere $S$ such that $B$ is the antistar of a vertex in $S$.
\end{enumerate}
\end{prop}

\noindent {\bf Proof.} Let $x$ be a new vertex (not in $B$), and
set $S := B \cup (x\ast \partial B)$. (Notice that, since $S$ is
to be a $d$-pseudomanifold without boundary and $B$ is a
$d$-pseudomanifold with boundary, this is the only choice of $S$
so that $B$ is the antistar of a vertex $x$ in $S$.) Clearly, $S =
\partial B_0$, where $B_0 = x \ast B$. Therefore, to prove the
result, it is enough to show that if $B$ is $k$-stacked
(respectively $k$-shelled) then so is $B_0$. But, this is trivial.
\hfill $\Box$

\bigskip

Next we present a characterization of $k$-stellated spheres of
dimension $\geq 2k-1$.

\begin{prop}\hspace{-1.8mm}{\bf .} \label{P4}
A triangulated sphere of dimension $\geq 2k -1$ is $k$-stellated
if and only if it is the boundary of a $k$-shelled ball. In
consequence, all $k$-stellated spheres of dimension $\geq 2k-1$
are $k$-stacked.
\end{prop}

\noindent {\bf Proof.} The ``if\," part is Corollary \ref{C1}
(which holds in all dimensions). We prove the ``only if\," part by
induction on the length $l(S)$ of a $k$-stellated sphere $S$ of
dimension $d\geq 2k-1$. If $l(S) = 0$ then $S = S^{\,d}_{d+2}$ is
the boundary of $B^{\,d+1}_{d+2}$. So, let $l(S) > 0$. Then $S$ is
obtained from a shorter member $S^{\hspace{.2mm} \prime}$ of
$\Sigma_k(d)$ by a single bistellar move $\alpha \mapsto \beta$ of
index $\leq k-1$. By induction hypothesis, there is a $k$-shelled
$(d+1)$-ball $B^{\hspace{.2mm}\prime}$ such that $\partial
B^{\hspace{.2mm}\prime} = S^{\hspace{.2mm}\prime}$. The induced
subcomplex of $S^{\hspace{.2mm} \prime}$ on the vertex set
$\alpha\sqcup \beta$ is $\overline{\alpha} \ast \partial \beta
\subseteq S^{\hspace{.2mm} \prime} \subseteq B^{\hspace{.2mm}
\prime}$. Since $\dim(\beta) \leq k-1 \leq d-k$, $\beta \not\in
S^{\hspace{.2mm} \prime} = \partial B^{\hspace{.2mm} \prime}$ and
(by Proposition \ref{P2}) $B^{\hspace{.2mm} \prime}$ is
$k$-stacked, it follows that $\beta \not\in B^{\hspace{.2mm}
\prime}$. Thus, the induced subcomplex of $B^{\hspace{.2mm}
\prime}$ on $\alpha \sqcup \beta$ is also $\overline{\alpha} \ast
\partial \beta$. So, $B^{\hspace{.2mm} \prime}$ admits the
shelling move $\alpha \leadsto \beta$ of index $\leq k-1$. Let $B$
be the $(d+1)$-ball obtained from $B^{\hspace{.2mm} \prime}$ by
this move. Since $B^{\hspace{.2mm} \prime}$ is $k$-shelled, so is
$B$. By Lemma \ref{L1}, $\partial B$ is obtained from
$S^{\hspace{.2mm} \prime} = \partial B^{\hspace{.2mm} \prime}$ by
the bistellar move $\alpha \mapsto \beta$. That is, $\partial B =
S$. This completes the induction. The second statement is now
immediate from the first statement and Proposition \ref{P2}.
\hfill $\Box$

\begin{prop}\hspace{-1.8mm}{\bf .} \label{P5}
Let $S$ be a $k$-stacked sphere of dimension $d\geq 2k+1$. Then there
is a unique $k$-stacked ball $\overline{S}$ such that $\partial
\overline{S} = S$.
\end{prop}

\noindent {\bf Proof.} Suppose $B_1$ and $B_2$ are two $k$-stacked
balls with $\partial B_1 = S = \partial B_2$. Put $S_i = \partial(x
\ast B_i)$, $i = 1, 2$, where $x$ is a new vertex. Then $S_1$ and
$S_2$ are two $(d+1)$-spheres with identical $(d-k)$-skeletons.
Since $d \geq 2k+1$, Theorems 1 and 2 in \cite{d} imply $S_1
=S_2$ and hence $B_1 = B_2$. \hfill $\Box$

\bigskip

For $k$-stellated spheres, we can improve the bound in
Proposition \ref{P5}. For $d \geq 2k$, such spheres are the
boundaries of uniquely determined $k$-shelled
balls. To describe this result, we introduce\,:

\medskip

\noindent {\bf Notation\,:} For a set $\alpha$ and a non-negative
integer $m$, ${\alpha \choose \leq \, m\,}$ will denote the
collection of all subsets of $\alpha$ of size $\leq m$. Also,
${\alpha \choose m\,}$ will denote the
collection of all subsets of $\alpha$ of size $= m$

\begin{prop}\hspace{-1.8mm}{\bf .} \label{P6}
Let $S$ be a $k$-stellated sphere of dimension $d \geq 2k$, say
with vertex set $V$. Then there is a unique $k$-stacked $(d+
1)$-ball $\overline{S}$ whose boundary is $S$. $($By Propositions
$\ref{P2}$ and $\ref{P4}$, $\overline{S}$ is actually $k$-shelled.$)$
It is given by the
formula \vspace{-2mm}
\begin{eqnarray} \label{eq1}
\overline{S} = \left\{\alpha \subseteq V \, : \, {\alpha \choose
\leq \, k+1\,} \subseteq S\right\}.
\end{eqnarray}
\end{prop}

\noindent {\bf Proof.} The existence of a  $k$-stacked ball $B$
with boundary $S$ is guaranteed by Proposition \ref{P4}. We prove
that $B = \overline{S}$ by induction on $l(S)$.

If $l(S)= 0$, then $S = S^{\,d}_{d+ 2}$ and trivially $B^{\,d+
1}_{d +2}$ is the unique $k$-stacked $(d+1)$-ball with boundary
$S^{\,d}_{d+2}$; it is indeed given by (\ref{eq1}). So, let $l(S)
> 0$. Then $S$ is obtained from a shorter member $S^{\hspace{.2mm}
\prime}$ of $\Sigma_k(d)$ by a bistellar move $\alpha \mapsto
\beta$ of index $\leq k - 1$. By induction hypothesis,
$\overline{S^{\hspace{.2mm}\prime}}$ (given by (\ref{eq1}) with
$S^{\hspace{.2mm}\prime}$ in place of $S$, and the vertex set
$V^{\prime}$ of $S^{\hspace{.2mm}\prime}$ in place of $V$) is the
unique $k$-stacked ball with boundary $S^{\hspace{.2mm}\prime}$.

Let $B$ be a $k$-stacked ball with boundary $S$.  We need to show
that $B = \overline{S}$. First we claim that $\alpha
\sqcup \beta \in B$. To prove this, fix a vertex $a \in \alpha$
and look at the boundary $d$-face $\alpha \sqcup \beta \setminus
\{a\}$ of $B$. Let $\sigma$ be the unique facet of $B$ containing
this $d$-face. Suppose, if possible, that $\sigma \neq \alpha
\sqcup \beta$. Let $b$ be the unique vertex in $\sigma \setminus
(\alpha \sqcup \beta)$. Then $\beta \cup \{b\} \in B$ and
$\dim(\beta \cup\{b\}) \leq k \leq d-k$, so that $\beta\cup \{b\}
\in S$. This is a contradiction since ${\rm lk}_S(\beta) =
\partial \alpha$ and $b \not \in \alpha$. This proves the claim\,:
$\alpha \sqcup \beta \in B$.

Put $B^{\hspace{.2mm}\prime} = (B \setminus \overline{\alpha
\sqcup\beta}) \cup (\overline{\alpha} \ast \partial \beta)$. Since
$\alpha\sqcup\beta$ is a facet of the triangulated $(d+1)$-ball
$B$, and $\overline{\alpha\sqcup \beta} \cap \partial B =
(\partial \alpha) \ast \overline{\beta}  $ is a $d$-ball, it
follows that $B^{\hspace{.2mm} \prime}$ is a triangulated
$(d+1)$-ball. Clearly, $\partial B^{\hspace{.2mm} \prime} =
S^{\hspace{.2mm} \prime}$. Let $\gamma$ be a face of
$B^{\hspace{.2mm}\prime}$ with $\dim(\gamma) \leq d-k$. Since
$B^{\hspace{.2mm}\prime} \subseteq B$, $\gamma\in B$. As $B$ is
$k$-stacked, it follows that $\gamma\in S = \partial B$. Since
$\beta \not\in B^{\hspace{.2mm}\prime}$, we have $\gamma
\not\supseteq \beta$. Therefore, $\gamma \in S^{\hspace{.2mm}
\prime} = \partial B^{\hspace{.2mm}\prime}$. Thus
$B^{\hspace{.2mm} \prime}$ is a $k$-stacked $(d+1)$-ball with
$\partial B^{\hspace{.2mm}\prime} = S^{\hspace{.2mm}\prime}$.
Therefore, by the induction hypothesis, we have
$B^{\hspace{.2mm}\prime} = \overline{S^{\hspace{.2mm}\prime}}$.
Hence $B = B^{\hspace{.2mm}\prime} \cup \overline{\alpha \sqcup
\beta} = \overline{S^{\hspace{.2mm}\prime}} \cup \overline{\alpha
\sqcup \beta}$.

So, to complete the proof, we need to show that $B =
\overline{S^{\hspace{.2mm}\prime}} \cup \overline{\alpha \sqcup
\beta}$ equals $\overline{S}$. Since $d \geq 2k$, it follows that
${\rm skel}_k(B) \subseteq {\rm skel}_{d - k}(B) \subseteq S$.
Therefore, $\sigma \in B$ $\Rightarrow$ ${\sigma \choose \leq \,
k+1} \subseteq S$ $\Rightarrow$ $\sigma \in \overline{S}$. Thus,
$B \subseteq \overline{S}$. To prove the reverse inclusion
$\overline{S} \subseteq B$, let $\sigma \in\overline{S}$. If
$\beta \not\subseteq \sigma$ then ${\sigma \choose \leq \,k+1}
\subseteq S$ $\Rightarrow$ ${\sigma \choose \leq \,k+1} \subseteq
S^{\hspace{.2mm}\prime}$ $\Rightarrow$ $\sigma \in
\overline{S^{\hspace{.2mm}\prime}} \subseteq B$ $\Rightarrow$
$\sigma \in B$. So, suppose $\beta \subseteq \sigma$. If $\sigma
\not\subseteq \alpha\sqcup \beta$ then take a vertex $x \in
\sigma\setminus (\alpha\sqcup\beta)$. Then $\beta\cup\{x\} \in
{\sigma \choose \leq \, k+1} \subseteq S \Rightarrow
\beta\cup\{x\} \in S$ with $x \not\in\alpha$. This is a
contradiction since ${\rm lk}_S(\beta) = \partial \alpha$. Thus,
in this case, $\sigma \subseteq \alpha\sqcup \beta \in B$. Hence
$\sigma \in B$ in either case. Thus $\overline{S} \subseteq B$ and
hence $B = \overline{S}$. \hfill $\Box$

\begin{cor}\hspace{-1.8mm}{\bf .} \label{C2}
For $k \leq e \leq d-k-1$, a $k$-stellated $d$-sphere does not
have any standard $e$-sphere as an induced subcomplex. In
consequence, such a $d$-sphere does not admit any bistellar move
of index $i$ for $k+1 \leq i \leq d-k$.
\end{cor}

\noindent {\bf Proof.} Notice that a triangulated sphere $S$
admits a bistellar move $\alpha \mapsto \beta$ of index $i$ if and
only if it has $\overline{\alpha}\ast \partial\beta$ as an induced
subcomplex. In this case, it has the standard $(i-1)$-sphere
$\partial \beta$ as the induced subcomplex on $\beta$. So, the
second statement is immediate from the first. The first statement
is vacuously true unless $d \geq 2k+1$. So, to prove it, we may
assume $d \geq 2k+1$. If $S$ contains a standard $e$-sphere as an
induced subcomplex on the vertex-set $\gamma$ (so, $\#(\gamma) =
e+2$), then all the proper subsets of $\gamma$ are faces of $S$.
In particular, if $e \geq k$, all the subsets of $\gamma$ of size
$\leq k+1$ are faces of $S$. Hence $\gamma \in \overline{S}$. If,
also, $e \leq d-k- 1$, then $\gamma \in {\rm skel}_{d-
k}(\overline{S}) = {\rm skel}_{d-k}(S)$ (by Proposition \ref{P6})
and hence $\gamma \in S$. Then the induced subcomplex of $S$ on
the vertex set $\gamma$ is the ball $\overline{\gamma}$, a
contradiction. \hfill $\Box$

\bigskip

If $S$ is a $k$-stellated $d$-sphere, other than the standard
sphere, then $S$ is obtained from a shorter $k$-stellated
$d$-sphere by a bistellar move of index $\leq k-1$. Hence such a
sphere admits the reverse move, which is a bistellar move of index
$\geq d-k+1$. In consequence, such a sphere always has an induced
subcomplex isomorphic to a standard sphere of some dimension $\geq
d-k$. In this sense, Corollary \ref{C2} is best possible. Indeed,
it is easy to prove by induction on the length that if $d \geq
2k-2$ and $S$ is a $k$-stellated $d$-sphere which is not
$(k-1)$-stellated, then $S$ has an $S^{\,d-k}_{d-k+2}$ as an
induced subcomplex.

In the following proof (and also later) we use the notation $V(X)$
for the vertex set of a simplicial complex $X$.

\begin{prop}\hspace{-1.8mm}{\bf .} \label{P7}
For a pseudomanifold $X$, the following are equivalent\,:
\vspace{-1mm}
\begin{enumerate}
\item[$(i)$] $X$ is a $1$-shelled ball,
\vspace{-1mm}
\item[$(ii)$] $X$ is a $1$-stacked ball,
\vspace{-1mm}
\item[$(iii)$] $\Lambda(X)$ is a tree.
\end{enumerate}
\end{prop}

\noindent {\bf Proof.} Let $X$ be of dimension $d+1\geq 1$.

$(i) \Rightarrow (ii)$\,: Follows from Proposition \ref{P2}.

$(ii) \Rightarrow (iii)$\,: The result is trivial for dimension 1.
So, assume that $d+ 1 \geq 2$. If $X$ has only one facet then the
result is trivial. So, assume that $X$ is a 1-stacked ball with at
least two facets. Since $X$ is a ball, $\Lambda(X)$ is connected.
To prove that $\Lambda(X)$ is a tree, it suffices to show that
each edge of $\Lambda(X)$ is a cut edge (i.e., deletion of any
edge from $\Lambda(X)$ disconnects the graph). Let $e_0 =
\sigma_1\sigma_2$ be an edge of $\Lambda(X)$. Then $\gamma :=
\sigma_1\cap \sigma_2$ is an interior $d$-face of $X$; i.e.,
$\gamma \not\in S := \partial X$. Since ${\rm skel}_{d-1}(X) =
{\rm skel}_{d-1}(S)$, $\partial\gamma\subseteq S$. Thus, $\partial
\gamma$ is an induced $S^{\,d-1}_{d+1}$ in the $d$-sphere $S$. By
Lemma 3.3 of \cite{bd9}, $S$ is obtained from a $d$-dimensional
weak pseudomanifold $\widetilde{S}$ (without boundary) by an
elementary handle addition. Since $S$ is simply connected, the
Seifert-Van Kampen theorem implies that $\widetilde{S}$ is
disconnected and hence has exactly two components (again by Lemma
3.3 of \cite{bd9}), say $S_1$ and $S_2$. Then $S = S_1 \# S_2$
(connected sum) and $V(S_1) \cap V(S_2) = \gamma$. For $1\leq i
\leq 2$, let $U_i$ be the set of facets of $X$ contained in
$V(S_i)$. Since $V(S_1) \cap V(S_2) = \gamma$, it follows that
$U_1 \cap U_2 = \emptyset$.

If the dimension $d+ 1 = 2$ then $\gamma$ is an edge and it
clearly divides the 2-disc $X$ into two parts and the triangles
(facets) in one part are in $U_1$ and the triangles in the other
part are in $U_2$. Now, assume that $d+1 \geq 3$. Let $uv$ be an
edge of $X$. Since $d+1 \geq 3$, $uv \in S$ and hence (since $S =
S_1\# S_2$) $uv \in S_1$ or $uv \in S_2$. Therefore, $u, v\in
V(S_1)$ or $u, v\in V(S_2)$. This implies that for any facet
$\sigma$ in $X$, either all the vertices of $\sigma$ are in
$V(S_1)$ or all the vertices of $\sigma$ are in $V(S_2)$. Thus,
any facet in $X$ is in $U_1$ or in $U_2$. Thus (for any dimension
$d+1 \geq 2$), $U_1 \sqcup U_2$ is a partition of the vertex-set
of the dual graph $\Lambda(X)$. Any facet $\sigma$ of $X$
containing a $d$-face $\alpha \neq \gamma$ of $S_1$ is in $U_1$.
So, $U_1 \neq \emptyset$. Similarly, $U_2 \neq \emptyset$.

Now, let $e = \alpha_1 \alpha_2$ be an edge of $\Lambda(X)$ with
$\alpha_i\in U_i$, $i = 1, 2$. Then $\alpha := \alpha_1\cap
\alpha_2 \subseteq V(S_i)$ for $i=1, 2$. Hence $\alpha \subseteq
V(S_1) \cap V(S_2) = \gamma$ and therefore $\alpha = \gamma$. So,
$e = e_0$. Thus, $e_0$ is the unique edge of $\Lambda(X)$ with one
end in $U_1$ and other end in $U_2$. So, $e_0$ is a cut edge of
$\Lambda(X)$. Since $e_0$ was an arbitrary edge of $\Lambda(X)$,
this proves that $\Lambda(X)$ is a tree.

$(iii) \Rightarrow (i)$\,: Suppose $\Lambda(X)$ is a tree. We
prove that $X$ is 1-shelled by induction on the number of facets
of $X$ (i.e., the number of vertices of $\Lambda(X)$). This is
trivial if $X$ has only one facet, i.e., $X = B^{\,d+1}_{d+2}$.
So, assume $\Lambda(X)$ is a tree with at least two vertices. Then
$\Lambda(X)$ has a vertex $\sigma$ of degree 1 (an end vertex).
Let $\sigma^{\hspace{.2mm}\prime}$ be the unique neighbour of
$\sigma$ in $\Lambda(X)$, and put $\gamma = \sigma \cap
\sigma^{\hspace{.2mm}\prime}$. Let $X^{\hspace{.2mm}\prime} = (X
\setminus \overline{\sigma}) \cup \overline{\gamma}$. Then
$X^{\hspace{.2mm} \prime}$ is a pseudomanifold and
$\Lambda(X^{\hspace{.2mm}\prime})$ is the tree obtained from the
tree $\Lambda(X)$ by deleting the end vertex $\sigma$ and the edge
$\sigma\sigma^{\hspace{.2mm}\prime}$. Therefore, by induction
hypothesis, $X^{\hspace{.2mm}\prime}$ is an 1-shelled ball. If $u$
is the vertex of $X$ in $\sigma\setminus \gamma$, then $X$ is
obtained from $X^{\hspace{.2mm}\prime}$ by the shelling move
$\gamma \leadsto \{u\}$ of index 0. Therefore, $X$ is also an
1-shelled ball. \hfill $\Box$

\bigskip

Thus a triangulated ball is 1-stacked if and only if it is
1-shelled. So, $\widehat{\Sigma}_1(d) =\widehat{\cal S}_1(d)$.
Now, Propositions \ref{P4} and \ref{P7} imply\,:

\begin{cor}\hspace{-1.8mm}{\bf .} \label{C3}
A triangulated sphere is $1$-stellated if and only if it is
$1$-stacked.
\end{cor}

Next we introduce\,:

\begin{defn}\hspace{-1.8mm}{\bf .} \label{D4}
{\rm For $0 \leq k \leq d$, ${\cal W}_k(d)$ consists of the
connected simplicial complexes of dimension $d$ all whose
vertex-links are $k$-stellated $(d-1)$-spheres, and ${\cal
K}_k(d)$ consists of the connected simplicial complexes of
dimension $d$ all whose vertex-links are $k$-stacked
$(d-1)$-spheres. }
\end{defn}

Thus, members of ${\cal W}_k(d)$ are combinatorial manifolds; the
members of ${\cal K}_k(d)$ are triangulated manifolds. In
consequence of Corollary \ref{C3}, we have\,:

\begin{cor}\hspace{-1.8mm}{\bf .} \label{C4}
${\cal W}_1(d) ={\cal K}_1(d)$.
\end{cor}

In consequence of Proposition \ref{P4}, we have\,:

\begin{cor}\hspace{-1.8mm}{\bf .} \label{C5}
${\cal W}_k(d) \subseteq {\cal K}_k(d)$ for $d \geq 2k$.
\end{cor}

\begin{prop}\hspace{-1.8mm}{\bf .} \label{P8}
$(a)$ All $k$-stellated $d$-spheres belong to the class ${\cal
W}_k(d)$. $(b)$ All $k$-stacked $d$-spheres belong to the class
${\cal K}_k(d)$.
\end{prop}

\noindent {\bf Proof.} Let $S$ be a $k$-stellated $d$-sphere. We
need to show that all the vertex-links of $S$ are $k$-stellated.
Again, the proof is by induction on the length $l(S)$ of $S$. If
$l(S) = 0$ then $S = S^{\,d}_{d +2}$, and all its vertex links are
$S^{\,d-1}_{d +1}$, so we are done. Therefore, let $l(S) >0$. Then
$S$ is obtained from a shorter $k$-stellated $d$-sphere
$S^{\hspace{.2mm}\prime}$ by a bistellar move $\alpha \mapsto
\beta$ of index $\leq k-1$. Let $x$ be a vertex of $S$. If $x
\not\in \alpha \sqcup\beta$ then ${\rm lk}_S(x) = {\rm
lk}_{S^{\hspace{.2mm}\prime}}(x)$ is $k$-stellated by induction
hypothesis. If $x \in \alpha$ then ${\rm lk} _S(x)$ is obtained
from the $k$-stellated sphere ${\rm lk}_{S^{\hspace{.2mm}
\prime}}(x)$ by the bistellar move $\alpha\setminus\{x\} \mapsto
\beta$ of index $\leq k-1$. If $x \in \beta$ and $\beta\neq \{x\}$
then ${\rm lk} _S(x)$ is obtained from the $k$-stellated sphere
${\rm lk}_{S^{\hspace{.2mm}\prime}}(x)$ by the bistellar move
$\alpha \mapsto \beta\setminus\{x\}$ of index $\leq k-2$. If
$\beta = \{x\}$ then ${\rm lk}_S(x)$ is the standard sphere
$\partial \alpha$. Thus, in all cases, ${\rm lk}_S(x)$ is
$k$-stellated. This proves part $(a)$.

Let $S$ be a $k$-stacked $d$-sphere. Let $B$ be a $k$-stacked
$(d+1)$-ball such that $\partial B = S$. If $x$ is a vertex of $S$
then $x$ is a vertex of $B$ and $B^{\hspace{.2mm}\prime} = {\rm
lk}_B(x)$ is a $d$-ball with $\partial B^{\hspace{.2mm}\prime} =
{\rm lk}_S(x)$. Therefore, it suffices to show that
$B^{\hspace{.2mm}\prime}$ is also $k$-stacked. Indeed, if $\gamma$
is a face of $B^{\hspace{.2mm}\prime}$ of codimension $\geq k+1$
then $\gamma \cup\{x\}$ is a face of $B$ of codimension $\geq
k+1$, and hence $\gamma \cup \{x\} \in \partial B = S$, so that
$\gamma \in {\rm lk}_S(x) = \partial B^{\hspace{.2mm}\prime}$.
\hfill $\Box$

\begin{prop}\hspace{-1.8mm}{\bf .} \label{P9}
Let $d \geq 2k+2$ and $M \in {\cal W}_k(d)$. Let $V(M)$ be the
vertex set of $M$. Then \vspace{-3mm}
\begin{eqnarray} \label{eq2}
\overline{M} := \left\{\alpha \subseteq V(M) \, : \, {\alpha
\choose \leq \, k+2\,} \subseteq M\right\}
\end{eqnarray}
is the unique combinatorial $\,(d+ 1)$-manifold such that $\,\partial
\overline{M} = M$ and $\,{\rm skel}_{d- k}(\overline{M}) = {\rm
skel}_{d- k}(M)$.
\end{prop}

\noindent {\bf Proof.} Fix $x \in V(\overline{M}) = V(M)$.

\smallskip

\noindent {\bf Claim\,:} ${\rm lk}_{\overline{M}}(x) =
\overline{{\rm lk}_M(x})$, where the right hand side is as defined
in Proposition \ref{P6}.

\smallskip

From the definition, we see that $\alpha \in {\rm
lk}_{\overline{M}}(x) \Rightarrow \alpha\sqcup\{x\}\in
\overline{M} \Rightarrow {\alpha\sqcup \{x\} \choose \leq k+2}
\subseteq M \Rightarrow {\alpha \choose \leq k+1} \subseteq {\rm
lk}_M(x) \Rightarrow \alpha \in \overline{{\rm lk}_M(x})$. Thus,
we have ${\rm lk}_{\overline{M}}(x)\subseteq \overline{{\rm
lk}_M(x})$.

Conversely, let $\alpha \in \overline{{\rm lk}_M(x})$. Then
${\alpha \choose \leq k+1} \subseteq {\rm lk}_M(x)$, so that each
$\gamma \subseteq \alpha \sqcup\{x\}$ such that $x \in \gamma$ and
$\#(\gamma) \leq k+ 2$ is in $M$. Therefore, to prove that $\alpha
\in {\rm lk}_{\overline{M}}(x)$, it suffices to show that each
$\gamma \subseteq \alpha$ with $\#(\gamma) \leq k+2$ is in $M$.
Since $\alpha \in \overline{{\rm lk}_M(x})$, such a set $\gamma$
is in $\overline{{\rm lk}_M(x})$, and hence $\gamma \in {\rm
skel}_{k+1}(\overline{{\rm lk}_{M}(x})) \subseteq {\rm skel}_{d-k-
1}(\overline{{\rm lk}_{M}(x})) = {\rm skel}_{d-k- 1}({\rm
lk}_{M}(x)) \subseteq {\rm lk}_M(x) \subseteq M$. (Here the first
inclusion holds since $k+1 \leq d-k-1$.) This proves that $\alpha
\in \overline{{\rm lk}_M(x}) \Rightarrow \alpha \in {\rm
lk}_{\overline{M}}(x)$, so that $\overline{{\rm lk}_M(x})
\subseteq {\rm lk}_{\overline{M}}(x)$. This proves the claim.

\smallskip

In view of Proposition \ref{P6}, the claim implies that
$\overline{M}$ is a combinatorial $(d+1)$-manifold with boundary,
and ${\rm lk}_{\partial \overline{M}}(x) = \partial({\rm
lk}_{\overline{M}}(x)) = \partial(\overline{{\rm lk}_M(x})) = {\rm
lk}_M(x)$ for every vertex $x$. Therefore, $\partial \overline{M}
= M$, and we have\,:
\begin{eqnarray*}
{\rm lk}_{{\rm skel}_{d-k}(\overline{M})}(x) & = & {\rm skel}_{d-
k-1}({\rm lk}_{\overline{M}}(x)) = {\rm skel}_{d-k-
1}(\overline{{\rm lk}_{M}(x})) = {\rm skel}_{d-k-1}({\rm
lk}_{M}(x)) \\ &=&  {\rm lk}_{{\rm skel}_{d-k}(M)}(x)
\end{eqnarray*}
for every vertex $x$. Thus, ${\rm skel}_{d-k}(\overline{M}) =
{\rm skel}_{d-k}(M)$.

Now, if $N$ is any $(d+1)$-manifold with $\partial N = M$ and
${\rm skel}_{d-k}(N) = {\rm skel}_{d-k}(M)$, then for any vertex
$x$, we have\,:
\begin{eqnarray*}
\partial({\rm lk}_N(x)) & = & {\rm lk}_{\partial N}(x) = {\rm
lk}_M(x), ~ \mbox{ and} \\
{\rm skel}_{d-k-1}({\rm lk}_N(x)) & = & {\rm lk}_{{\rm
skel}_{d-k}(N)}(x) = {\rm lk}_{{\rm skel}_{d-k}(M)}(x)  = {\rm
skel}_{d-k-1}({\rm lk}_{M}(x)).
\end{eqnarray*}
Therefore, the uniqueness assertion in Proposition \ref{P6}
implies that ${\rm lk}_N(x) = \overline{{\rm lk}_M(x}) = {\rm
lk}_{\overline{M}}(x)$ for every vertex $x$ and hence $N =
\overline{M}$. This completes the proof. \hfill $\Box$

\begin{remark}\hspace{-1.8mm}{\bf .} \label{R1}
{\rm If $M$ is a $k$-stellated sphere of dimension $d \geq 2k+2$
then $M \in {\cal W}_k(d)$ by Proposition \ref{P8}. In this case,
the uniqueness statements in Propositions \ref{P6} and \ref{P9}
show that the two definitions of $\overline{M}$ (given in (1) and
(2)) agree. Also, if we define $\overline{\cal W}_k(d+1)$ to be
the class of all $(d+1)$-dimensional simplicial complexes all whose
vertex links are $k$-shelled $d$-balls, then by Propositions
\ref{P6} and \ref{P9}, for $d\geq 2k+2$, $M \mapsto \overline{M}$
is a bijection from ${\cal W}_k(d)$ onto $\overline{\cal W}_k(d+1)$.
The boundary map provides its inverse.}
\end{remark}


\section{Polytopal spheres and balls\,: a diversion}

For a subset $A$ of an Euclidean space, we write ${\rm conv}(A)$
(respectively ${\rm aff}(A)$) for the convex (respectively affine)
hull of $A$. For a convex set $C$, the topological interior
(respectively the topological boundary) of $C$ in ${\rm aff}(C)$
is called the {\em relative interior} (respectively the {\em
boundary}) of $C$ and we denote it by $\intC$ (respectively
$\bdC$).

Recall that a (convex) {\em polytope} in the Euclidean space
$\RR^{\hspace{.1mm}n}$ is the convex hull of a finite set of
points. Equivalently, a polytope in $\RR^{\hspace{.1mm}n}$ is a
compact subset of $\RR^{\hspace{.1mm}n}$ which may be obtained as
the intersection of finitely many closed half-spaces of
$\RR^{\hspace{.1mm}n}$. As general references on polytopes, cf
\cite{g, z}. The {\em dimension} of a polytope $P$ is defined to
be the dimension of the affine space ${\rm aff}(P)$. A (geometric)
{\em simplex} is a polytope which is the convex hull of a set of
affinely independent points. A {\em face} of a polytope $P$ in
$\RR^{\hspace{.1mm}n}$ is either $P$ itself or is the intersection
of $P$ with a hyperplane $H$ of $\RR^{\hspace{.1mm}n}$ such that
$P$ is contained in one of the two closed half spaces determined
by $H$. The zero-dimensional faces of a polytope are called its
{\em vertices}, and the $d$-dimensional faces (i.e., maximal
proper faces) of a $(d+1)$-dimensional polytope are called its
{\em facets}. Notice that any polytope is the convex hull of its
vertex set. It is also the disjoint union of the relative interiors
of its faces.

Recall that a {\em geometric simplicial complex} $X$ is a
collection of geometric simplices such that the intersection of
any two members of $X$ is again a member of $X$ and any face of a
member of $X$ is again a member of $X$. If $X$ is a geometric
simplicial complex with vertex set $V(X)$, then $X_{\rm abs} :=
\{A \subseteq V(X) \, : \, {\rm conv}(A) \in X\}$ is an abstract
simplicial complex and is called the {\em abstract scheme} of $X$.
We sometimes identify $X$ with $X_{\rm abs}$.

A polytope is {\em simplicial} if all its proper faces are
simplices. If $P$ is a simplicial polytope then all the proper
faces of $P$ form a geometric simplicial complex ${\rm Bd}(P)$.
The abstract scheme of ${\rm Bd}(P)$ is called the {\em boundary
complex} of $P$ and is denoted by $\partial P$. Thus, $\partial P
= \{A \subseteq V(P) \, : \, {\rm conv}(A)$ is a proper face of
$P\}$. Clearly, the union of all the proper faces of $P$ is the
topological boundary $\bdP$ of $P$. Thus, the boundary complex
$\partial P$ of $P$ triangulates the topological sphere $\bdP$.
We identify ${\rm Bd}(P)$ with $\partial P$.

\begin{defn}\hspace{-1.8mm}{\bf .} \label{polytopal_sphere}
{\rm  A triangulated sphere is said to be a {\em polytopal sphere}
if it is isomorphic to the boundary complex of a simplicial
polytope. }
\end{defn}

\begin{defn}\hspace{-1.8mm}{\bf .} \label{simplicial_subdivision}
{\rm  A {\em simplicial subdivision} $P^{\hspace{0.2mm}\prime}$ of
a simplicial polytope $P$ is a geometric simplicial complex such
that $V(P^{\hspace{0.2mm}\prime}) = V(P)$ and $P$ is the
union of all the simplices in $P^{\hspace{0.2mm}\prime}$. Let
$\overline{P}$ be the abstract scheme of $P^{\hspace{0.2mm}
\prime}$.  Then $\overline{P}$ triangulates $P$ and hence
$\overline{P}$ is a triangulated ball. We identify $\overline{P}$
with $P^{\hspace{0.2mm}\prime}$ and also say that $\overline{P}$
is a {\em simplicial subdivision} of $P$. A {\em polytopal
$d$-ball} is a triangulated $d$-ball which is isomorphic to a
simplicial subdivision $\overline{P}$ of some $d$-polytope $P$. }
\end{defn}
(Warning\,: Most authors do not include the hypothesis
$V(P^{\hspace{0.2mm}\prime}) = V(P)$ in the definition of
simplicial subdivision.)

\begin{lemma}\hspace{-1.8mm}{\bf .} \label{L2}
Let $P$ be a polytope with vertex set $V(P)$. Let $A \subseteq
V(P)$. Then \vspace{-2mm}
\begin{enumerate}
\item[{\rm (a)}] either ${\rm conv}(A) \subseteq \bdP$ or $\intcA
\subseteq \intP$; and \vspace{-2mm} \item[{\rm (b)}] if, further,
$P$ is simplicial and ${\rm conv}(A) \subseteq \bdP$ then ${\rm
conv}(A)$ is a proper face of $P$.
\end{enumerate}
\end{lemma}

\noindent {\bf Proof.} (a) Suppose ${\rm conv}(A) \not\subseteq
\bdP$. Then there is a point $u_0$ in ${\rm conv}(A) \cap \intP$.
For any point $x \in \intcA$, $x \neq u_0$, the line $L$ joining $x$
and $u_0$ meets ${\rm conv}(A)$ in a line segment $[a, b]$.
Since $x$ belongs to $\intcA$, it belongs to the relative interior
$(a, b)$ of the segment $[a, b]$. Also, if $L \cap P = [c, d]$ then
$[a, b] \subseteq [c, d]$. Since $u_0 \in L \cap \intP$, it follows
that $(c, d) \subseteq \intP$. Thus, $x \in (a, b) \subseteq (c, d)
\subseteq \intP$. So, $\intcA \subseteq \intP$.

(b) We may assume that ${\rm aff}(P) = \RR^{\hspace{0.2mm}d}$. In
this case ${\rm conv}(A)$ and $\intP$ are disjoint convex sets, of
which the first one is compact and the second one is open in
$\RR^{\hspace{0.2mm}d}$. So, there is a hyperplane $H$ in
$\RR^{\hspace{0.2mm}d}$ strictly separating ${\rm conv}(A)$ and
$\intP$. Then ${\rm conv}(A)$ is contained in the proper face $H
\cap P$ of $P$. Since $P$ is simplicial, it follows that ${\rm
conv}(A)$ is a face of $P$. \hfill $\Box$

\bigskip

Notice that, as a consequence of Lemma \ref{L2}, if
$P^{\hspace{0.2mm} \prime}$ is a simplicial subdivision (in the sense
of Definition \ref{simplicial_subdivision}, which is stronger
than the usual definition) of a simplicial
$d$-polytope $P$, then each simplex in $\partial P^{\hspace{0.2mm}
\prime}$ is a proper face of $P$ and hence (since both $\partial
P^{\hspace{0.2mm}\prime}$ and $\partial P$ are $(d-
1)$-pseudomanifolds without boundary) $\partial P^{\hspace{0.2mm}
\prime} = \partial P$. Thus, if $B$ is a polytopal ball
triangulating a simplicial polytope $P$, then $\partial B$ is
isomorphic to the boundary complex $\partial P$ of $P$.

The following proposition is essentially Theorem 4.1 in \cite{mc2}.

\begin{prop}\hspace{-1.8mm}{\bf .} \label{P10}
Let $B$ be a $k$-stacked triangulated ball of dimension $d+1\geq
2k+ 1$. If $\partial B$ is a polytopal $d$-sphere then $B$ is a
polytopal ball.
\end{prop}

\noindent {\bf Proof.} Since $\partial B$ is polytopal, there is a
$(d+1)$-polytope $P$ in $\RR^{\hspace{0.2mm}d+1}$ such that
$\partial B$ is the boundary complex $\partial P$ of $P$. Thus, we
may identify the vertices of $B$ with those of $P$. For any face
$\alpha \in B$, let $|\alpha|$ denote the convex hull of $\alpha$.
Note that, for $\alpha \in \partial B$, $|\alpha|$ is a proper
face of $P$. It follows that the simplices $|\alpha|$, $\alpha
\in \partial B$,  have pairwise disjoint relative interiors. Indeed,
Lemma \ref{L2} implies that, for $\alpha \in B$ and $\beta\in
\partial B$, $|\alpha|$ and $|\beta|$ have disjoint relative
interiors, whenever $\alpha \neq \beta$.

\smallskip

\noindent {\sf Claim 1\,:} If $\alpha$ is an $i$-face of $B$ then
$|\alpha|$ is a geometric (non-singular) $i$-simplex. That is,
$\dim(|\alpha|) = \dim(\alpha)$ for all $\alpha\in B$.

Suppose there exists an $i$-face $\alpha$ of $B$ such that
$|\alpha|$ is not a geometric $i$-simplex. Then $\alpha$ is a set
of $i+1$ points in the affine space ${\rm aff}(\alpha)$ of
dimension $\leq i-1$. Then, by Radon's Theorem  (cf. \cite[Page
124]{g}), there exist disjoint proper subsets $\beta, \gamma
\subseteq \alpha$ such that $|\beta|^{{}^{\circ}} \cap
|\gamma|^{{}^{\circ}} \neq \emptyset$. Let $\#(\gamma) \leq
\#(\beta)$. Then, $2\#(\gamma) \leq \#(\alpha) \leq d+2 \leq
 2d-2k+2$. Thus, $\#(\gamma) \leq d-k+1$. So, $\dim(\gamma)
\leq d-k$. Therefore, $\gamma \in \partial B$ and hence
(by the comment preceding Claim 1) $|\beta|^{{}^{\circ}} \cap
|\gamma|^{{}^{\circ}} = \emptyset$, a contradiction. This
proves Claim 1.

\smallskip

\noindent {\sf Claim 2\,:} $X := \{|\alpha| \, : \, \alpha\in B\}$
is a geometric simplicial complex.

We have to show that for any two faces $\alpha, \beta$ in $B$,
$|\alpha| \cap |\beta|$ is a common face of both $|\alpha|$ and
$|\beta|$. Otherwise, $I := \{(\alpha, \beta) \in B \times B \, :
\, |\alpha| \cap |\beta|$ is not a common face of $|\alpha|$ and
$|\beta|\}$ is a non-empty set. Fix $(\alpha, \beta) \in I$ such
that $\dim(\alpha)+\dim(\beta)$ is minimum. If $|\alpha| \cap
|\beta|^{{}^{\circ}} = \emptyset$ then there is a proper face
$\beta_1$ of $\beta$ such that $|\alpha| \cap |\beta_1| =
|\alpha|\cap |\beta|$ and hence $(\alpha, \beta_1) \in I$,
contradicting the choice of $(\alpha, \beta)$. So, $|\alpha| \cap
|\beta|^{{}^{\circ}} \neq \emptyset$. Similarly,
$|\alpha|^{{}^{\circ}} \cap |\beta| \neq \emptyset$. Therefore
$|\alpha|^{{}^{\circ}} \cap |\beta|^{{}^{\circ}} \neq \emptyset$
(if $x \in |\alpha| \cap |\beta|^{{}^{\circ}}$, $y \in
|\alpha|^{{}^{\circ}} \cap |\beta|$ then $\frac{x+y}{2} \in
|\alpha|^{{}^{\circ}} \cap |\beta|^{{}^{\circ}}$). Hence
$\alpha, \beta \in B \setminus
\partial B$. Since $B$ is $k$-stacked, it follows that
$\dim(\alpha) \geq d-k+1$, $\dim(\beta) \geq d-k+1$. Hence
$\dim(\alpha) + \dim(\beta) \geq 2d-2k+2\geq d+2$. Since
$|\alpha|, |\beta|$ are simplices in $\RR^{\hspace{0.2mm}d+1}$, it
follows that there is a line $L$, through any given point $x \in
|\alpha|^{{}^{\circ}} \cap |\beta|^{{}^{\circ}}$, contained in
${\rm aff}(\alpha) \cap {\rm aff}(\beta)$. Since the line segments
$L \cap |\alpha|$ and $L \cap |\beta|$ have the interior point $x$
in common, it follows that $[a, b] := L \cap |\alpha| \cap
|\beta|$ is a non-trivial line segment. Thus, $a\neq b$ are points
in the boundary of $|\alpha| \cap |\beta|$. If $a \in
|\alpha|^{{}^{\circ}}$ then $a \in |\beta|^{{}^{\bullet}}$ and
hence $a \in |\beta_1|^{{}^{\circ}}$ for some proper face
$\beta_1$ of $\beta$. Then $(\alpha, \beta_1) \in I$,
contradicting the choice of $(\alpha, \beta)$. So, $a \in
|\alpha|^{{}^{\bullet}}$. Similarly, $b \in
|\alpha|^{{}^{\bullet}}$. Since $x \in [a, b] \cap
|\alpha|{{}^{\circ}}$, $a$, $b$ are not both vertices of $\alpha$.
Assume that $a$ is not a vertex of $\alpha$. Since $a \in
|\alpha|^{{}^{\bullet}}$, $a \in |\alpha_1|^{{}^{\circ}}$ for some
proper face $\alpha_1$ (of dimension $\geq 1$) of $\alpha$. Then
$(\alpha_1, \beta) \in I$, contradicting the choice of $(\alpha,
\beta)$. This completes the proof of Claim 2.

\smallskip

Claim 1 shows that $B$ is the abstract scheme of $X$. Let $|B|$
denote the union of the simplices in
$X$. Since $B$ is a triangulated ball, it follows that $|B|$ is a
topological $(d+1)$-ball. Clearly $|B| \subseteq P$. Since $|B|$
and $P$ are topological $(d+1)$-balls with the same boundary, $|B|
= P$. So, $X$ is a simplicial subdivision of $P$ and is abstractly
isomorphic to $B$.  \hfill $\Box$

\begin{prop}\hspace{-1.8mm}{\bf .} \label{P11}
Let $B$ be a polytopal $k$-stacked ball. Then $B$ does not contain
any standard sphere of dimension $\geq k$ as an induced subcomplex.
\end{prop}

\noindent {\bf Proof.} Let $\dim(B) = d$. We may assume that $B$ is
(the abstract scheme of) a simplicial subdivision of a simplicial
$d$-polytope $P$. For any set $A$ of vertices of $B$, we let
$\langle A\rangle$ denote the convex hull of $A$.

Suppose, if possible, that $\alpha$ is the vertex set of an
induced $S^{\,m}_{m+2}$ in $B$, $m \geq k$. Since the
$d$-pseudomanifold $B$ can't properly contain a $d$-pseudomanifold
without boundary, we must have $k \leq m<d$. Let $\beta$ be an $m$-face
of $B$ contained in $\alpha$. Clearly, there is a facet $\sigma$
of $B$ containing $\beta$ such that $\intcosi \cap \intcoal \neq
\emptyset$. Write $\sigma =\beta\sqcup\mu$. Take a point $x$ in
$\intcosi \cap \intcoal$. Since $x \in \langle\beta \sqcup
\mu\rangle^{{}^{\circ}}$, there are points $b \in \intcobe$, $c
\in \intcomu$ such that $x$ belongs to the open line segment $(b,
c)$. Let the line $bc$ meet $\langle \alpha \rangle$ in the line
segment $[b, b^{\hspace{.2mm}\prime}]$. Since
$b^{\hspace{.2mm}\prime}$ is in the boundary of $\langle
\alpha\rangle$, there is a face $\gamma \subset \alpha$, $\gamma
\not\subseteq \beta$, such that $b^{\hspace{.2mm} \prime} \in
\intcoga$. Since we have the point $x$ in $(b, c)\cap (b,
b^{\hspace{.2mm}\prime})$, it follows that either
$b^{\hspace{.2mm}\prime}\in (b, c)$ or $b^{\hspace{.2mm} \prime} =
c$ or $c\in (b, b^{\hspace{.2mm} \prime})$. If
$b^{\hspace{.2mm}\prime} \in (b, c)$ then $b^{\hspace{.2mm}\prime}
\in \intcosi$ and hence $b^{\hspace{.2mm}\prime} \in \intcosi \cap
\intcoga$. So, $\intcosi \cap \intcoga \neq \emptyset$. But, this
is not possible since $\sigma\neq \gamma$ are faces of $B$. If
$b^{\hspace{.2mm}\prime} = c$ then $b^{\hspace{.2mm} \prime} \in
\intcomu$ and hence $\intcomu \cap \intcoga \neq \emptyset$. This
is also not possible since $\mu \neq \gamma$ are faces of $B$.
Therefore, $c \in (b, b^{\hspace{.2mm} \prime})$.
By Lemma \ref{L2}, $\intcoal \subseteq \intP$. Now, $c$ belongs
to $\intcoal \subseteq \intP$ as well as to $\intcomu$. Thus,
$\intcomu \cap \intP \neq \emptyset$. Then, by Lemma \ref{L2},
$\intcomu \subseteq \intP$ and hence  $\mu \not\in \partial P =
\partial B$. This is a contradiction since $\mu \in B$, $B$ is
a $k$-stacked $d$-ball and $\dim(\mu) = d-m-1 \leq d-k-1$. \hfill
$\Box$

\begin{cor}\hspace{-1.8mm}{\bf .} \label{C6}
Let $S$ be a $k$-stacked polytopal sphere of dimension $d \geq
2k$. Then there is a unique $k$-stacked $(d+1)$-ball
$\overline{S}$ such that $S = \partial \overline{S}$. Further,
$\overline{S}$ is given by the formula $(\ref{eq1})$.
\end{cor}

\noindent {\bf Proof.} Let $B$ be a $k$-stacked
$(d+1)$-ball such that $\partial B = S$. Then, by Proposition
\ref{P10}, $B$ is a polytopal ball.
Therefore, by Proposition \ref{P11}, $B$ contains no induced
standard sphere of dimension $\geq k$. Since $B$ is a $k$-stacked
ball of dimension $\geq 2k+1$, ${\rm skel}_k(B) \subseteq \partial
B = S$. Therefore, $B \subseteq \overline{S}$, where
$\overline{S}$ is defined by formula (\ref{eq1}). If $B \neq
\overline{S}$ then take a minimal face $\alpha \in \overline{S}
\setminus B$. Since $S \subseteq B$, we have $\alpha \in
\overline{S} \setminus S$, and hence $\dim(\alpha) \geq k+ 1$.
Therefore, $\alpha$ induces a standard sphere of dimension $\geq
k$ in $B$, a contradiction. Thus $B = \overline{S}$. \hfill $\Box$

\begin{remark}\hspace{-1.8mm}{\bf .} \label{R2}
{\rm  The uniqueness statement in Corollary \ref{C6} is due to
McMullen (\cite[Theorem 3.3]{mc2}). However, the explicit
description of the ball $\overline{S}$ given above appears to be
new. Indeed, after his Theorem 4.1, McMullen remarks\,: ``The
implication of Theorem 4.1 for the equality case (Conjecture 3.1)
of GLBC is obvious - all we need is an appropriate combinatorial
triangulation of our polytope $P$. However, there are no
corresponding pointers to finding such a triangulation."
Corollary \ref{C6} provides such ``pointers", answering the
question implicit in McMullen's remark. }
\end{remark}

\begin{prop}\hspace{-1.8mm}{\bf .} \label{P12}
Let $S$ be a $(k+1)$-neighbourly polytopal sphere of dimension
$d$. Then $S$ is $(d-k)$-stellated.
\end{prop}

\noindent {\bf Proof.} Fix a vertex $x$ of $S$. Let $A$ be the
antistar of $x$ in $S$. By Bruggesser-Mani (cf. \cite[Theorem
8.12]{z}) $A$ is a shellable $d$-ball. Hence, $x \ast A$ is a
shellable $(d+1)$-ball. Clearly, $\partial(x\ast A) = S$. Since
$S$ is $(k+1)$-neighbourly, $x \ast A$ is $(d-k)$-stacked. Hence,
by Proposition \ref{P2}, $x \ast A$ is $(d-k)$-shelled. Therefore,
by Corollary \ref{C1}, $S$ is $(d-k)$-stellated. \hfill $\Box$

\bigskip

Notice that, as a particular case of Proposition \ref{P12}, any
$(k+1)$-neighbourly polytopal sphere of dimension $2k+1$ is $(k+
1)$-stellated. Also, every polytopal $d$-sphere is $d$-stellated.

The case $k=1$ of the following result is due to M. A. Perles (cf.
\cite[Theorem 1]{as}).

\begin{prop}\hspace{-1.8mm}{\bf .} \label{P13}
Let $S$ be a $(k+1)$-neighbourly polytopal sphere of dimension
$2k+1$. Then $S \in {\cal W}_k(2k+1)$.
\end{prop}

\noindent {\bf Proof.} Let $S$ be the boundary complex of a
simplicial polytope $P$. Then $P$ is a $(k+ 1)$-neighbourly 
$(2k+2)$-polytope. Fix a vertex $v$
of $S$, and let $L = {\rm lk}_S(v)$. We need to prove that $L$ is
$k$-stellated. This is trivial if $S$ is a standard sphere. So,
assume that $P$ is not a simplex. It
follows that $Q := {\rm conv}(V(P) \setminus \{v\})$ is also a
$(2k + 2)$-dimensional polytope. Clearly, $Q$ is also
$(k+1)$-neighbourly and hence, by Radon's Theorem, $Q$ is also
simplicial. Let $B$ be the pure simplicial complex of dimension
$2k+1$ whose facets are those facets of the
polytope $Q$ which are visible from the point $v$. By
Bruggesser-Mani (cf. \cite[Theorem 8.12]{z}), $B$ is a shellable
ball. Clearly, $\partial B = L = S \cap B$. Let $\alpha$ be a $k$-face of
$B$. Then (as $S$ is $(k+1)$-neighbourly) $\alpha \in S \cap B = \partial B$.
Thus, $B$ is $k$-stacked.
Hence $L = {\rm lk}_S(v)$ is $k$-stellated by Propositions
\ref{P2} and \ref{P4}. Since $v$ was an arbitrary vertex of $S$,
it follows that $S \in {\cal W}_k(2k+1)$. \hfill $\Box$

\section{The \boldmath{$g$}-, beta- and mu-vectors and tightness}

For a $d$-dimensional simplicial complex $X$, $f_i = f_i(X)$
denotes the number of $i$-dimensional faces of $X$ ($-1 \leq i\leq
d$). Thus, $f_{-1} = 1$, corresponding to the empty face of $X$.
The vector {\boldmath $f$}$(X) =(f_0, f_1, \dots, f_d)$ is called
the {\em face-vector} (or {\em $f$-vector}) of $X$.

We recall that a pure $d$-dimensional simplicial complex $X$ is
shellable if it may be obtained from the standard $d$-ball by a
finite sequence of shelling moves. Thus, $X$ is shellable if there
is a shelling sequence $B^{\hspace{.2mm}d}_{d+1} = X_0 \subset
X_1 \subset \cdots \subset X_n = X$ of (necessarily pure)
simplicial complexes such that, for $0\leq i \leq n-1$, $X_{i+1}$
is obtained from $X_i$ by a single shelling move. Clearly, each
shelling move of index $j-1$ increases the number of $i$-faces of
a $d$-dimensional simplicial complex by ${d-j+1 \choose i-j+1}$.
Therefore, if, amongst a sequence of shelling moves used to
obtain $X$ from $B^{\hspace{.2mm}d}_{d+1}$, exactly $h_j$ are of
index $j-1$ ($0\leq j \leq d+1$), then the face-vector of $X$ is
given by
$$
f_i(X) = \sum_{j=0}^{i+1} {d-j+1 \choose i-j+1}h_j, ~~ -1 \leq i
\leq d.
$$
(Here, by convention, $h_0 =1$, and the term with $j=0$ gives the
number of $i$-faces in the initial standard $d$-ball.)

Inverting this system of linear equations, we find that the
numbers $h_j$ are given in terms of the face-vector of $X$ by the
formula
\begin{equation} \label{eq3}
h_j = \sum_{i=-1}^{j-1} (-1)^{j-i-1}{d-i \choose j-i-1}f_i, ~~ 0
\leq j \leq d+1.
\end{equation}
This formula shows that the vector {\boldmath $h$}$(X) =(h_0,
\dots, h_{d+1}) = (h_0(X), \dots, h_{d+1}(X))$ depends only on the
simplicial complex $X$, and not on the particular sequence of
shelling moves used to obtain $X$. It is called the {\em
$h$-vector} of $X$. More generally, for any simplicial complex $X$
of dimension $d$, the $h$-vector of $X$ is defined in terms of its
$f$-vector by the formula (\ref{eq3}).

The {\em $g$-vector} {\boldmath $g$}$(X) =(g_0, g_1, \dots,
g_{d+1}) = (g_0(X), g_1(X), \dots, g_{d+1}(X))$ of a simplicial
complex $X$ of dimension $d$ is defined in terms of its $h$-vector
by the formula
$$
g_j(X) = h_j(X) - h_{j-1}(X), ~~ 0\leq j \leq d+1
$$
(where $h_{-1}(X) \equiv 0$). In view of (\ref{eq3}), the
$g$-vector of $X$ is given in terms of its $f$-vector by:
\begin{equation} \label{eq4}
g_j(X) = \sum_{i=-1}^{j-1} (-1)^{j-i-1}{d-i+1 \choose j-i-1}f_i(X),
~~ 0 \leq j \leq d+1.
\end{equation}
(\ref{eq4})  may be inverted to obtain the $f$-vector of $X$ in terms
of its $g$-vector\,:
\begin{equation} \label{eq5}
f_i(X) = \sum_{j=0}^{i+1} {d-j+2 \choose i-j+1}g_j(X), ~~ -1 \leq i
\leq d.
\end{equation}

Let $B^{\hspace{.2mm}d}_{d+1} = X_0 \subset X_1 \subset \cdots
\subset X_n = S$ be a shelling sequence for a shellable $d$-sphere
$S$. For $0 \leq j \leq n$, let $\widetilde{X}_j$ be the pure
simplicial complex of dimension $d$ whose facets are those facets
of $S$ which are not in $X_{j-1}$ (with $X_{-1} = \emptyset$).
Then $B^{\hspace{.2mm}d}_{d+1} = \widetilde{X}_n \subset
\widetilde{X}_{n -1} \subset \cdots \subset \widetilde{X}_0 = S$
is another shelling sequence for $S$. Indeed, if $X_j$ is obtained
from $X_{j-1}$ by the shelling move $\alpha_j \leadsto \beta_j$
then $ \widetilde{X}_{j-1}$ is obtained from $ \widetilde{X}_{j}$
by the shelling move $\beta_j \leadsto \alpha_j$. Therefore, if
the original shelling sequence for $S$ involves $h_i$ shelling
moves of index $i-1$, then the reverse sequence involves $h_i$
shelling moves of index $d-i$. Since the $h$-vector of $S$ is
independent of the particular shelling sequence used, this shows
that the $h$-vector of any shellable $d$-sphere satisfies $h_{d+1-
i} = h_i$, $0\leq i \leq d+1$. From the definition $g_i := h_i
-h_{i-1}$, it follows that the $g$-vector of any shellable
$d$-sphere satisfies
\begin{equation} \label{eq6}
g_{d+2-i} = - \, g_i, ~~ 1 \leq i \leq d+1.
\end{equation}
In fact, the
$g$-vector of any triangulated $d$-sphere satisfies (\ref{eq6}).
Indeed, this is equivalent to the famous Dehn-Sommerville
equations for triangulated spheres. Even more generally, the
$g$-vector of any triangulated closed $d$-manifold with Euler
characteristic $\chi$ satisfies Klee's formula (cf. \cite{k})\,:
\begin{equation} \label{eq7}
g_{d+2-i} + g_i = (-1)^{i-1} {d+2 \choose i}(\chi -
\chi(S^{\hspace{.1mm}d})), ~~ 1 \leq i \leq d+1.
\end{equation}
(In particular, any triangulated closed manifold of odd dimension
$d$ satisfies (\ref{eq6}).) However, for $k$-stellated
$d$-spheres, the $g$-vector has a geometric significance which is
lacking in the more general situations. This geometric meaning of
the $g$-vector stems from\,:

\begin{lemma}\hspace{-1.8mm}{\bf .} \label{L3}
Let $X$, $Y$ be two simplicial complexes of dimension $d$. If $Y$
is obtained from $X$ by a single bistellar move of index $l$ then,
for $0\leq j \leq d$,
$$
g_{j+1}(Y) - g_{j+1}(X) = \left\{
\begin{array}{rl}
+1 & \mbox{ if $j = l \neq d/2$} \\
-1 & \mbox{ if $j = d-l \neq d/2$} \\
0  & \mbox{ otherwise}.
\end{array}
\right.
$$
\end{lemma}

\noindent {\bf Proof.} Notice that a bistellar move of index $l$
creates ${d+1-l \choose i-l}$ new $i$-faces, and destroys ${l+1
\choose d-i+1}$ old $i$-faces. Thus $f_i(Y) - f_i(X) = {d+1-l
\choose i-l} - {l+1 \choose d-i+1}$ for $-1 \leq i \leq d$. Hence
the formula (\ref{eq4}) for the $g$-vector yields
$$
g_{j+1}(Y) - g_{j+1}(X) = \sum_{i=-1}^{j} (-1)^{j-i}{d-i+1 \choose
j-i}\left[{d+1-l \choose i-l} - {l+1 \choose d-i+1}\right].
$$
Now,
\begin{eqnarray*}
\sum_{i=-1}^{j} (-1)^{j-i}{d-i+1 \choose j-i}{l+1 \choose d-i+1} &
= & \sum_{i=-1}^{j}(-1)^{j-i}{l+1\choose d-j+1}{l+j-d\choose j-i} \\
& = &  {l+1 \choose d-j+1}\sum_{i=0}^{j+1} (-1)^{i}{l+j-d \choose
i} \\ & = &  (-1)^{j+1}{l+1 \choose d-j+1}{l+j-d -1\choose j+1}.
\end{eqnarray*}
Since $0 \leq j \leq d$ and $0 \leq l \leq d$, this means that
$$
\sum_{i=-1}^{j} (-1)^{j-i}{d-i+1 \choose j-i}{l+1 \choose d-i+1}
= \left\{
\begin{array}{cl}
+1 & \mbox{ if $j = d-l$} \\
0  & \mbox{ otherwise}.
\end{array}
\right.
$$
Replacing $l$ by $d-l$ in this formula, we get
$$
\sum_{i=-1}^{j} (-1)^{j-i}{d-i+1 \choose j-i}{d+1-l \choose i-l} =
\left\{
\begin{array}{cl}
+1 & \mbox{ if $j = l$} \\
0  & \mbox{ otherwise}.
\end{array}
\right.
$$
Hence the result. \hfill $\Box$

\bigskip

The following result actually holds for any triangulated $(d+1)$-ball
(cf. \cite[Corollary 2]{mc2}). Note that this result implies (\ref{eq6}).

\begin{cor}\hspace{-1.8mm}{\bf .} \label{C7}
If $B$ is a shellable ball of dimension $d+1$, then
$g_{j}(\partial B) = h_j(B) - h_{d+2-j}(B)$, $0 \leq j \leq d+1$.
\end{cor}

\noindent {\bf Proof.} Induction on the number of facets of $B$.
This is trivial if $B = B^{\hspace{.1mm}d+1}_{d+2}$. Else $B$ is
obtained from a shellable ball $B^{\hspace{.2mm}\prime}$ with one
less facet by a shelling move, say of index $l$ ($0\leq l \leq
d$). From the definition of the $h$-vector (for a shellable ball),
we have $h_j(B) - h_j(B^{\hspace{.2mm}\prime}) = \delta_{j, l+1}$,
$h_{d + 2 - j}(B)- h_{d+2-j}(B^{\hspace{.2mm}\prime}) = \delta_{j,
d-l+1}$. (Here $\delta_{i,j}= \delta_{ij}$ is Kronecker delta\,:
$\delta_{ij} = 1$ if $i=j$ and $\delta_{ij}=0$ if $i\neq j$.)
Also, $\partial B$ is obtained from $\partial
B^{\hspace{.2mm}\prime}$ by a bistellar move of index $l$ (by
Lemma \ref{L1}), and hence, by Lemma \ref{L3}, $g_j(\partial B) -
g_j(\partial B^{\hspace{.2mm}\prime}) = \delta_{j, l+1} -
\delta_{j, d-l+1}$. Therefore, we get\,: $g_j(\partial B) - h_j(B)
+ h_{d + 2 - j}(B) = g_j(\partial B^{\hspace{.2mm}\prime}) -
h_j(B^{\hspace{.2mm}\prime})  + h_{d + 2 - j}(B^{\hspace{.2mm}
\prime}) = 0$ by induction hypothesis. \hfill $\Box$

\bigskip

Now, if $S$ is a $k$-stellated sphere of dimension $d \geq 2k-1$,
then, in a sequence of bistellar moves of index $<k$ used to
obtain $S$ from $S^{\hspace{.1mm}d}_{d+2}$, the contribution to
the $g$-vector by an $l$-move is never cancelled by a $(d-l)$-move
(and there is no $(d/2)$-move in the sequence). Therefore, as an
immediate consequence of Lemma \ref{L3}, we have\,:

\begin{prop}\hspace{-1.8mm}{\bf .} \label{P14}
Let $S$ be a $k$-stellated sphere of dimension $d \geq 2k-1$. Then
any sequence of bistellar moves of index $<k$ used to obtain $S$
from $S^{\hspace{.1mm}d}_{d+2}$ contains exactly $g_{j+1}(S)$
moves of index $j$ $(0\leq j <k)$. Hence the length of $S$ is
given by the formula $l(S) = h_k(S) -1$.
\end{prop}

\begin{cor}\hspace{-1.8mm}{\bf .} \label{C8}
Let $S$ be a $k$-stellated $d$-sphere. Then the $g$-vector of $S$
satisfies $g_{j} = 0$ for $k+1 \leq j \leq d-k+1$.
\end{cor}

\noindent {\bf Proof.} This is vacuous unless $d\geq 2k$. So, we
may assume that $d \geq 2k$. First suppose $k+1\leq j <\frac{d}{2}
+1$. Fix a sequence of bistellar moves of index $< k$ used to
obtain $S$ from $S^{\hspace{.1mm}d}_{d+2}$. Since $S$ is
$k$-stellated and $j > k$, $S$ is also $j$-stellated. Since $d
\geq 2j- 1$, Proposition \ref{P14} implies that this given
sequence contains exactly $g_j$ moves of index $j-1$. But, as $j-1
\geq k$, it contains no move of index $j-1$. Thus, $g_j = 0$ for
$k+1\leq j <\frac{d}{2} +1$. If $\frac{d}{2} +1 < j \leq d-k+1$,
then $k+1 \leq d+2-j < \frac{d}{2} +1$ and by formula (\ref{eq6}),
we have $g_j = - g_{d+2-j} = 0$ in this case. Finally, if $j =
\frac{d}{2} +1$ then, by formula (\ref{eq6}), $g_j = 0$. This
completes the proof. \hfill $\Box$

\begin{remark}\hspace{-1.8mm}{\bf .} \label{R3}
{\rm (a) More generally, Corollary \ref{C8} holds for $k$-stacked
spheres. This may be deduced from Proposition 9.1 in \cite{bd9}
using (\ref{eq4}). Therefore, the formula (\ref{eq6}) implies that
the entire $g$-vector of a $k$-stacked sphere of dimension $d \geq
2k$ is determined by the $k$ numbers $g_1, g_2, \dots, g_k$.
(Notice that $g_0 =1$ for any triangulated sphere.) Equivalently,
the $f$-vector of such a sphere is determined by the $k$ numbers
$f_0, f_1, \dots, f_{k-1}$.

(b) In particular, any $k$-stacked sphere of dimension $\geq 2k$
has $g_{k+1}=0$. In \cite{mw}, McMullen and Walkup posed the
famous generalized lower bound conjecture (GLBC)\,: any
triangulated sphere of dimension $\geq 2k+1$ satisfies $g_{k+1}
\geq 0$ with equality (if and) only if the sphere is $k$-stacked.
Actually, McMullen and Walkup originally posed this conjecture
only for polytopal spheres. In this case, the inequality of GLBC
was proved by Stanley \cite{st}. Later McMullen gave a simpler
proof in \cite{mc1}. The equality case of GLBC remains open even
for polytopal spheres. Notice that, by (\ref{eq6}), any
triangulated sphere of dimension $d=2k$ satisfies $g_{k+1} = 0$.
So, the hypothesis $d \geq 2k+1$ in the GLBC is essential.

(c) Recall that a simplicial complex $X$ is said to be {\em
$l$-neighbourly} if any $l$ vertices of $X$ form a face of $X$,
i.e., if $f_{l-1} = {n \choose l}$ with $n = f_0(X)$
(equivalently, if $f_i(X) = {n \choose i+1}$ for all $i \leq
l-1$). Using (\ref{eq4}), it is easy to see that if a simplicial
complex $X$ of dimension $d$ is $l$-neighbourly then $g_l(X) =
{n+l-d-3 \choose l}$, where $n = f_0(X)$. If an $n$-vertex
triangulated $d$-sphere $S$ is $l$-neighbourly, where $l \geq
\frac{d}{2} +1$ and $S$ is not a standard sphere (so that $l \leq
n-2$, $n \geq d+3$), then (as $d+2-l \leq l$) $S$ is also
$(d+2-l)$-neighbourly, so that we get $g_{l} + g_{d+2-l} =
{n+l-d-3 \choose l} + {n-l-1 \choose d+2 -l} > 0$, contradicting
(\ref{eq6}). Thus, if a triangulated $d$-sphere is not the
standard sphere $S^{\hspace{.1mm}d}_{d+2}$ then it can be at most
$\lfloor\frac{d+1}{2}\rfloor$-neighbourly. The $\lfloor\frac{d +
1}{2}\rfloor$-neighbourly triangulated $d$-spheres are called the
{\em upper bound spheres} since they attain the componentwise
maximum among all the face-vectors of triangulated $d$-spheres
with a given number of vertices. }
\end{remark}

\begin{cor}\hspace{-1.8mm}{\bf .} \label{C9}
Let $S$ be a $k$-stellated $d$-sphere which is not the standard
$d$-sphere. Then $S$ is at most $k$-neighbourly.
\end{cor}

\noindent {\bf Proof.} Suppose $S$ is $(k+1)$-neighbourly. Then,
by Remark \ref{R3} (c) above, $d \geq 2k+1$. Hence, by Corollary
\ref{C8} and Remark \ref{R3} (c), we get ${n+k-d-2 \choose k+1} =
g_{k+1} = 0$, where $n = f_0(S)$. Hence $n+k-d-2 < k+1$, i.e., $n
< d+3$. Hence $n = d+2$ and therefore $S =
S^{\hspace{.1mm}d}_{d+2}$. \hfill $\Box$

\bigskip

The following lemma is perhaps well known. For a more conceptual
proof in the case of polytopal spheres, see \cite[Theorem
5.1]{mc2}.

\begin{lemma}\hspace{-1.8mm}{\bf .} \label{L4}
If $X$ is a simplicial complex of dimension $d$, then
$$
\sum_{x \in V(X)} g_j({\rm lk}_X(x)) = (d+2-j)g_j(X) + (j+1)
g_{j+1}(X) ~ \mbox{ for } ~ 0\leq j \leq d.
$$
\end{lemma}

\noindent {\bf Proof.} A simple two-way counting yields
$\sum_{x \in V(X)} f_i({\rm lk}_X(x)) = (i+2)f_{i+1}(X)$.
Therefore, we get\,:
\begin{eqnarray*}
\sum_{x \in V(X)} g_j({\rm lk}_X(x)) & = & \sum_{i=-1}^{j-1}
(-1)^{j-i-1}{d-i \choose j-i-1}\sum_{x\in V(X)} f_i({\rm lk}_X(x)) \\
& = & \sum_{i=-1}^{j-1}(-1)^{j-i-1}{d-i \choose j-i-1} (i+2) f_{i+1}(X) \\
& = & \sum_{i=0}^{j}(-1)^{j-i}{d-i +1\choose j-i} (i+1) f_{i}(X) \\
& = & \sum_{i=-1}^{j}(-1)^{j-i}{d-i +1\choose j-i} (i+1) f_{i}(X) \\
& = & \sum_{i=-1}^{j}(-1)^{j-i}\left[(j+1){d-i+1 \choose j-i} -
(d+2-j){d-i +1\choose j-i-1}\!\right]\!f_{i}(X) \\
& = & (d+2-j)\sum_{i=-1}^{j-1}(-1)^{j-i-1}{d-i +1\choose j-i-1}f_{i}(X) \\
& & ~~~~ + (j+1)\sum_{i=-1}^{j}(-1)^{j-i}{d-i+1 \choose j-i}f_{i}(X) \\
& = & (d+2-j)g_j(X) + (j+1)g_{j+1}(X).
\end{eqnarray*}
\hfill $\Box$

\begin{prop}\hspace{-1.8mm}{\bf .} \label{P15}
The $g$-vector of any member $M$ of ${\cal W}_k(d)$ satisfies
$$
\frac{g_j}{{d+2 \choose j}} = (-1)^{j-k-1} \frac{g_{k+1}}{{d+2
\choose k+1}} ~ \mbox{ for } ~ k+1 < j \leq d-k+1.
$$
If, further, $d\geq 2k$ and $d$ is even, then the Euler
characteristic $\chi$ of $M$ is given by the formula
$$
(-1)^{k}(\chi - 2) =  \frac{2\,g_{k+1}}{{d+2 \choose k+1}}.
$$
\end{prop}

\noindent {\bf Proof.} The link of any vertex $x$ in $M$ is a
$k$-stellated sphere of dimension $d-1$. Therefore, by Corollary
\ref{C8}, we have $g_j({\rm lk}_M(x)) = 0$ for $k+1 \leq j \leq
d-k$. Hence Lemma \ref{L4} yields $(d+2-j)g_j + (j+1) g_{j+1} =
0$, i.e.,
$$
\frac{g_{j+1}}{{d+2 \choose j+1}} = - \frac{g_j}{{d+2 \choose j}}
 ~ \mbox{ for } ~ k+1 \leq j \leq d-k.
$$
Hence, by finite induction on $j$, we get the formula
$$
\frac{g_j}{{d+2 \choose j}} = (-1)^{j-k-1} \frac{g_{k+1}}{{d+2
\choose k+1}} ~ \mbox{ for } ~ k+1 \leq j \leq d-k+1.
$$

Now, suppose $d\geq 2k$ is even. Then $k +1 \leq \frac{d}{2}+1
\leq d-k+1$. Hence we have
$$
\frac{g_{\frac{d}{2}+1}}{{d+2 \choose \frac{d}{2}+1}} =
(-1)^{\frac{d}{2}-k} \frac{g_{k+1}}{{d+2 \choose k+1}}.
$$
On the other hand, (\ref{eq7}) contains the formula
$$
\frac{2\,g_{\frac{d}{2}+1}}{{d+2 \choose \frac{d}{2}+1}}
= (-1)^{{\frac{d}{2}}}(\chi -2).
$$
Comparing these two, we get the formula for $\chi$. \hfill $\Box$

\begin{cor}\hspace{-1.8mm}{\bf .} \label{C10}
Let $d\geq 2k$ and $M \in {\cal W}_k(d)$. Then the entire
$g$-vector of $M$ is determined by the $k+1$ numbers $g_i$, $1
\leq i \leq k+1$. Equivalently, the face-vector of $M$ is
determined by the $k+1$ numbers $f_i$, $0\leq i \leq k$.
\end{cor}

\noindent {\bf Proof.} Note that $g_0 =1$. Also, Proposition
\ref{P15} determines $g_{k+2}, \dots, g_{d-k+1}$ in terms of $g_{k
+1}$. Klee's formula (\ref{eq7}) determines $g_{d-k+2}, \dots,
g_{d+1}$ in terms of $g_{1}, \dots, g_{k}$ and the Euler
characteristic $\chi$ of $M$. But $\chi = 0$ if $d$ is odd, and
Proposition \ref{P15} determines $\chi$ in terms of $g_{k+1}$ when
$d$ is even. \hfill $\Box$

\bigskip

\noindent {\bf Notation\,:} For any simplicial complex $X$ with
vertex set $V(X)$ and any set $A$, $X[A]$ will denote the induced
subcomplex of $X$ on the vertex set $A\cap V(X)$. Thus, $X[A] =
\{\alpha\in X \, :\, \alpha \subseteq A\}$.

\begin{defn}\hspace{-1.8mm}{\bf .} \label{beta-vector}
{\rm  Let $X = X^d_m$ be a simplicial complex of dimension $d$ on
$m$ vertices. Let $\FF$ be a field. Then we define the beta-,
sigma-, and mu-vector of $X$ (with respect to $\FF$) as follows.
The {\em beta-vector} of $X$ is the vector $(\beta_0, \beta_1,
\dots, \beta_d)$, where $\beta_i=\beta_i(X)=\beta_i(X; \FF)$ is
the $i^{\rm th}$
Betti number of $X$ with $\FF$-coefficients. That is, $\beta_i(X)
= \dim_{\F}H_i(X; \FF)$, $0\leq i\leq d$. More generally, if $Y$
is a subcomplex of $X$, we use $\beta_i(X, Y)$ to denote
$\dim_{\F}H_i(X, Y; \FF)$, $0\leq i\leq d$. As usual,
$\widetilde{\beta}_i$ will denote the corresponding reduced Betti
numbers. Thus $\widetilde{\beta}_i = \beta_i$ if $i \neq 0$ and
 $\widetilde{\beta}_0 = \beta_0-1$.

The {\em sigma-vector} $(\sigma_0, \sigma_1, \dots, \sigma_d)$
of $X$ (with respect to $\FF$) is defined by
$$
\sigma_i = \sigma_i(X; \FF) = \sum_{j=0}^m \frac{1}{{m\choose j}}
\sum_{A\in {V(X) \choose j}} \widetilde{\beta}_i(X[A]), ~~
0\leq i\leq d.
$$

We define the {\em mu-vector} $(\mu_0, \dots, \mu_d)$ of $X$
(with respect to $\FF$) by\,:
\begin{eqnarray*}
\mu_0 & = & \mu_0(X; \FF) = 1 \\
\mu_i & = & \mu_i(X; \FF) = \delta_{i1} + \frac{1}{m}\sum_{x \in
V(X)} \sigma_{i-1}({\rm lk}_X(x)), ~~ 1 \leq i\leq d.
\end{eqnarray*}
(By a slight abuse of notation, let's write $\emptyset$ for the
trivial simplicial complex whose only face is the empty set. Here,
we have adopted the convention $\widetilde{\beta}_0(\emptyset) = -
1$ and $\widetilde{\beta}_i(\emptyset) = 0$ if $i \neq 0$.
This convention accounts for the Kronecker delta in the definition
of the mu-vector.) }
\end{defn}

\noindent {\bf Notation\,:} We write $\hoop$ for the covering
relation for set inclusion. Thus, for sets $A$ and $B$, $A \hoop
B$ means that $A \subseteq B$ and $\#(B\setminus A) =1$.

\medskip

With this notation, we have\,:

\begin{lemma}\hspace{-1.8mm}{\bf .} \label{L5}
Let $X$ be a $2$-neighbourly simplicial complex of dimension $d$
on $m$ vertices. Then, the mu-vector of $X$ $($with respect to any
field$)$ is given by\,:
$$
\mu_i= \frac{1}{m}\sum_{j=1}^m \frac{1}{{m-1\choose j-1}}
\sum_{\scriptsize
\begin{array}{c}
A, B \subseteq V(X),\\
\#(B) = j, \\ A \hoops B
\end{array}
}
\hspace{-4mm}\beta_i(X[B], X[A]),
~~~ 0\leq i \leq d.
$$
\end{lemma}

\noindent {\bf Proof.} Since $X$ is 2-neighbourly, so is any
induced subcomplex of $X$. Hence, for $A \hoop B \subseteq V(X)$,
$\beta_0(X[B], X[A]) =1$ if $\#(B)=1$, $A = \emptyset$ and
$\beta_0(X[B], X[A]) =0$ otherwise. Hence the formula holds for
$i=0$. So, let $1\leq i \leq d$.  For any $x\in V(X)$, let $L_x$
be the link of $x$ in $X$. Then, each $L_x$ is a simplicial
complex of dimension $d-1$ with exactly $m-1$ vertices. Let $V=
V(X)$ be the vertex set of $X$. Thus, the vertex set of $L_x$ is
$V \setminus\{x\}$. We have
$$
\sigma_{i-1}(L_x) = \sum_{j=1}^m \frac{1}{{m-1 \choose j-1}}
\sum_{A \in {V\setminus\{x\} \choose j-1}} \widetilde{\beta}_{i
-1}(L_x[A]).
$$
But the exact homology sequence for pairs and the excision
theorem yield
$$
\widetilde{\beta}_{i-1}(L_x[A]) = \left\{
\begin{array}{ll}
\beta_i(x \ast L_x[A], L_x[A]) = \beta_i(X[A \sqcup\{x\}], X[A]) &
\mbox{ when } ~  A \neq \emptyset ~ \mbox{ and}\\
 \beta_i(X[A \sqcup\{x\}], X[A])
- \delta_{i1} & \mbox{ when } ~  A = \emptyset.
\end{array}
\right.
$$
Therefore, from the definition of the mu-vector of $X$, we have
$$
\mu_i= \frac{1}{m}\sum_{x\in V}\sum_{j=1}^m \frac{1}{{m-1\choose j-1}}
\sum_{A \in {V\setminus\{x\} \choose j-1}}
\beta_i(X[A\sqcup\{x\}], X[A]).
$$ \hfill $\Box$

\bigskip

The following linear algebra lemma must be well known. But,
we could not find a reference to it in the required form.

\begin{lemma}\hspace{-1.8mm}{\bf .} \label{L6}
Let $V_1 \stackrel{T_1}{\longrightarrow} V_2 \rightarrow \cdots
\rightarrow V_{2m} \stackrel{T_{2m}}{\longrightarrow}  V_{2m+1}$
be an exact sequence of linear transformations between finite
dimensional vector spaces $($involving an even number $2m$ of
arrows$)$. Then $\sum_{i=1}^{2m+1} (-1)^i \dim(V_i) \leq 0$.
Equality holds here if and only if $T_1$ is injective and $T_{2m}$
is surjective.
\end{lemma}

\noindent {\bf Proof.} From the assumed exactness, we have, for $1
< i < 2m+1$, $\dim(V_i) = {\rm rank}(T_i) + {\rm nullity}(T_i) =
{\rm rank}(T_i) + {\rm rank}(T_{i-1})$. Therefore,
\begin{eqnarray*}
\sum_{i=1}^{2m+1} (-1)^i \dim(V_i) &=&
-\dim(V_1) - \dim(V_{2m+1}) + \sum_{i= 2}^{2m} (-1)^i
\, [{\rm rank}(T_i) + {\rm rank}(T_{i -1})] \\
&=& -[\dim(V_1) - {\rm rank}(T_1)] - [\dim(V_{2m+1})- {\rm
rank}(T_{2m})] \\
& \leq & 0 + 0 =0.
\end{eqnarray*}
From this argument, it is immediate that the necessary and
sufficient condition for equality here is that $T_1$ should be
injective and $T_{2m}$ surjective. \hfill $\Box$

\bigskip

The following proposition is our version of the usual
combinatorial Morse theory. In particular, the parts (a) and (b)
of this proposition are the strong and weak Morse inequalities
(averaged over all possible ``regular simplexwise linear"
functions; compare \cite{k95}). For an alternative combinatorial
version of Morse theory, consult \cite{fo}. As we shall see, the
version developed  here is specially suited to the study of
$\FF$-tightness of $\FF$-orientable 2-neighbourly triangulated
closed manifolds.

\begin{prop}\hspace{-1.8mm}{\bf .}  \label{P16}
Let $X$ be a $2$-neighbourly simplicial complex of dimension
$d$. Then the mu- and beta-vectors of $X$ are related as follows
$($with respect to any given field $\FF\,)$. \vspace{-1mm}
\begin{enumerate}
\item[{\rm (a)}] For $0\leq j \leq d$,\, $\sum_{i=0}^j (-1)^{j-i}
\mu_i \geq \sum_{i=0}^j (-1)^{j-i} \beta_i$, with equality for
$j=d$.
 \vspace{-1mm}
\item[{\rm (b)}]  For $0\leq j \leq d$,\, $\mu_j \geq \beta_j$.
 \vspace{-1mm}
\item[{\rm (c)}] The following are equivalent for any fixed index
$j$ $(0\leq j\leq d)$\,:
\begin{enumerate}
\item[$(i)$] $\sum_{i=0}^j (-1)^{j-i} \mu_i = \sum_{i=0}^j (-1)^{j
- i} \beta_i$, and \vspace{-1mm}
\item[$(ii)$] for any induced subcomplex $Y$ of $X$, the morphism
$H_{j}(Y; \FF) \to H_{j}(X; \FF)$ induced by the inclusion map $Y
\hookrightarrow X$ is injective.
\end{enumerate}
 \vspace{-1mm}
\item[{\rm (d)}] The following are equivalent for any fixed index
$j$ $(0\leq j\leq d)$\,: \vspace{-1mm}
\begin{enumerate}
\item[$(iii)$] $\mu_j = \beta_j$, and
\item[$(iv)$] for any induced subcomplex $Y$ of $X$, both the
morphisms $H_{j-1}(Y; \FF) \to H_{j-1}(X; \FF)$ and $H_{j}(Y; \FF)
\to H_{j}(X; \FF)$ induced by the inclusion map $Y \hookrightarrow
X$ are injective.
\end{enumerate}
 \vspace{-1mm}
\item[{\rm (e)}] If, further, $X$ is an $\FF$-orientable closed
manifold, then $\beta_{d-j} = \beta_{j}$ and $\mu_{d-j} = \mu_{j}$
for $0\leq j\leq d$.
\end{enumerate}
\end{prop}

\noindent {\bf Proof.} (a) Fix an index $j$ and subsets $A
\subseteq B$ of $V(X)$. We have the following exact sequence of
relative homology\,: $H_{j}(X[A]) \to H_{j}(X[B]) \to H_{j}(X[B],
X[A]) \to H_{j-1}(X[A]) \to \cdots \to H_{0}(X[A]) \to H_{0}(X[B])
\to H_{0}(X[B], X[A]) \to 0$. If necessary, we may append an extra
$0\to 0$ at the extreme right, to ensure that this exact sequence
has an even number of arrows. Applying Lemma \ref{L6} to this
sequence, we get\,:
\begin{eqnarray} \label{eq8}
\sum_{i=0}^j (-1)^{j-i}\beta_i(X[B], X[A]) & \geq &
\sum_{i=0}^j (-1)^{j-i}\left(\beta_i(X[B]) - \beta_i(X[A])\right)
\end{eqnarray}
for all pairs $A \subseteq B$ of subsets of $V(X)$.

Since the extreme right arrow in the above sequence is trivially a
surjection, Lemma \ref{L6} says that, for any given pair
$A\subseteq B$, equality in (\ref{eq8}) holds if and only if the
morphism $H_{j}(X[A]) \to H_{j}(X[B])$ induced by the inclusion
map $A \hookrightarrow B$ is an injection.

Now, in view of Lemma \ref{L5}, taking the appropriate weighted
sum of the inequalities (\ref{eq8}) over all pairs $(A, B)$ with
$A \hoop B$, we get
\begin{eqnarray} \label{eq9}
\sum_{i=0}^j (-1)^{j-i}\mu_i & \geq & \frac{1}{m} \sum_{l=1}^m
\frac{1}{{m-1 \choose l-1}} \sum_{\scriptsize
\begin{array}{c}
A, B \subseteq V(X),\\
\#(B) = l, \\ A \hoops B
\end{array}
} \sum_{i=0}^j(-1)^{j-i}\left(\beta_i(X[B]) - \beta_i(X[A])\right).
\end{eqnarray}
Here $m = \#(V(X))$.

Equality holds in (\ref{eq9}) if and only if $H_j(X[A]) \to
H_j(X[B])$ is injective for all pairs $(A, B)$ with $A \hoops
B\subseteq V(X)$. In particular, since $X$ is $d$-dimensional
and each $d$-cycle of $X[A]$ is a $d$-cycle of $X[B]$, equality
holds in (\ref{eq9}) for $j = d$. The right hand side of
(\ref{eq9}) may be written as $\sum_{i= 0}^j (-1)^{j-i}
\sum_{C \subseteq V(X)} \alpha(C) \beta_i(X[C])$, where the
coefficients $\alpha(C)$ are given in terms of $n :=\#(C)$ by the
formula
$$
\alpha(C) \, = \, \frac{1}{m} \frac{1}{{m-1 \choose n-1}}
\sum_{\scriptsize
\begin{array}{c}
A \subseteq V(X),\\
 A \hoops C
\end{array}}\hspace{-4mm}1 ~
- ~ \frac{1}{m} \frac{1}{{m-1 \choose n}}
\sum_{\scriptsize
\begin{array}{c}
B \subseteq V(X),\\
C \hoops B
\end{array}}\hspace{-4mm}1
 \, = \, \frac{1}{m} \left(\frac{n}{{m-1 \choose n-1}} -
\frac{m-n}{{m-1 \choose n}}\right),
$$
where the first term occurs only for $n>0$, and the second term
occurs only for $n<m$. This simplifies to
$$
\alpha(C) = \left\{\begin{array}{rl}
+1 & \mbox{if } ~ C =V(X) \\
-1 & \mbox{if } ~ C =\emptyset \\
0 & \mbox{otherwise}.
\end{array}
\right.
$$
Therefore, the right hand side of (\ref{eq9}) simplifies to
${\displaystyle \sum_{i=0}^j(-1)^{j-i}(\beta_i(X) -
\beta_i(\emptyset)) = \sum_{i=0}^j(-1)^{j-i}\beta_i.}$

(b) We have,
$$
\mu_j = \sum_{i=0}^{j-1} (-1)^{j-1-i} \mu_i +
\sum_{i=0}^{j} (-1)^{j-i} \mu_i \geq
\sum_{i=0}^{j-1} (-1)^{j-1-i} \beta_i +
\sum_{i=0}^{j} (-1)^{j-i} \beta_i = \beta_j.
$$

(c) By the proof of part (a), we see that the equality $(i)$ holds
if and only if $H_{j}(X[A]) \to H_{j}(X[B])$ is injective for all
pairs $A \hoop B \subseteq V(X)$. Now, let $Y$ be an induced
subcomplex of $X$, say with vertex set $A$. Take a sequence $A =
A_1 \hoop A_2 \hoop \cdots \hoop A_n = V(X)$. Then the morphism
$H_{j}(Y) \to H_{j}(X)$ is the composition of the morphisms
$H_{j}(X[A_1]) \to H_{j}(X[A_2]) \to \cdots \to H_{j}(X[A_n])$. If
$(i)$ holds then $H_{j}(Y) \to H_{j}(X)$, being a composition of
injective morphisms, is itself injective. Thus, $(ii)$ holds.
Conversely, if $(ii)$ holds, then for any pair $A_1 \hoop A_2
\subseteq V(X)$, choose $Y=X[A_1]$. Then the composition of the
above sequence of morphisms is injective. So, the first morphism
$H_{j}(X[A_1]) \to H_{j}(X[A_2])$ in this sequence must be
injective. Therefore $(i)$ holds.

(d) From the proof of part (b), we see that $(iii)$ holds if and
only if $\sum_{i=0}^{j-1} (-1)^{j-1-i} \mu_i = \sum_{i=0}^{j-1}
(-1)^{j-1-i} \beta_i$ and $\sum_{i=0}^{j} (-1)^{j-i} \mu_i =
\sum_{i=0}^{j} (-1)^{j-i} \beta_i$. Therefore, part (d) follows
from part (c).

(e) In this case, Alexander duality yields $\beta_{d-i}= \beta_i$.
For any $A \subseteq V(X)$, let $A^c$ denote the complement of $A$
with respect to $V(X)$. Observe that $(A, B) \mapsto (B^c, A^c)$
is a permutation of the set of all pairs $(A, B)$ with $A \hoop B
\subseteq V(X)$. Also, by Alexander duality, we have $\beta_{d-
i}(X[B], X[A]) = \beta_i(X[A^c], X[B^c])$. Therefore, Lemma
\ref{L5} yields $\mu_{d-i} = \mu_i$. \hfill $\Box$

\bigskip

Now we recall\,:

\begin{defn}\hspace{-1.8mm}{\bf .} \label{tight}
{\rm Let $X$ be a $d$-dimensional simplicial complex and $\FF$ be
a field. We say that $X$ is {\em tight with respect to} $\FF$ (or,
in short, {\em $\FF$-tight}) if (i) $X$ is connected, and (ii) for
all induced subcomplexes $Y$ of $X$ and for all $0\leq j \leq d$,
the morphism $H_{j}(Y; \FF) \to H_{j}(X; \FF)$ induced by the
inclusion map $Y \hookrightarrow X$ is injective.

Note that, for fields $\FF_1 \subseteq \FF_2$, $X$ is
$\FF_1$-tight if and only if $X$ is $\FF_2$-tight. Therefore, in
studying $\FF$-tightness, we may, without loss of generality,
restrict to prime fields $\FF$, i.e., $\FF=\QQ$ or $\FF=\ZZ_p$,
$p$ prime. Moreover, for any simplicial complex $X$, the following
are equivalent\,: (a) $X$ is $\FF$-tight for all fields $\FF$, (b)
$X$ is $\ZZ_p$-tight for all primes $p$, and (c) $X$ is
$\QQ$-tight. In view of this observation, we shall say that $X$ is
{\em tight} if it is $\QQ$-tight. }
\end{defn}

Clearly, if $X$ is $\FF$-tight then so is every connected induced
subcomplex of $X$. From this trivial observation, it is easy to
see that the standard sphere $S^{\hspace{.1mm}d}_{d + 2}$ is the
only tight triangulated $d$-sphere, and the standard ball
$B^{\hspace{.1mm}d}_{d +1}$ is the only tight triangulated
$d$-ball. We also have\,:

\begin{prop}\hspace{-1.8mm}{\bf .}  \label{P17}
Let $X$ be an $\FF$-tight simplicial complex $($for some field
$\FF)$.
 \vspace{-1mm}
\begin{enumerate}
\item[{\rm (a)}] If $X$ is $(k-1)$-connected in the sense of
homotopy $($for some $k \geq1)$ then $X$ is $(k+1)$-neighbourly.
 \vspace{-1mm}
\item[{\rm (b)}] If $X$ is a triangulated closed manifold, then
$X$ is $\FF$-orientable.
\end{enumerate}
\end{prop}

\noindent {\bf Proof.} (a) Suppose not. Let $l$ be the smallest
integer such that $X$ is not $(l+1)$-neighbourly. We have $1\leq l
\leq k$. The induced subcomplex $X[\alpha]$ of $X$ on any missing
$l$-face $\alpha$ is an $S^{\hspace{.2mm}l-1}_{l+1}$. Since
$H_{l-1}(X[\alpha]; \FF) \to H_{l-1}(X; \FF)$ is injective, it
follows that $\widetilde{H}_{l-1}(X; \FF) \neq 0$. This is a
contradiction since $l\leq k$ and $X$ is $(k-1)$-connected.

(b) This is trivial if $\dim(X) = 1$. So, assume $d := \dim(X)
\geq 2$. Fix a vertex $x$ of $X$, and let $L$ be the link of $x$ in
$X$. Then, $L$ is a homology $(d-1)$-sphere. In particular,
$\beta_{d-1}(L) = 1 >0$. Thus, all the terms in the sum defining
$\sigma_{d-1}(L)$ (cf. Definition \ref{beta-vector}) are
non-negative, and at least one is positive. Hence $\sigma_{d-1}(L)
> 0$ for any vertex link $L$ of $X$. Thus, all the terms in the
sum defining $\mu_d(X)$ are $>0$. Therefore, the mu-vector of $X$
satisfies $\mu_d >0$. On the other hand, if $X$ is not
$\FF$-orientable, then $\beta_d =0$. Thus $\mu_d \neq \beta_d$.
Also, by part (a), $X$ is 2-neighbourly. Therefore, by Proposition
\ref{P16}\,(d), there is an induced subcomplex $Y$ of $X$ such
that $H_{d-1}(Y; \FF) \to H_{d-1}(X; \FF)$ is not injective (note
that $H_{d}(Y; \FF) \to H_{d}(X; \FF)$ is always injective). This
contradicts the $\FF$-tightness of $X$. \hfill $\Box$

\bigskip

Another way of stating Proposition \ref{P17}\,(b) is that, if $X$
is a triangulated closed manifold which is not orientable (over
$\ZZ$), then $\FF= \ZZ_2$ is essentially the only choice for a
field for which $X$ has a chance of being $\FF$-tight. Note that,
as a special case ($k=1$) of Proposition \ref{P17}\,(a), any
$\FF$-tight simplicial complex is necessarily 2-neighbourly
(whatever the field $\FF$).

Now we have\,:

\begin{prop}\hspace{-1.8mm}{\bf .}  \label{P18}
A simplicial complex $X$ of dimension $d$ is $\FF$-tight if and
only if $X$ is $2$-neighbourly and $\mu_i(X; \FF) = \beta_i(X;
\FF)$ for all indices $i$, $0\leq i \leq d$.
\end{prop}

\noindent {\bf Proof.} This is immediate from Proposition
\ref{P16}\,(d) and Proposition \ref{P17}\,(a). \hfill $\Box$

\section{A tightness criterion for members of \boldmath{${\cal
W}_k(d)$}}

In the following result, $B^{\hspace{.2mm}e}$ stands for an
arbitrary triangulated ball of dimension $e$, and
$S^{\hspace{.2mm}e -1}$ stands for an arbitrary triangulated
sphere of dimension $e-1$. Here $e \geq 0$, and $S^{-1} =
\emptyset$.

\begin{lemma}\hspace{-1.8mm}{\bf .} \label{L7}
Let $X$, $Y$ be simplicial complexes such that $Y = X \cup
B^{\hspace{.2mm}e}$. Suppose $X \cap B^{\hspace{.2mm}e} =
B^{\hspace{.2mm}e-1}$ $($with $e\geq 1)$ or $X \cap
B^{\hspace{.1mm}e} = S^{\hspace{.2mm}e-1}$ $($with $e\geq 0)$.
Then the reduced Betti numbers $($with respect to any field $\FF\,)$
of $X$ and $Y$ are related as follows. \vspace{-1mm}
\begin{enumerate}
\item[{\rm (a)}] If $X \cap B^{\hspace{.2mm}e} = B^{\hspace{.2mm}e-1}$
then $\widetilde{\beta}_i(Y) = \widetilde{\beta}_i(X)$ for all $i$.
\item[{\rm (b)}] If $X \cap B^{\hspace{.2mm}e} = S^{\hspace{.2mm}e-1}$
then
\begin{enumerate}
\item[\mbox{either}] ~  $\widetilde{\beta}_i(Y) -
\widetilde{\beta}_i(X) = \left\{\begin{array}{rll}
+1 & \mbox{if } & i=e \\
0 & \mbox{if } & i\neq e,
\end{array}
\right.$ \item[\mbox{or}] ~ $\widetilde{\beta}_i(Y) -
\widetilde{\beta}_i(X) = \left\{\begin{array}{rll}
-1 & \mbox{if } & i=e-1 \\
0 & \mbox{if } & i\neq e-1.
\end{array}
\right.$
\end{enumerate}
\end{enumerate}
\end{lemma}

\noindent {\bf Proof.} When $e\neq 0$, this is immediate from
Mayer-Vietoris theorem for reduced (simplicial) homology. When
$e=0$, the hypothesis says that $Y$ is the disjoint union of $X$
and a point, so that the result is trivial in this case (and the
first alternative holds). \hfill $\Box$

\begin{lemma}\hspace{-1.8mm}{\bf .} \label{L8}
Let $X$, $Y$ be simplicial complexes of dimension $d$ such that
$Y$ is obtained from $X$ by a single bistellar move $\alpha
\mapsto \beta$, say of index $t$ $(0 \leq t \leq d\,)$. Then, for
any subset $A$ of $V(Y)$, the reduced Betti numbers $($with respect
to any field $\FF\,)$ of $X[A]$ and $Y[A]$ are related as follows.
\vspace{-1mm}
\begin{enumerate}
\item[$(i)$] If $A \supseteq \beta$, $A \cap \alpha = \emptyset$
then \vspace{-1mm}
\begin{enumerate}
\item[\mbox{either}] ~  $\widetilde{\beta}_i(Y[A]) -
\widetilde{\beta}_i(X[A]) = \left\{\begin{array}{rll}
+1 & \mbox{if } & i=t \\
0 & \mbox{if } & i\neq t,
\end{array}
\right.$ \item[\mbox{or}] ~ $\widetilde{\beta}_i(Y[A]) -
\widetilde{\beta}_i(X[A]) = \left\{\begin{array}{rll}
-1 & \mbox{if } & i=t-1 \\
0 & \mbox{if } & i\neq t-1.
\end{array}
\right.$
\end{enumerate}
\vspace{-2mm} \item[$(ii)$] If $A \supseteq \alpha$, $A \cap \beta
= \emptyset$ then \vspace{-2mm}
\begin{enumerate}
\item[\mbox{either}] ~  $\widetilde{\beta}_i(Y[A]) -
\widetilde{\beta}_i(X[A]) = \left\{\begin{array}{rll}
-1 & \mbox{if } & i= d-t \\
0 & \mbox{if } & i\neq d-t,
\end{array}
\right.$
\item[\mbox{or}] ~ $\widetilde{\beta}_i(Y[A]) -
\widetilde{\beta}_i(X[A]) = \left\{\begin{array}{rll}
+1 & \mbox{if } & i=d-t-1 \\
0 & \mbox{if } & i\neq d-t-1.
\end{array}
\right.$
\end{enumerate}
\vspace{-2mm}
\item[$(iii)$] In all other cases, $\widetilde{\beta}_i(Y[A]) =
\widetilde{\beta}_i(X[A])$ for all $i$.
\end{enumerate}
\end{lemma}

\noindent {\bf Proof.} If $A \supseteq \beta$, $A \cap \alpha =
\emptyset$ then we have $Y[A] = X[A] \cup \overline{\beta}$ and
$X[A] \cap \overline{\beta} = \partial \beta$. If $A \supseteq
\alpha$, $A \cap \beta = \emptyset$ then we have $X[A] = Y[A] \cup
\overline{\alpha}$ and $Y[A] \cap \overline{\alpha} = \partial
\alpha$. So, the result is immediate from Lemma \ref{L7} in these
cases. (Remember that $\dim(\beta) =t$ and $\dim(\alpha)=d-t$.)

If $A \supseteq \alpha\cup\beta$ then $Y[A]$ is obtained from
$X[A]$ by the bistellar move $\alpha \mapsto \beta$. Hence
$Y[A]$ and $X[A]$ are homeomorphic in this case, and the
result follows.

If $A$ contains neither $\alpha$ nor $\beta$, then $Y[A] = X[A]$,
and the result is trivial.

If $A \supseteq \beta$ and $\alpha_0 := A\cap \alpha$ is a proper
non-empty subset of $\alpha$, then $Y[A] = X[A]\cup
B^{\hspace{.2mm}e}$ and $X[A]\cap B^{\hspace{.2mm}e} =
B^{\hspace{.2mm}e -1}$, where $e=t + \#(\alpha_0) > 0$,
$B^{\hspace{.2mm}e} = \overline{\alpha_0\cup\beta}$,
$B^{\hspace{.2mm}e-1} = \overline{\alpha_0} \ast \partial \beta$.
Hence the result follows from Lemma \ref{L7}.

If $A \supseteq \alpha$ and $\beta_0 := A\cap \beta$ is a proper
non-empty subset of $\beta$, then $X[A] = Y[A]\cup
B^{\hspace{.2mm}e}$ and $Y[A]\cap B^{\hspace{.2mm}e} =
B^{\hspace{.2mm}e-1}$, where $e=d-t + \#(\beta_0) > 0$,
$B^{\hspace{.2mm}e} = \overline{\alpha\cup\beta_0}$,
$B^{\hspace{.2mm}e-1} = \overline{\beta_0} \ast \partial \alpha$.
Hence the result follows from Lemma \ref{L7} in this case also.
\hfill $\Box$

\bigskip

Now, we are in a position to prove a crucial result on the
sigma-vectors of $k$-stellated spheres\,:

\begin{prop}\hspace{-1.8mm}{\bf .}  \label{P19}
For $k\geq 1$, let $S$ be an $m$-vertex $k$-stellated triangulated
sphere of dimension $d \geq 2k-1$. Then, with respect to any
field, the sigma-vector of $S$ is related to its $g$-vector by\,:
\vspace{-2mm}
\begin{enumerate}
\item[{\rm (a)}] $\sigma_i =0$ \, for \, $k \leq i \leq d-k-1$,
\vspace{-1mm}
\item[{\rm (b)}] ${\displaystyle \sum_{i=0}^{l}(-1)^{l-i}
\sigma_i \leq \frac{m+1}{d+3} \sum_{i=0}^{l+1}(-1)^{l+1-i}
\frac{g_i}{{d+2 \choose i}}}$ \, for \, $0 \leq l \leq k-2$,
\, and \vspace{-1mm}
\item[{\rm (c)}] ${\displaystyle \sum_{i=0}^{l}(-1)^{l-i}
\sigma_i = \frac{m+1}{d+3} \sum_{i=0}^{l+1}(-1)^{l+1-i}
\frac{g_i}{{d +2 \choose i}}}$ \, for \, $k-1 \leq l \leq d-k-1$.
\end{enumerate}
\end{prop}

\noindent {\bf Proof.} Induction on the length $l(S)$ of $S$. If
$l(S) = 0$, then $S= S^{\hspace{.2mm}d}_{d+2}$. In this case,
$\sigma_i(S) = -\delta_{i0}$ for $0\leq i<d$, $g_i(S) =
\delta_{i0}$ for $0\leq i\leq d+1$, and $m=d+2$. So, the result is
trivial in this case. Now, assume $l(S) >0$. Then $S$ is obtained
from a shorter $k$-stellated $d$-sphere $S^{\hspace{.2mm}\prime}$
by a single bistellar move $\alpha \mapsto \beta$, say of index
$t$ ($0\leq t <k$).

For $k \leq i\leq d-k-1$, we have $t < i < d-t-1$ and hence by
Lemma \ref{L8}, $\widetilde{\beta}_i(S[A]) =
\widetilde{\beta}_i(S^{\hspace{.2mm}\prime}[A])$ for all subsets
$A$ of $V(S)$. Taking an appropriate weighted sum of these
equalities over all sets $A$, we get $\sigma_i(S) =
\sigma_i(S^{\hspace{.2mm}\prime}) = 0$ for $k\leq i \leq d-k-1$.
Here the last equality is by induction hypothesis. This proves
part (a).

Fix an index $l$ such that $0\leq l\leq d-k-1$ and let $0\leq i
\leq l$. Then $i \leq d-k-1 < d-t-1$, and hence by Lemma \ref{L8},
$\widetilde{\beta}_i(S[A]) =
\widetilde{\beta}_i(S^{\hspace{.2mm}\prime}[A])$ unless
$A\supseteq \beta$, $A\cap \alpha = \emptyset$. So, let ${\cal A}$
be the collection of all subsets $A$ of $V(S)$ such that
$A\supseteq \beta$ and $A\cap \alpha = \emptyset$. For $A \in
{\cal A}$, we have $\widetilde{\beta}_i(S[A]) =
\widetilde{\beta}_i(S^{\hspace{.2mm}\prime}[A])$ for $i\neq t,
t-1$, and either $\widetilde{\beta}_t(S[A])=
\widetilde{\beta}_t(S^{\hspace{.2mm}\prime}[A]) +1$,
$\widetilde{\beta}_{t-1}(S[A])=
\widetilde{\beta}_{t-1}(S^{\hspace{.2mm}\prime}[A])$, or else
$\widetilde{\beta}_t(S[A])=
\widetilde{\beta}_t(S^{\hspace{.2mm}\prime}[A])$,
$\widetilde{\beta}_{t-1}(S[A])=
\widetilde{\beta}_{t-1}(S^{\hspace{.2mm}\prime}[A])-1$. Let ${\cal
A}^{\hspace{.2mm}+}$ be the set of all $A\in {\cal A}$ for which
the first alternative holds. Then we get\,:

\begin{eqnarray} \label{eq10}
\sum_{i=0}^l (-1)^{l-i}\left[\,\widetilde{\beta}_i(S[A]) -
\widetilde{\beta}_i(S^{\hspace{.2mm}\prime}[A])\right] =
\left\{\begin{array}{cl}
0 & \mbox{if } ~ A \not\in {\cal A} ~ \mbox{or if } ~ l < t-1 \\
0 & \mbox{if } ~ A \in {\cal A}^{\hspace{.2mm}+} \, \mbox{ and  }
~ l = t-1 \\ (-1)^{l-t} & \mbox{otherwise.}
\end{array}
\right.
\end{eqnarray}

First consider the case $l\leq t-1$  (which can occur only for
$l < k-1$ and $t>0$). Then (\ref{eq10}) implies
$$
\sum_{i=0}^l (-1)^{l-i} \widetilde{\beta}_i(S[A]) \leq
\sum_{i=0}^l (-1)^{l-i} \widetilde{\beta}_i(S^{\hspace{.2mm}\prime}[A])
~~ \mbox{for all} ~~ A \subseteq V(S).
$$
Taking the appropriate weighted sum of these inequalities over all
$A$, we get

\begin{eqnarray*}
\sum_{i=0}^l (-1)^{l-i}\sigma_i(S) ~ \leq ~ \sum_{i=0}^l
(-1)^{l-i} \sigma_i(S^{\hspace{.2mm}\prime})
&\leq &\frac{m+1}{d+3} \sum_{i=0}^{l+1}(-1)^{l+1-i}\,
\frac{g_i(S^{\hspace{.2mm}\prime})}{{d+2 \choose i}} \\
&= &\frac{m+1}{d+3} \sum_{i=0}^{l+1}(-1)^{l+1-i}\,
\frac{g_i(S)}{{d+2 \choose i}},
\end{eqnarray*}
where the second inequality is by induction hypothesis and the
final equality holds since by Lemma \ref{L3}, we have $g_i(S) =
g_i(S^{\hspace{.2mm}\prime})$ for $i\leq l+1 \leq t$. This
completes the induction step in this case.

Next consider the case $0< t \leq l \leq d-k-1$. In this case,
(\ref{eq10}) says\,:
$$
\sum_{i=0}^l (-1)^{l-i} \left[\,\widetilde{\beta}_i(S[A]) -
\widetilde{\beta}_i(S^{\hspace{.2mm}\prime}[A])\right] =
\left\{\begin{array}{cl}
(-1)^{l-t} & \mbox{if } ~ A \in {\cal A} \\
0 & \mbox{otherwise.}
\end{array}
\right.
$$
Notice that there are exactly ${m-d-2 \choose j-t-1}$ $j$-sets in
${\cal A}$. Therefore, adding these equations over all $j$-subsets
of $V := V(S)$ we get (for $0\leq j \leq m$)\,:
\begin{eqnarray} \label{eq11}
\sum_{i=0}^l (-1)^{l-i} \sum_{A\in{V \choose j}} \left[\,\widetilde{\beta}_i(S[A]) -
\widetilde{\beta}_i(S^{\hspace{.2mm}\prime}[A])\right] =
(-1)^{l-t}{m-d-2 \choose j-t-1}.
\end{eqnarray}
Dividing this equation by ${m\choose j}$ and adding over all $j$, we get\,:
$$
\sum_{i=0}^l (-1)^{l-i} \left[\sigma_i(S) - \sigma_i(S^{\hspace{.2mm}\prime})\right] =
(-1)^{l-t}\sum_{j=0}^m\frac{{m-d-2 \choose j-t-1}}{{m\choose j}}.
$$
But, we have the computation
\begin{eqnarray*}
\sum_{j=0}^{m} \frac{{m-d-2 \choose j-t-1}}{{m\choose j}}
&= &\sum_{j=t+1}^{m+t-d-1} \frac{{m-d-2 \choose j-t-1}}{{m \choose j}}
= \sum_{j=0}^{m-d-2} \frac{{m-d-2 \choose j}}{{m \choose j+t+1}} \\
&=& (m+1)\sum_{j=0}^{m-d-2} {m-d-2 \choose j} \int_0^1 x^{j+t+1}(1-x)^{m-j-t-1}dx \\
&=& (m+1) \int_0^1 x^{t+1}(1-x)^{d+1-t}dx = \frac{m+1}{d+3}
\frac{1}{{d+2 \choose t+1}},
\end{eqnarray*}
where we have made two uses of Euler's famous identity\,:
$$
\int_0^1 x^a(1-x)^b dx= \frac{1}{(a+b+1){a+b \choose a}}~,
$$
for non-negative integers $a, b$.

So, we have
\begin{eqnarray} \label{eq12}
\sum_{j=0}^{m} \frac{{m-d-2 \choose j-t-1}}{{m\choose j}} =
\frac{m+1}{d+3}\frac{1}{{d+2 \choose t+1}}.
\end{eqnarray}
Thus we get\,:
$$
\sum_{i=0}^l (-1)^{l-i} \left[\sigma_i(S) -
\sigma_i(S^{\hspace{.2mm}\prime})\right] = (-1)^{l-t}\,
\frac{m+1}{d+3}\frac{1}{{d+2 \choose t+1}}.
$$
Therefore, induction hypothesis gives the following inequality for
$0<t\leq l\leq d-k-1$ (with equality for $k-1\leq l\leq d-k-1$)\,:
\begin{eqnarray*}
\sum_{i=0}^l (-1)^{l-i}\sigma_i(S) &\leq & (-1)^{l-t} \,
\frac{m+1}{d+3}\frac{1}{{d+2 \choose t+1}}
+ \frac{m+1}{d+3} \sum_{i=0}^{l+1} (-1)^{l+1-i} \,
\frac{g_i(S^{\hspace{.2mm}\prime})}{{d+2 \choose i}} \\
&= &\frac{m+1}{d+3} \sum_{i=0}^{l+1}(-1)^{l+1-i}\,
\frac{g_i(S^{\hspace{.2mm}\prime}) + \delta_{i, t+1}}{{d+2 \choose
i}}  \\
&= & \frac{m+1}{d+3} \sum_{i=0}^{l+1}(-1)^{l+1-i} \,
\frac{g_i(S)}{{d+2 \choose i}},
\end{eqnarray*}
where the last equality holds since by Lemma \ref{L3}, we have
$g_i(S) = g_i(S^{\hspace{.2mm}\prime}) + \delta_{i, t+1}$ as
$d\geq 2k-1$ (so that $d\neq 2t$) and $i-1 \leq l \leq d-k < d-t$.
This completes the induction step in the second case.

Finally, consider the case $0=t \leq l\leq k-d-1$. In this case,
$\beta$ is a vertex of $S$ not in $S^{\hspace{.2mm}\prime}$. Let
$V^{\hspace{.2mm}\prime} = V\setminus \{\beta\}$ be the vertex set
of $S^{\hspace{.2mm}\prime}$. (Thus, $S^{\hspace{.2mm}\prime}$ has
$m-1$ vertices in this case.) Dividing equation (\ref{eq11}) (with
$t=0$) by ${m \choose j}$ and adding over all $j$ ($0\leq j\leq
m$) we get (in view of (\ref{eq12}))

\begin{eqnarray*}
\sum_{i=0}^l (-1)^{l-i}\sigma_i(S) = (-1)^{l} \, \frac{m+
1}{(d+3)(d+2)} + \sum_{i=0}^l (-1)^{l-i} \sum_{j=0}^m \frac{1}{{m
\choose j}} \sum_{A\in{V \choose j}}
\widetilde{\beta}_i(S^{\hspace{.2mm}\prime}[A]).
\end{eqnarray*}
Now, for each $i\leq l$,
\begin{eqnarray*}
\sum_{j=0}^m \frac{1}{{m \choose j}} \sum_{A\in{V \choose j}}
\widetilde{\beta}_i(S^{\hspace{.2mm}\prime}[A])
&= & \sum_{j=0}^{m-1} \frac{1}{{m \choose j}} \sum_{A\in{V \choose j}}
\widetilde{\beta}_i(S^{\hspace{.2mm}\prime}[A])
\, = \, \sum_{j=0}^{m-1} \frac{1}{{m
\choose j}} \sum_{B\in {V^{\hspace{.1mm}\prime} \choose j} \sqcup
{V^{\hspace{.1mm}\prime} \choose j-1}}
\hspace{-2mm}\widetilde{\beta}_i(S^{\hspace{.2mm}\prime}[B]) \\
&= & \sum_{j=0}^{m-1}
\left[\frac{1}{{m \choose j}} + \frac{1}{{m \choose j+1}}\right]
\sum_{B\in {V^{\hspace{.1mm}\prime} \choose j}} \hspace{-2mm}
\widetilde{\beta}_i(S^{\hspace{.2mm}\prime}[B])
\, = \, \frac{m+1}{m} \, \sigma_i(S^{\hspace{.2mm}\prime}).
\end{eqnarray*}
Since $i\leq l\leq d-k-1 <d$,
$\widetilde{\beta}_i(S^{\hspace{.2mm} \prime}) = 0$. This
justifies the first and third equalities above. The second
equality holds since, by definition, $S^{\hspace{.2mm} \prime}[A]
= S^{\hspace{.2mm}\prime}[B]$, where $B = A \cap V^{\hspace{.1mm}
\prime}$, and the map $A\mapsto A\cap V^{\hspace{.1mm} \prime}$ is
a bijection between ${V \choose j}$ and ${V^{\hspace{.1mm}\prime}
\choose j} \sqcup {V^{\hspace{.1mm}\prime} \choose j-1\,}$. The
last equality is because of the trivial identity $\frac{1}{{m
\choose j}} + \frac{1}{{m \choose j+1}} = \frac{m+1}{m}
\frac{1}{{m -1 \choose j}}$. So, we have (when $t=0$)\,:
$$
\sum_{i=0}^l (-1)^{l-i}\sigma_i(S) = (-1)^{l} \, \frac{m+
1}{(d+3)(d+2)} + \frac{m+1}{m}\sum_{i=0}^l (-1)^{l-i}
\sigma_i(S^{\hspace{.2mm}\prime}), ~ 0\leq l\leq d-k-1.
$$
Now, since $S^{\hspace{.2mm}\prime}$ has $m-1$ vertices, induction
hypothesis gives
$$
\sum_{i=0}^l (-1)^{l-i} \sigma_i(S^{\hspace{.2mm}\prime}) \leq
\frac{m}{d+3}\sum_{i=0}^{l+1} (-1)^{l+1-i} \,
\frac{g_i(S^{\hspace{.2mm}\prime})}{{d+2
\choose i}}\,,
$$
with equality for $k-1 \leq l \leq d-k-1$. Therefore, we get
\begin{eqnarray*}
\sum_{i=0}^l (-1)^{l-i}\sigma_i(S) & \leq & (-1)^{l} \, \frac{m+
1}{(d+3)(d+2)} + \frac{m+1}{d+3}\sum_{i=0}^{l+1} (-1)^{l+1-i} \,
\frac{g_i(S^{\hspace{.2mm}\prime})}{{d+2 \choose i}} \\
&=& \frac{m+1}{d+3}\sum_{i=0}^{l+1} (-1)^{l+1-i} \,
\frac{g_i(S^{\hspace{.2mm}\prime})+ \delta_{i1}}{{d+2 \choose i}}  \\
&=& \frac{m+1}{d+3}\sum_{i=0}^{l+1} (-1)^{l+1-i} \,
\frac{g_i(S)}{{d+2 \choose i}}\,,
\end{eqnarray*}
with equality for $k-1 \leq l\leq d-k-1$. This completes the
induction in the last case, thus proving (b) and (c). \hfill
$\Box$

\bigskip

Now, the following key result on the mu-vector of 2-neighbourly
members of ${\cal W}_k(d)$ is more or less immediate.

\begin{prop}\hspace{-1.8mm}{\bf .}  \label{P20}
Let $M \in {\cal W}_k(d)$ be $2$-neighbourly with $d\geq 2k\geq
2$. Then the mu-vector of $M$ $($with respect to any field$\,)$ is
related to its $g$-vector as follows\,: \vspace{-1mm}
\begin{enumerate}
\item[{\rm (a)}] $\mu_i =0$ \, for \, $k+1 \leq i \leq d-k-1$,
\vspace{-1mm}
\item[{\rm (b)}]  ${\displaystyle
\sum_{i=1}^{l}(-1)^{l-i}\mu_i \leq \frac{g_{l+1}}{{d+2 \choose
l+1}}}$ \, for \, $1 \leq l \leq k-1$, \, and  \vspace{-2mm}
\item[{\rm (c)}] ${\displaystyle \sum_{i=1}^{l}(-1)^{l-i}
\mu_i = \frac{g_{l+1}}{{d+2 \choose l+1}}}$, \,  for \,
$k \leq l \leq d-k-1$.
\end{enumerate}
\end{prop}

\noindent {\bf Proof.} Let $m$ be the number of vertices of $M$
and let $L_x$ be the link of $x$ in $M$ for $x\in V(M)$. Then each
$L_x$ is a $k$-stellated sphere, of dimension $d-1\geq 2k-1$, on
exactly $m-1$ vertices. Therefore, for $k+1 \leq i \leq d-k-1$,
$\sigma_{i-1}(L_x) = 0$ by Proposition \ref{P19}\,(a). Taking the
sum of these equations over all $x\in V(M)$, we get $\mu_i(M) = 0$
for $k+1 \leq i \leq d-k-1$. This proves part (a).

Also, by Proposition \ref{P19}\,(b) and (c), we have, for $1\leq
l\leq d-k-1$,
$$
\sum_{i=1}^l (-1)^{l-i} \, \sigma_{i-1}(L_x) \leq \frac{m}{d+2}
\sum_{i=0}^l (-1)^{l-i} \, \frac{g_i(L_x)}{{d+1 \choose i}},
$$
with equality for $k\leq l \leq d-k-1$. Adding these over all
$x\in V(M)$, and dividing the result by $m$, we get (in view of
Lemma \ref{L4})\,:
$$
\sum_{i=1}^l (-1)^{l-i} \, (\mu_{i} - \delta_{i1}) \leq
\frac{1}{d+2} \sum_{i=0}^l \frac{(-1)^{l-i}}{{d+1 \choose
i}}((d+2-i)g_i + (i+1)g_{i+1}),
$$
with equality for $k\leq l \leq d-k-1$. That is,
$$
\sum_{i=1}^l (-1)^{l-i} \, \mu_{i} \leq (-1)^{l-1} + \sum_{i=0}^l
(-1)^{l-i} \left(\frac{g_i}{{d+2 \choose i}} + \frac{g_{i+1}}{{d+2
\choose i+1}}\right) = \frac{g_{l+1}}{{d+2 \choose l+1}} ~
(\mbox{since } ~ g_0 =1),
$$
with equality for $k\leq l \leq d-k-1$. This proves (b) and (c).
\hfill $\Box$

\bigskip

Now, we can prove one of the main results of this paper.

\begin{prop}[A lower bound theorem for {\boldmath ${\cal
W}_k(d)$}]\hspace{-2mm}{\bf .}  \label{P21} Let $M \in {\cal
W}_k(d)$ be  $2$-\hspace{-1.5mm} neighbourly. Then the $g$-vector
of $M$ is related to its Betti numbers $($with respect to any
field$\,)$ as follows\,: \vspace{-2mm}
\begin{enumerate}
\item[{\rm (a)}] if \, $d = 2k$,  \, then \, $g_{l+1} \geq {d+2
\choose l+1}  {\displaystyle \sum_{i=1}^{l}(-1)^{l-i}\beta_i }$ \,
for \, $1 \leq l \leq k-1$,
 \vspace{-2mm}
\item[{\rm (b)}] if \, $d \geq 2k+1$, \, then \,  $g_{l+1} \geq {d+2
\choose l+1}  {\displaystyle \sum_{i=1}^{l}(-1)^{l-i}\beta_i }$ \,
for \, $1 \leq l \leq k$,
 \vspace{-2mm}
\item[{\rm (c)}] if \, $d \geq 2k+2$, \, then \, $g_{l+1} = {d+2
\choose l+1}  {\displaystyle \sum_{i=1}^{l}(-1)^{l-i}\beta_i }$ \,
for \, $k \leq l \leq d-k-1$, \, and
\item[{\rm (d)}] if \, $d\geq 2k+2$, \, then \, $\beta_{i} = 0$ \,
for \, $k+1 \leq i \leq d-k-1$, \vspace{-2mm}
\end{enumerate}
\end{prop}

\noindent {\bf Proof.} By Propositions \ref{P16}\,(b) and
\ref{P20}\,(a), we have $0\leq \beta_i \leq \mu_i=0$ for $k+1 \leq
i \leq d-k-1$. This proves part (d).

Since $M$ is 2-neighbourly, it is connected. Therefore, $\beta_0 =
1 = \mu_0$. Hence, Proposition \ref{P16}\,(a) yields $\sum_{i=
1}^{l}(- 1)^{l-i} \beta_i \leq  { \sum_{i=1}^{l}(-1)^{l-i} \mu_i
}$ for $l \geq 1$. Therefore, parts (a) and (b) are immediate from
Proposition \ref{P20}\,(b) and (c).

If $k < l \leq d-k-1$, then $\beta_l = \mu_l$ by part (d), and
hence by Proposition \ref{P16}\,(d), $H_l(Y) \to H_l(M)$ is
injective for any induced subcomplex $Y$ of $M$. Hence, by
Proposition \ref{P16}\,(c), ${ \sum_{i= 1}^{l}(- 1)^{l-i} \beta_i
} = {\sum_{i= 1}^{l}(-1)^{l-i}\mu_i }$. But, ${\sum_{i= 1}^{l}(-
1)^{l-i}\mu_i } = g_{l+1}/{d+2 \choose l+1}$ by Proposition
\ref{P20}\,(c). This proves part (c) in case $k<l\leq d-k-1$.
Finally, when $k+1 \leq d-k-1$, we have $\beta_{k+1} =0= \mu_{k+1}$ by
part (d) and Proposition \ref{P20}\,(a), hence $H_k(Y) \to H_k(M)$ is injective for any induced
subcomplex $Y$ of $M$ (by Proposition \ref{P16}\,(d)). Therefore,
by Propositions \ref{P16}\,(c) and \ref{P20}\,(c), ${
\sum_{i=1}^{k}(- 1)^{k-i} \beta_i } = {\sum_{i= 1}^{k}(-1)^{k
-i}\mu_i } = g_{k+1}/{d+2 \choose k+1}$, completing the proof of
part (c).  \hfill $\Box$

\bigskip

We also need the following elementary result.

\begin{lemma}\hspace{-1.8mm}{\bf .} \label{L9}
Let $X$ be an $(l+1)$-neighbourly simplicial complex. Then the
beta- and mu-vectors of $X$ $($with respect to any field\,$)$
satisfy $\beta_i = 0 = \mu_i$ for $1\leq i \leq l-1$.
\end{lemma}

\noindent {\bf Proof.} The $l$-skeleton of $X$ agrees with that of
the standard ball of dimension $f_0(X)-1$. Since the ball is
homologically trivial and the $i^{\rm th}$ homology of a
simplicial complex is the same as the $i^{\rm th}$ homology of
its $(i+1)$-skeleton, it follows that
$\beta_i(X) = 0$ for $1\leq i\leq l-1$.

Also, for any vertex $x$ of $X$, the link $L_x$ of $x$ in $X$ is
$l$-neighbourly. Therefore, by the same argument, we have, for
$1\leq i \leq l-1$, $\widetilde{\beta}_{i-1}$ of any induced
subcomplex of $L_x$ is $=0$, except that $\widetilde{\beta}_{0} =
-1$ for the empty subcomplex. Therefore, taking an appropriate
weighted sum, we get $\sigma_{i-1}(L_x) = - \delta_{i1}$ for
$1\leq i \leq l-1$ and $x\in V(X)$. Adding over all $x\in V(X)$,
we get $\mu_i(X) =0$ for $1\leq i \leq l-1$. \hfill $\Box$

\bigskip

The following result (which is the first known combinatorial
criterion for tightness) is due to K\"{u}hnel \cite{k95}.

\begin{lemma}[K\"{u}hnel]\hspace{-1.5mm}{\bf .} \label{L10}
For $k\geq 1$, let $M$ be a $(k+1)$-neighbourly triangulation of
an $\FF$-orientable closed manifold of dimension $2k$. Then $M$ is
$\FF$-tight.
\end{lemma}
(Note that, when $k\geq 2$, $M$ is by assumption at least
3-neighbourly, and hence simply connected. Thus, the hypothesis of
orientability is automatic for $k\geq 2$.)

\bigskip

\noindent {\bf Proof.} Since $M$ is at least 2-neighbourly, it is
connected. Therefore, $\mu_0 = 1 = \beta_0$. By Lemma \ref{L9},
$\mu_i = 0=\beta_i$ for $1\leq i\leq k-1$. So, by duality
(Proposition \ref{P16}\,(e)), $\mu_i = 0 = \beta_i$ for $k+1 \leq
i \leq 2k-1$ and $\mu_{2k} = 1 = \beta_{2k}$. Thus $\mu_i =
\beta_i$ for all $i$, except possibly for $i=k$. But then, the
equality $\sum_{i=0}^{2k} (-1)^i \mu_i = \sum_{i=0}^{2k} (-1)^i
\beta_i$ from Proposition \ref{P16}\,(a) implies $\mu_i = \beta_i$
for $i=k$ as well. Therefore, by Proposition \ref{P18}, $M$ is
$\FF$-tight. \hfill $\Box$

\bigskip

The ``if\," part of Proposition \ref{P22}\,(a) below is
essentially due to Effenberger \cite{ef}. This paper was largely
motivated by a desire to understand and generalize Effenberger's
result. (In this connection, recall that ${\cal W}_1(d) =  {\cal
K}_1(d)$ by Corollary \ref{C4}.)

\begin{prop}\hspace{-1.8mm}{\bf .}  \label{P22}
Let $M \in {\cal W}_1(d)$. Then we have the following.
\vspace{-2mm}
\begin{enumerate}
\item[{\rm (a)}]  If $d \neq 3$, then $M$ is $\FF$-tight if and
only if $M$ is $2$-neighbourly  and $\FF$-orientable.
\vspace{-2mm}
\item[{\rm (b)}] If $d = 3$, then $M$ is $\FF$ tight if and only
if $M$ is $2$-neighbourly, $\FF$-orientable, and satisfies
$\beta_1(M; \FF) = {(n-4)(n-5)}/{20}$, where $n = f_0(M)$.
\end{enumerate}
\end{prop}

\noindent {\bf Proof.} By Proposition \ref{P17}, to be
$\FF$-tight, $M$ must be 2-neighbourly and $\FF$-orientable. Also,
if $d = 3$ and $M$ is 2-neighbourly on $n$ vertices, then (from
(\ref{eq4})) $g_2={n-4 \choose 2}$. By Proposition \ref{P20}\,(c),
$\mu_1(M) = g_2(M)/10 = (n-4)(n-5)/20$. Therefore, for $M$ to be
$\FF$-tight, Proposition \ref{P18} requires $\beta_1(M) = (n-4)(n-
5)/20$. Thus, we have the ``only if\," part of (a) and (b).

Now, we prove the ``if\," parts. If $d=1$, the result is trivial
since $S^{\hspace{.1mm}1}_3$ is the only 2-neighbourly closed
1-manifold. If $d=2$, the result is immediate from Lemma
\ref{L10}. If $d=3$, we have $\mu_0 = 1 =\beta_0$ because of
connectedness, and hence $\mu_3 = 1 = \beta_3$ by duality. Also,
$\mu_1 = g_2/{d+2 \choose 2} = \beta_1$ by Proposition
\ref{P20}\,(c) and hypothesis. Therefore, $\mu_2 = g_2/{d+2
\choose 2} = \beta_2$ by duality. Hence, by Proposition \ref{P18},
$M$ is $\FF$-tight. So, assume that $d\geq 4$. Then $\mu_0 = 1
=\beta_0$ and hence $\mu_d = 1 = \beta_d$. By Propositions
\ref{P20}\,(c) and \ref{P21}\,(c), $\mu_1 = g_2/{d+2 \choose 2} =
\beta_1$, and hence $\mu_{d-1} = g_2/{d +2 \choose 2} =
\beta_{d-1}$ by duality. Also, by Propositions \ref{P20}\,(a) and
\ref{P21}\,(d), $\mu_i = 0 = \beta_i$ for $2\leq i \leq d-2$.
Thus, $\mu_i = \beta_i$ for all $i$. Hence $M$ is $\FF$-tight by
Proposition \ref{P18}. \hfill $\Box$

\bigskip

The ``$d=4$\," case of the following result is due to Walkup
\cite{wa} and K\"{u}hnel \cite{k95} (cf.
\cite[Proposition 2]{bd10}). Part (b) of this proposition is
Theorem 5 of Lutz, Sulanke and Swartz \cite{lss}.

\begin{prop}[A lower bound theorem for triangulated
manifolds]\hspace{-1.5mm}{\bf .}  \label{P23}
Let $M$ be a connected closed triangulated manifold of dimension
$d\geq 3$. Let $\beta_1 = \beta_1(M; \ZZ_2)$. Then the face vector
of $M$ satisfies\,:
\begin{enumerate}
\item[{\rm (a)}] $f_j \geq \left\{\begin{array}{ll}
{\,d+1\, \choose j}f_0 +j{\,d+2\, \choose j+1}(\beta_1 -1), & \mbox{if } ~ 1 \leq j < d,  \\[2mm]
df_0 + (d-1)(d+2)(\beta_1-1), & \mbox{if } ~ j = d.
\end{array}
\right.
$
\item[{\rm (b)}] ${\,f_0-d-1\, \choose 2} \geq {\,d+2\, \choose 2}\beta_1$.
\end{enumerate}
When $d\geq 4$, equality holds in {\rm (a)} $($for some $j \geq 1$,
equivalently, for all $j\,)$ if and only if $M\in {\cal W}_1(d)$,
and equality holds in {\rm (b)} if and only if $M$ is a $2$-neighbourly
member of ${\cal W}_1(d)$.
\end{prop}

\noindent {\bf Proof.} Using Corollary \ref{C8} and formulae
(\ref{eq5}) and (\ref{eq6}), it is easy to see that the face
vector of any 1-stellated $(d-1)$-sphere $S$ is given by
$$
f_{j-1}(S) = \left\{\begin{array}{ll}
{d \choose j-1\,}f_0(S) -(j-1){\,d+1\, \choose j}, & \mbox{if }
~ 1 \leq j < d,  \\[2mm]
(d-1)f_0(S) - (d-2)(d+1), & \mbox{if } ~ j = d.
\end{array}
\right.
$$
The lower bound theorem of Kalai (\cite{ka}) says that if $L$ is 
any triangulated closed manifold of
dimension $d-1$ with $f_0(L) = f_0(S)$, then $f_{j-1}(L) \geq
f_{j-1}(S)$ for $1\leq j\leq d$. Also, when $d\geq 4$, equality
holds here (for some $j>1$, equivalently, for all $j$) if and only if
$L$ is an 1-stellated sphere. Applying this result to the vertex links
$L_x$, $x \in V(M)$, of $M$, we get
$$
f_{j-1}(L_x) \geq \left\{\begin{array}{ll}
{d \choose \,j-1\,}f_0(L_x) -(j-1){\,d+1\, \choose j}, & \mbox{if } ~ 1 \leq j < d,  \\[2mm]
(d-1)f_0(L_x) - (d-2)(d+1), & \mbox{if } ~ j = d.
\end{array}
\right.
$$
But $\sum_{x\in V(M)}f_{j-1}(L_x) = (j+1)f_j$, $1\leq j \leq d$.
Therefore, adding the above inequalities over all vertices $x$ of $M$,
we get that the face vector of $M$ satisfies\,:

\begin{eqnarray} \label{eq13}
f_{j} \geq \left\{\begin{array}{ll}
\frac{2}{j+1}{d \choose \,j-1\,}f_1 -\frac{j-1}{j+1}{\,d+1\,
\choose j}f_0, & \mbox{if } ~ 1 \leq j < d,  \\[2mm]
\frac{2d-2}{d+1}f_1 - (d-2)f_0, & \mbox{if } ~ j = d.
\end{array}
\right.
\end{eqnarray}
Also, when $d\geq 4$, equality holds in (\ref{eq13}) for some $j> 1$
(equivalently, for all $j$) if and only if all the vertex links
of $M$ are 1-stellated, i.e., if and only if $M\in {\cal W}_1(d)$.

Now, Theorem 5.2 of Novik and Swartz \cite{ns} says that $g_2(M) \geq
{d+2 \choose 2}\beta_1$, i.e.,
\begin{eqnarray} \label{eq14}
f_1 \geq (d+1)f_0 + {d+2 \choose 2}(\beta_1 -1).
\end{eqnarray}
Also, when $d\geq 4$, equality holds in (\ref{eq14}) if and only if
$M\in {\cal W}_1(d)$. Notice that (\ref{eq14}) is just the case
$j=1$ of part (a). Now, combining (\ref{eq13}) and (\ref{eq14}), we get
all cases of part (a), after a little simplification.

We also have
$
(d+1)f_0 + {d+2\, \choose 2}(\beta_1-1) \leq f_1 \leq \frac{f_0(f_0-1)}{2}\,,
$
where the first inequality is from (\ref{eq14}) (with equality for
$d\geq 4$ if and only if $M \in {\cal W}_1(d)$) and the second inequality is
trivial (with equality if and only if $M$ is 2-neighbourly).
Hence we have $(d+1)f_0 + {d+2\, \choose 2}(\beta_1-1) \leq \frac{f_0(f_0-1)}{2}$,
which simplifies to the inequality in part (b). \hfill $\Box$

\bigskip

Next we introduce\,:

\begin{defn}\hspace{-1.8mm}{\bf .} \label{strongly_minimal}
{\rm A $d$-dimensional simplicial complex $X$ is said to be {\em
minimal} if $f_0(X) \leq f_0(Y)$ for every triangulation $Y$ of
the geometric carrier $|X|$ of $X$. We shall say that $X$ is {\em
strongly minimal} if $f_i(X) \leq f_i(Y)$, $0\leq i \leq d$,
for all such $Y$.}
\end{defn}

In \cite{kl}, K\"{u}hnel and Lutz conjectured that all $\FF$-tight
triangulated  manifolds are strongly minimal. Our next result is
a powerful evidence in favour of this conjecture. It also generalizes
a result of Swartz \cite[Theorem 4.7]{s}, who proved that the
tight triangulations $K^{\,d}_{2d+3}\in {\cal W}_1(d)$ (cf. Example
\ref{E3}\,(b) below) are strongly minimal for all $d\geq 2$.

\begin{cor}\hspace{-1.8mm}{\bf .} \label{C11}
Every $\FF$-tight member of\,  ${\cal W}_1(d)$ is strongly minimal.
\end{cor}

\noindent {\bf Proof.} Let $M_0 \in {\cal W}_1(d)$ be $\FF$-tight.
Proposition \ref{P22} implies that $M_0$ is $\ZZ_2$-tight. By the same
proposition, $M_0$ is 2-neighbourly. Now, for $d\leq 2$, any
2-neighbourly closed $d$-manifold is strongly minimal. (This is
entirely trivial for $d=1$ since $S^1_3$ is the only 2-neighbourly
closed manifold in this case. It is only a little less trivial for
$d=2$\,: the face vector is determined by $f_0$ and $\beta_1$ in
this case.) So, assume $d\geq 3$.

We claim that $M_0$ attains all the bounds in Proposition \ref{P23}.
This is immediate from the proposition itself when $d\geq 4$. If $d=3$,
then - as $M_0\in {\cal W}_1(3)$ is $\ZZ_2$-tight, we have $g_2(M_0)
= 10\beta_1$ by Proposition \ref{P22}. Hence, following the proof
of Proposition \ref{P23}, one sees that $M_0$ attains the bounds in
Proposition \ref{P23} in this case also.

Now, let $M$ be any triangulation of $|M_0|$. Since the Betti numbers
are topological invariants, we have $\beta_1(M; \ZZ_2) = \beta_1(M_0;
\ZZ_2) = \beta_1$ (say). By Proposition \ref{P23} and the above claim,
$$
{f_0(M) -d-1 \choose 2} \geq {d+2 \choose 2}\beta_1 = {f_0(M_0) -d-1
\choose 2}.
$$
Since, trivially, $f_0(M)$, $f_0(M_0) \geq d+2$, this implies $f_0(M)
\geq f_0(M_0)$. Therefore, we get\,:
$
f_j(M) \geq a_jf_0(M) + b_j \geq a_jf_0(M_0) + b_j = f_j(M_0) ~ \mbox{
for } ~ 0\leq j\leq d,
$
where $a_j >0$ and $b_j$ are constants (depending only on $d$, $j$ and
$\beta_1$) given by Proposition \ref{P23}. \hfill $\Box$

\bigskip

2-neighbourly members of ${\cal W}_1(d)$ were called ``{\em tight
neighbourly}\," by Lutz, Sulanke and Swartz \cite{lss}. By
Proposition \ref{P22}, tight neighbourly manifolds of dimension
$\neq 3$ are $\ZZ_2$-tight if non-orientable, and tight if
orientable. Part (b) of Proposition \ref{P22} gives a criterion
for the tightness of a tight neighbourly 3-manifold in terms of
its first Betti number. Corollary \ref{C11} shows that tight
neighbourly triangulations of dimension $\geq 4$ are strongly minimal.

In consequence of Corollary \ref{C9}, a member of ${\cal W}_k(d)$
can be at most $(k+1)$-neighbourly, unless it is a standard
sphere. This observation leads us to introduce\,:

\bigskip

\noindent {\bf Notation\,:} ${\cal W}_k^{\ast}(d)$ will denote the
subclass of ${\cal W}_k(d)$ consisting of all the
$(k+1)$-neighbourly members of the latter class.

\medskip

Thus, Proposition \ref{P22} was about the class ${\cal
W}_1^{\ast}(d)$, the class of tight neighbourly $d$-manifolds.
By a result of Kalai (cf. \cite[Corollary 8.4]{ka} or
\cite[Proposition 3]{bd10}), for $d\geq 4$, any member of ${\cal
W}_1(d)$ triangulates the connected sum of finitely many copies of
$S^{1}\times S^{\hspace{.2mm}d-1}$ or of $\TPSSD$. The following
result may be compared with Kalai's result.

\begin{prop}\hspace{-1.8mm}{\bf .}  \label{P24}
Let $M \in {\cal W}_k^{\ast}(d)$, where $k\geq 2$ and $d\geq
2k+2$. Suppose $M$ is not a $\ZZ$-homology sphere. Then $M$ has
the same integral homology group as the connected sum of $\beta$
copies of $S^{\hspace{.2mm}k}\times S^{\hspace{.2mm}d-k}$, where
the positive integer $\beta$ is given in terms of the number $m$
of vertices of $M$ by the formula $\beta = {m+k-d-2 \choose
k+1}/{d+2 \choose k+1}$. In consequence, we must have $m \geq
2d+4-k$ and ${d+2 \choose k+1}$ divides ${m+k-d-2 \choose k+1}$.
\end{prop}
(As to the inequality $m\geq 2d+4-k$, note that by a result of
Brehm and K\"{u}hnel \cite{bk}, this lower bound on the number of
vertices holds, more generally, for any triangulation of a closed
$d$-manifold which is not $k$-connected.)

\bigskip

\noindent {\bf Proof.} Since $M$ is at least 3-neighbourly, it is
simply connected and hence orientable. Therefore, by Poincar\'e
duality, the Betti numbers of $M$ with respect to any field $\FF$
satisfy $\beta_{d-i} = \beta_i$, $0\leq i \leq d$. Since $M$ is
connected, we have $\beta_0 =1$ and hence $\beta_d =1$. By Lemma
\ref{L9}, $\beta_i =0$ for $1\leq i \leq k-1$. Hence, by duality,
$\beta_i=0$ for $d-k+1 \leq i \leq d-1$. By Proposition \ref{P21}
and duality, $\beta_{d-k} = \beta_k = g_{k+1}/{d+2 \choose k+1}$
and $\beta_i =0$ for $k+1 \leq i \leq d-k-1$. Thus, all the Betti
numbers of $M$ are independent of the choice of the field $\FF$.
Therefore, by the universal coefficients theorem, the
$\ZZ$-homologies of $M$ are torsion free, and the $\ZZ$-Betti
numbers are given by the same formulae as above. Since all the
middle Betti numbers except possibly $\beta_k = \beta_{d-k}$ are
zero, it follows that if $M$ is not a homology sphere, then the
value $\beta = g_{k+1}/{d+2 \choose k+1} = {m+k-d-2 \choose
k+1}/{d+2 \choose k+1}$ of $\beta_k$ must be a strictly positive
integer. The last statement follows from this observation. \hfill
$\Box$

\begin{remark}\hspace{-1.8mm}{\bf .} \label{R4}
{\rm If $M \in {\cal W}_k^{\ast}(d)$ is a homology sphere $(k\geq
2, \, d\geq 2k+2)$, then - by the above proof - we must have
$g_{k+1}(M) = 0$. In this case, we expect (GLBC) $M$ to be a
$k$-stacked $d$-sphere. Since $S^{\hspace{.2mm}d}_{d+2}$ is the
only $(k+1)$-neighbourly $k$-stacked $d$-sphere, we should
therefore have $M= S^{\hspace{.2mm}d}_{d+2}$ in this case. Thus,
it should be possible to replace the hypothesis ``$M$ is not an
integral homology sphere" in Proposition \ref{P24} by the simpler
hypothesis $M \neq S^{\hspace{.2mm}d}_{d+2}$. }
\end{remark}

The next result is our generalization of Proposition \ref{P22} to
the case $k\geq 2$.

\begin{prop}[A combinatorial criterion for tightness]\hspace{-1.5mm}{\bf .}
\label{P25} Let $M \in {\cal W}_k^{\ast}(d)$, where $k\geq 2$.
Then we have\,: \vspace{-2mm}
\begin{enumerate}
\item[{\rm (a)}]  if $d \neq 2k+1$ then $M$ is tight, and
\vspace{-2mm}
\item[{\rm (b)}] if $d = 2k+1$, then $M$ is $\FF$-tight if and
only if \, $\beta_k(M; \FF) = \frac{{n-k-3 \choose k+1}}{{2k+3
\choose k+1}}$, where $n = f_0(M)$.
\end{enumerate}
\end{prop}

\noindent {\bf Proof.} This is trivial if $M =
S^{\hspace{.2mm}d}_{d +2}$. So, assume $M \neq
S^{\hspace{.2mm}d}_{d+2}$. Then there is a vertex $x$ of $M$ such
that the link of $x$ in $M$ (is a $k$-neighbourly $(d-1)$-sphere
which) is not a standard sphere. Therefore, by Remark
\ref{R3}\,(c), $d-1 \geq 2k-1$. Thus $d\geq 2k$. If $d=2k$, then
the result follows from Lemma \ref{L10}. So, assume that $d\geq
2k+1$.

Notice that $M$ is at least 3-neighbourly, hence orientable. Thus
the duality result of Proposition \ref{P16}\,(e) applies. Since
$M$ is connected, we have $\beta_0 = 1 = \mu_0$ and hence $\beta_d
= 1 = \mu_d$. By Lemma \ref{L9}, $\beta_i = 0 = \mu_i$ for $1\leq
i \leq k-1$, hence also for $d-k+1\leq i \leq d-1$. Since $M$ is
$(k+1)$-neighbourly, we have $g_{k+1}(X) = {n+k-d-2 \choose k+1}$, where $n=f_0(M)$.
We also have $\beta_k = g_{k+1}/{d+2 \choose k+1} = \mu_k$ (by
hypothesis and Proposition \ref{P20}\,(c) when $d = 2k+1$; by
Propositions \ref{P20}\,(c) and \ref{P21}\,(c) when $d\geq 2k+2$).
Hence, $\beta_{d-k} = g_{k+1}/{d+2 \choose k+1} = \mu_{d-k}$. 
By Propositions \ref{P20}\,(a) and \ref{P21}\,(d), we
also have $\beta_i = 0 = \mu_i$ for $k+1 \leq i \leq d-k-1$. Thus,
$\beta_i = \mu_i$ for all $i$. Hence $M$ is $\FF$-tight by
Proposition \ref{P18} (and, when $d\neq 2k+1$, this argument
applies to all fields $\FF$).

For the converse statement in part (b), note that - more generally
- for any $M \in {\cal W}_k^{\ast}(d)$ with $d \geq 2k+1$, Lemma
\ref{L9} and Proposition \ref{P20}\,(c) imply that $\mu_k =
g_{k+1}/{d+2 \choose k+1}$. Therefore, for $M$ to be $\FF$-tight,
we must have (by Proposition \ref{P18}) $\beta_k = g_{k+1}/{d+2
\choose k+1}$ as well. Thus, for $d=2k+1$, $\beta_k = {n+k-d-2
\choose k+1}/{d+2 \choose k+1} = {n-k-3 \choose k+1}/{2k+3 \choose
k+1}$. \hfill $\Box$

\bigskip

Propositions \ref{P20} and \ref{P21} may be used to prove many
more tightness criteria. For instance, we have\,:

\begin{prop}[Another tightness criterion?]\hspace{-1.5mm}{\bf .}
\label{P26} Let $k\geq 2$, and let $M \in {\cal W}_k(d)$ be
$k$-neighbourly of dimension $d=2k$ or $d\geq 2k+2$, on $n$
vertices. If $M$ is $\FF$-orientable and $\beta_{k-1}(M; \FF) =
{n+k-d-3 \choose k}/{d+2 \choose k}$ then $M$ is $\FF$-tight.
\end{prop}
(Note that the requirement of $\FF$-orientability is automatically
fulfilled if $\FF = \ZZ_2$ or $k\geq 3$.)

\bigskip

\noindent {\bf Proof.} Since $M$ is $k$-neighbourly and $k\geq 2$,
we have $\mu_0 = 1 = \beta_0$ and $\mu_i = 0 = \beta_i$ for $1\leq
i \leq k-2$. Hence, by duality, $\mu_d = 1 = \beta_d$ and
$\mu_i=0=\beta_i$ for $d-k+2 \leq i \leq d-1$. Also, as $M$ is a
$k$-neighbourly $d$-manifold on $n$ vertices, it has $g_k =
{n+k-d-3 \choose k}$. Thus we have $g_k/{d+2 \choose k} =
\beta_{k-1} \leq \mu_{k-1} \leq g_k/{d+2 \choose k}$, where the
equality is by hypothesis and the inequalities are from
Propositions \ref{P16}\,(b) and \ref{P20}\,(b). Thus, $\mu_{d-k+1}
= \mu_{k-1} = \beta_{k-1} = \beta_{d-k+1}$. So, we have $\mu_i =
\beta_i$ for $0\leq i \leq k-1$ and for $d-k+1 \leq i \leq d$. So,
in view of Proposition \ref{P18}, to complete the proof it is
sufficient to show that $\mu_i = \beta_i$ for $k\leq i\leq d-k$ as
well.

If $d=2k$, then we have $\mu_i = \beta_i$ for all $i\neq k$.
Hence, by the equality statement in Proposition \ref{P16}\,(a), we
have $\mu_i = \beta_i$ for $i=k$ also.

Now, suppose $d\geq 2k+2$. Then, by Propositions \ref{P20}\,(c)
and \ref{P21}\,(c), $\mu_k - \mu_{k-1} = g_{k+1}/{d+2 \choose k+1}
= \beta_k-\beta_{k-1}$. Hence $\mu_k = \mu_{k-1} + g_{k+1}/{d+2
\choose k+1} = g_{k}/{d+2 \choose k} + g_{k+1}/{d+2 \choose k+1} =
\beta_{k-1} + g_{k+1}/{d+2 \choose k+1} =  \beta_{k}$. Therefore,
by duality, $\mu_{d-k} = \beta_{d-k}$. Also, by Propositions
\ref{P20}\,(a) and \ref{P21}\,(d), $\mu_i = 0 =\beta_i$ for
$k+1\leq i \leq d-k-1$. \hfill $\Box$

\begin{remark}\hspace{-1.8mm}{\bf .} \label{R5}
{\rm Note that ${\cal W}_{k-1}(d) \subseteq {\cal W}_k(d)$. If
$M\in {\cal W}_{k-1}(d)$ is $k$-neighbourly and $d\geq 2k$, then
$M$ is $\FF$-tight (whenever it is $\FF$-orientable) by
Propositions \ref{P22} and \ref{P25}. Also, in this case, we
automatically have (by Lemma \ref{L9} and Proposition
\ref{P21}\,(c)) $\beta_{k-1} = {n+k-d-3 \choose k}/{d+2 \choose
k}$. Thus, Proposition \ref{P26} has new content only for $M \in
{\cal W}_{k}(d) \setminus {\cal W}_{k-1}(d)$. }
\end{remark}

\section{Examples, counterexamples, questions and conjectures}

\begin{eg}[Stellated versus stacked spheres]\hspace{-1.5mm}{\bf .}
\label{E1}
\begin{itemize}
\item[{\bf (a)}] {\rm Let $S^{\hspace{.2mm}d}_{2d+2} =
(S^{\hspace{.2mm}0}_{2})^{\ast\,d+1}$, the join of $d+1$ copies of
$S^{\hspace{.2mm}0}_{2}$. Being the boundary complex of the
$(d+1)$-dimensional cross polytope, $S^{\hspace{.2mm}d}_{2d+2}$ is
a polytopal $d$-sphere. Therefore, by Proposition \ref{P12} (with
$k=0$), it is $d$-stellated. Also, by Proposition \ref{P1}, it is
$d$-stacked. Since $S^{\hspace{.2mm}d}_{2d+2}$ is the clique
complex of its edge graph (1-skeleton), it is not $(d-
1)$-stacked. (If there was a $(d-1)$-stacked $(d+1)$-ball $B$ such
that $\partial B = S^{\hspace{.2mm}d}_{2d+2}$, then all the faces
of $B$ would be cliques of the edge graph of
$S^{\hspace{.2mm}d}_{2d +2}$. But, all such cliques are in
$S^{\hspace{.2mm}d}_{2d+2}$ itself.) For the same reason,
$S^{\hspace{.2mm}d}_{2d+2}$ does not contain any induced standard
sphere except $S^{\hspace{.2mm}0}_{2}$. Therefore, it does not
admit any bistellar move of index $\geq 2$. Therefore,
$S^{\hspace{.2mm}d}_{2d+2}$ is not $(d-1)$-stellated.
(By the comment following Corollary \ref{C2}, any $k$-stellated
$d$-sphere, excepting $S^{\hspace{.2mm}d}_{d +2}$, admits a
bistellar move of index $>d-k$.)

It is easy see that all the induced subcomplexes of
$S^{\hspace{.2mm}d}_{2d +2}$ are spheres and balls. Also, for
$-1 \leq i\leq d$, $S^{\hspace{.2mm}d}_{2d +2}$ has exactly ${d+1
\choose i+1}$ induced $i$-spheres, each isomorphic to
$S^{\hspace{.2mm}i}_{2i +2}$. Hence one computes\,:
$$
\sigma_i(S^{\hspace{.2mm}d}_{2d +2}) = \left\{\begin{array}{rl}
-\frac{2d}{2d+1}, & i=0 \\[2mm]
\frac{{d+1 \choose i+1}}{{2d+2 \choose 2i+2}}, & 1\leq i\leq d.
\end{array}
\right.
$$
Hence one finds (since all the vertex links of
$S^{\hspace{.2mm}d}_{2d +2}$ are isomorphic to
$S^{\hspace{.2mm}d-1}_{2d}$) that the mu-vector of
$S^{\hspace{.2mm}d}_{2d +2}$ is given by
\vspace{-2mm}
$$
\mu_i(S^{\hspace{.2mm}d}_{2d +2}) = \frac{{d \choose i}}{{2d
\choose 2i}}  ~~ \mbox{for} ~~ 0\leq i\leq d.
$$
Surprisingly the mu-vector of $S^{\hspace{.2mm}d}_{2d +2}$
satisfies the duality relation $\mu_{d-i} \equiv \mu_i$, even
though it is not 2-neighbourly. Also, as a curiosity, we find
${\displaystyle \sum_{i=0}^d}
(-1)^{d-i}\mu_i(S^{\hspace{.2mm}d}_{2d +2}) =
\frac{2d+1}{2d+2}\chi(S^{\hspace{.2mm}d})$. Thus,
$S^{\hspace{.2mm}d}_{2d +2}$ fails the strong Morse inequalities
(Proposition \ref{P16}\,(a)).

}

\item[{\bf (b)}] {\rm It is more difficult to find examples of
$(d+1)$-stellated $d$-spheres (i.e., combinatorial $d$-spheres)
which are not $d$-stellated. The following example is due to
Dougherty, Faber and Murphy \cite{dfm}.

Let $S^{\,3}_{16}$ be the pure 3-dimensional simplicial complex
with vertex set $\ZZ_{16} = \ZZ/16\ZZ$ and an automorphism $i
\mapsto i+1$ (mod 16). Modulo this automorphism, the basic facets
of $S^{\,3}_{16}$ are\,:
$$
\{0,1,4,6\}, \{0,1,4,9\}, \{0,1,6,14\}, \{0,1,8,9\}, \{0,1,8,10\},
\{0,1,10,14\}, \{0,2,9,13\}.
$$
Of these, the fourth facet generates an orbit of length 8, while
each of the other facets generates an orbit of length 16. Thus,
$S^{\,3}_{16}$ has $1 \times 8 + 6 \times 16 = 104$ facets. The
face vector of $S^{\,3}_{16}$ is $(16, 120, 208, 104)$. Since $120
= {16 \choose 2}$, $S^{\,3}_{16}$ is 2-neighbourly and hence it
does not allow any bistellar 1-move. Also, it is easy to verify
that $S^{\,3}_{16}$ has no edge of (minimum) degree 3, so that it
does not allow any bistellar move of index 2 or 3 either. (So,
$S^{\,3}_{16}$ is an {\em unflippable} $3$-sphere in the sense of
\cite{dfm}\,: it does not allow any bistellar move of positive
index.) Thus, $S^{\,3}_{16}$ is not 3-stellated. (Being a
combinatorial 3-sphere, it is of course 4-stellated.) Following
the proof of Proposition \ref{P1}, fix a vertex $x$ of
$S^{\,3}_{16}$, and let $B^{\,4}_{16} = \{\{x\} \sqcup\alpha \, :
\, x \not\in \alpha\in S^{\,3}_{16}\}$. Then $ B^{\,4}_{16}$ is a
4-ball with $\partial B^{\,4}_{16} = S^{\,3}_{16}$. Since
$S^{\,3}_{16}$ is 2-neighbourly, $ B^{\,4}_{16}$ is a 2-stacked
ball, and hence $S^{\,3}_{16}$ is an example of a 2-stacked
3-sphere which is not even 3-stellated. If $ B^{\,4}_{16}$ was
shellable, then (by Proposition \ref{P2}) it would be 2-shelled
and hence (by Corollary \ref{C1}) $S^{\,3}_{16}$ would be
2-stellated. Thus, $ B^{\,4}_{16}$ is an example of a
non-shellable 2-stacked ball. }

\item[{\bf (c)}] {\rm It is even more difficult to find examples
of triangulated $d$-spheres which are not $(d+1)$-stellated (i.e.,
not combinatorial $d$-spheres). Trivially, all triangulated
spheres of dimension $d\leq 3$ are combinatorial spheres. In
\cite{ed} and \cite{f}, Edwards and Freedman proved that a
triangulated homology manifold of dimension $d\geq 3$ is a
triangulated manifold if and only if all its vertex links are
simply connected. In conjunction with Perelman's theorem
(3-dimensional Poincar\'{e} conjecture) this shows that all
triangulated 4-manifolds are combinatorial manifolds. The (non-)
existence of triangulated 4-spheres which are not combinatorial
spheres is equivalent to the still unresolved 4-dimensional smooth
Poincar\'{e} conjecture. (According to \cite{bd5}, any such
4-sphere would require at least 13 vertices.) Thus, $d=5$ is the
smallest dimension in which we may reasonably expect triangulated
spheres which are not combinatorial spheres. The following
16-vertex triangulation $\Sigma^{\,3}_{16}$ of the Poincar\'{e}
(integral) homology 3-sphere was found by Bj\"{o}rner and Lutz
\cite{bl}. The vertices of $\Sigma^{\,3}_{16}$ are $1, \dots, 9,
1^{\hspace{.1mm}\prime}, \dots, 7^{\hspace{.2mm}\prime}$. Its
facets are\,: $1  2  4  9$,    $1  2  4  6^{\hspace{.2mm}\prime}
$, $1265^{\hspace{.2mm}\prime} $, $1266^{\hspace{.2mm}\prime}$,
$1295^{\hspace{.2mm}\prime} $, $1343^{\hspace{.2mm}\prime}$,
$1346^{\hspace{.2mm}\prime}$, $1371^{\hspace{.1mm}\prime} $,
$1373^{\hspace{.2mm}\prime} $, $ 1 3 1^{\hspace{.1mm}\prime}
6^{\hspace{.2mm}\prime} $, $1  4 9 3^{\hspace{.2mm}\prime} $,
$1564^{\hspace{.2mm}\prime} $, $ 15  6 5^{\hspace{.2mm}\prime} $,
$158 2^{\hspace{.2mm}\prime} $, $ 1584^{\hspace{.2mm}\prime} $,
$152^{\hspace{.2mm}\prime}5^{\hspace{.2mm}\prime}$,
$164^{\hspace{.2mm}\prime}6^{\hspace{.2mm}\prime}$,
$1781^{\hspace{.1mm}\prime}$, $1782^{\hspace{.2mm}\prime}$,
$172^{\hspace{.2mm}\prime}3^{\hspace{.2mm}\prime}$,
$181^{\hspace{.1mm}\prime}4^{\hspace{.2mm}\prime}$,
$192^{\hspace{.2mm}\prime}3^{\hspace{.2mm}\prime}$,
$192^{\hspace{.2mm}\prime}5^{\hspace{.2mm}\prime}$,
$11^{\hspace{.1mm}\prime}4^{\hspace{.2mm}\prime}6^{\hspace{.2mm}\prime}$,
$2351^{\hspace{.1mm}\prime}$, $2352^{\hspace{.2mm}\prime}$,
$2371^{\hspace{.1mm}\prime}$, $2374^{\hspace{.2mm}\prime}$,
$232^{\hspace{.2mm}\prime}4^{\hspace{.2mm}\prime}$,
$2494^{\hspace{.2mm}\prime}$,
$242^{\hspace{.2mm}\prime}4^{\hspace{.2mm}\prime}$,
$242^{\hspace{.2mm}\prime}6^{\hspace{.2mm}\prime}$,
$2582^{\hspace{.2mm}\prime} $, $ 2 5 8  3^{\hspace{.2mm}\prime} $,
$ 2 5 1^{\hspace{.1mm}\prime} 3^{\hspace{.2mm}\prime} $,
$261^{\hspace{.1mm}\prime} 3^{\hspace{.2mm}\prime} $,
$261^{\hspace{.1mm}\prime} 5^{\hspace{.2mm}\prime} $,
$263^{\hspace{.2mm}\prime} 6^{\hspace{.2mm}\prime} $,
$2794^{\hspace{.2mm}\prime}$, $ 2 7  9  5^{\hspace{.2mm}\prime}$,
$ 2 7 1^{\hspace{.1mm}\prime} 5^{\hspace{.2mm}\prime} $,
$282^{\hspace{.2mm}\prime} 6^{\hspace{.2mm}\prime} $,
$283^{\hspace{.2mm}\prime} 6^{\hspace{.2mm}\prime} $,
$3455^{\hspace{.2mm}\prime}$, $ 3 4  5  6^{\hspace{.2mm}\prime} $,
$ 3  4 3^{\hspace{.2mm}\prime} 5^{\hspace{.2mm}\prime} $,
$351^{\hspace{.1mm}\prime} 6^{\hspace{.2mm}\prime} $,
$352^{\hspace{.2mm}\prime} 5^{\hspace{.2mm}\prime} $,
$373^{\hspace{.2mm}\prime} 4^{\hspace{.2mm}\prime} $,
$32^{\hspace{.2mm}\prime}4^{\hspace{.2mm}\prime}5^{\hspace{.2mm}\prime}$,
$33^{\hspace{.2mm}\prime}4^{\hspace{.2mm}\prime}5^{\hspace{.2mm}\prime}$,
$    4  5  6  7 $, $ 4  5  6 5^{\hspace{.2mm}\prime} $,
$4576^{\hspace{.2mm}\prime}$, $4 6  7  2^{\hspace{.2mm}\prime} $,
$ 4 6 1^{\hspace{.1mm}\prime} 2^{\hspace{.2mm}\prime} $,
$461^{\hspace{.1mm}\prime} 5^{\hspace{.2mm}\prime} $,
$472^{\hspace{.2mm}\prime} 6^{\hspace{.2mm}\prime} $,
$4893^{\hspace{.2mm}\prime}$, $4894^{\hspace{.2mm}\prime}$,
$481^{\hspace{.1mm}\prime} 4^{\hspace{.2mm}\prime} $,
$481^{\hspace{.1mm}\prime} 5^{\hspace{.2mm}\prime} $,
$483^{\hspace{.2mm}\prime} 5^{\hspace{.2mm}\prime} $,
$41^{\hspace{.1mm}\prime}2^{\hspace{.2mm}\prime}4^{\hspace{.2mm}\prime}$,
$5674^{\hspace{.2mm}\prime} $, $5794^{\hspace{.2mm}\prime}$,
$5796^{\hspace{.2mm}\prime} $, $5893^{\hspace{.2mm}\prime}$,
$5894^{\hspace{.2mm}\prime}$,
$591^{\hspace{.1mm}\prime}3^{\hspace{.2mm}\prime}$,
$591^{\hspace{.1mm}\prime} 6^{\hspace{.2mm}\prime}$,
$672^{\hspace{.2mm}\prime} 3^{\hspace{.2mm}\prime}$,
$673^{\hspace{.2mm}\prime} 4^{\hspace{.2mm}\prime}$,
$61^{\hspace{.1mm}\prime}2^{\hspace{.2mm}\prime}3^{\hspace{.2mm}\prime}$,
$63^{\hspace{.2mm}\prime}4^{\hspace{.2mm}\prime}6^{\hspace{.2mm}\prime}$,
$781^{\hspace{.1mm}\prime}5^{\hspace{.2mm}\prime}$,
$782^{\hspace{.2mm}\prime} 6^{\hspace{.2mm}\prime}$,
$785^{\hspace{.2mm}\prime} 6^{\hspace{.2mm}\prime}$,
$795^{\hspace{.2mm}\prime} 6^{\hspace{.2mm}\prime}$,
$83^{\hspace{.2mm}\prime}5^{\hspace{.2mm}\prime}6^{\hspace{.2mm}\prime}$,
$91^{\hspace{.1mm}\prime}2^{\hspace{.2mm}\prime}3^{\hspace{.2mm}\prime}$,
$91^{\hspace{.1mm}\prime}2^{\hspace{.2mm}\prime}7^{\hspace{.2mm}\prime}$,
$91^{\hspace{.1mm}\prime}6^{\hspace{.2mm}\prime}7^{\hspace{.2mm}\prime}$,
$92^{\hspace{.2mm}\prime}5^{\hspace{.2mm}\prime}7^{\hspace{.2mm}\prime}$,
$95^{\hspace{.2mm}\prime}6^{\hspace{.2mm}\prime}7^{\hspace{.2mm}\prime}$,
$1^{\hspace{.1mm}\prime}2^{\hspace{.2mm}\prime}
4^{\hspace{.2mm}\prime}7^{\hspace{.2mm}\prime}$,
$1^{\hspace{.1mm}\prime}4^{\hspace{.2mm}\prime}
6^{\hspace{.2mm}\prime}7^{\hspace{.2mm}\prime}$,
$2^{\hspace{.2mm}\prime}4^{\hspace{.2mm}\prime}
5^{\hspace{.2mm}\prime}7^{\hspace{.2mm}\prime} $,
$3^{\hspace{.2mm}\prime}4^{\hspace{.2mm}\prime}
5^{\hspace{.2mm}\prime}6^{\hspace{.2mm}\prime} $,
$4^{\hspace{.2mm}\prime}5^{\hspace{.2mm}\prime}
6^{\hspace{.2mm}\prime}7^{\hspace{.2mm}\prime}$. The face vector
of $\Sigma^{\,3}_{16}$ is $(16,106, 180, 90)$. Bj\"{o}rner and
Lutz conjectured that it is strongly minimal.

Note that the vertex $6^{\hspace{.2mm}\prime}$ is adjacent with
all other vertices in $\Sigma^{\,3}_{16}$. Let $D^{\,4}_{16}$ be
the 4-dimensional simplicial complex whose facets are $\alpha \cup
\{6^{\hspace{.2mm}\prime}\}$, where $\alpha$ ranges over all
facets of $\Sigma^{\,3}_{16}$ not containing the vertex
$6^{\hspace{.2mm} \prime}$. Define $S^{\,5}_{18}=
\partial(D^{\,4}_{16} \ast B^{1}_{2})$, the boundary of the join
of $D^{\,4}_{16}$ and an edge. Observe that, $|S^{\,5}_{18}|$ is
the double suspension of the Poincar\'{e} homology sphere
$|\Sigma^{\,3}_{16}|$. Therefore, by Cannon's double suspension
theorem (cf. \cite{ca}, actually Cannon's theorem is a
straightforward consequence of the result of Edwards and Freedman
quoted above), $S^{\,5}_{18}$ is a triangulated 5-sphere. Since it
has $\Sigma^{\,3}_{16}$ as the link of an edge, $S^{\,5}_{18}$ is
not a combinatorial sphere.

Let $D^{\,6}_{18} = D^{\,4}_{16} \ast B^{1}_{2}$, $D^{\,7}_{19} =
D^{\,4}_{16} \ast B^{\hspace{.1mm}2}_{3}$ and $S^{\,6}_{19} =
\partial D^{\,7}_{19}$. By the above logic, $S^{\,6}_{19}$ is a
triangulated 6-sphere. Since $D^{\,6}_{18}$ is the antistar of a
vertex in $S^{\,6}_{19}$, it follows from Lemma 4.1 in \cite{bd9}
that $D^{\,6}_{18}$ is a triangulated $6$-ball. Since the vertex
$6^{\hspace{.2mm}\prime}$ is adjacent to all the vertices in
$\Sigma^{\,3}_{16}$, the construction of $D^{\,6}_{18}$ shows that
all the 3-faces of $D^{\,6}_{18}$ lie in its boundary. Thus,
$D^{\,6}_{18}$ is a 2-stacked triangulated ball. As $S^{\,5}_{18}
= \partial D^{\,6}_{18}$, it follows that $S^{\,5}_{18}$ is an
example of 2-stacked 5-sphere which is not even 6-stellated. }

\item[{\bf (d)}] {\rm Let $S$ be a triangulated $d$-sphere and $B$
is a $k$-stacked ball such that $\partial B = S$. Then, for any
$e\geq 0$, $B \ast B^{\hspace{.2mm}e}_{e+1}$ is a $k$-stacked
ball, and hence $\partial(B \ast B^{\hspace{.2mm}e}_{e+1})$ is a
$k$-stacked $(d+e+1)$-sphere. Also, $S$ is a combinatorial sphere
if and only if $\partial(B \ast B^{\hspace{.2mm}e}_{e +1})$ is so.
Applying this construction to the pair $(S^{\hspace{.2mm}5}_{18},
D^{\hspace{.2mm}6}_{18})$ in example (c) above, we find that
for each $d \geq 5$, there are 2-stacked triangulated $d$-spheres
which are not even $(d+1)$-stellated.

\noindent {\bf Claim.} If $B= B^{\,4}_{16}$ is as in example (b)
above then $\partial(B \ast B^{\hspace{.1mm}e}_{e+1})$
is unflippable.

For $e \geq 0$, let $\widetilde{B}^{\,e+5}_{e+17} : = B^{\,4}_{16}
\ast B^{\,e}_{e+1}$ and $\widetilde{S}^{\,e+4}_{e+17} := \partial
\widetilde{B}^{\,e+5}_{e+17}$. Thus, $\widetilde{S}^{\,e+4}_{e+17}
= (S^{\,3}_{16} \ast B^{\,e}_{e+1}) \cup (B^{\,4}_{16} \ast
S^{\,e-1}_{e+1})$. Since $S^{\,3}_{16}$ is 2-neighbourly, so is
$\widetilde{S}^{\,e +4}_{e+17}$. Therefore, $\widetilde{S}^{\,e
+4}_{e+17}$ does not admit any bistellar 1-move. Suppose, if
possible, that $\alpha \mapsto \beta$ is a bistellar move of index
$\geq 2$ on $\widetilde{S}^{\,e +4}_{e + 17}$. Thus, ${\rm
lk}_{\widetilde{S}^{\,e +4}_{e+17}}(\alpha) = \partial \beta$ and
$\dim(\beta) \geq 2$, $\beta \not\in \widetilde{S}^{\,e +4}_{e
+17}$. Write $\alpha = \alpha_1\sqcup \alpha_2$, where $\alpha_1$
is a face of $B^{\,4}_{16}$ and $\alpha_2$ is a face of
$B^{\,e}_{e +1}$. If $\alpha_1$ is an interior face of
$B^{\,4}_{16}$, then $\alpha_2 \in S^{\,e- 1}_{e+1}$ and $\partial
\beta = {\rm lk}_{\widetilde{S}^{\,e +4}_{e+17}}(\alpha) =  {\rm
lk}_{B^{\,4}_{16}}(\alpha_1)\ast {\rm lk}_{S^{\,e-1}_{e+
1}}(\alpha_2)$. Since the standard sphere $\partial \beta$ can't
be written as the join of two spheres, it follows that either
$\alpha_1$ is a facet of $B^{\,4}_{16}$ or $\alpha_2$ is a facet
of $S^{\,e-1}_{e+1}$. If $\alpha_1$ is a facet of $B^{\,4}_{16}$,
then $\partial \beta = {\rm lk}_{S^{\,e -1}_{e+1}}(\alpha_2)$ and
hence $\beta \in B^{\,e}_{e+1} \subseteq \widetilde{S}^{\,e +4}_{e
+ 17}$. This is a contradiction since $\overline{\alpha} \ast
\partial\beta$ is an induced subcomplex of
$\widetilde{S}^{\,e +4}_{e+17}$. So, $\alpha_2$ is a facet of
$S^{\,e-1}_{e+1}$ and hence $\partial \beta = {\rm
lk}_{B^{\,4}_{16}}(\alpha_1)$. Then, $2 \leq \dim(\alpha_1)= 4 -
\dim(\beta) \leq 2$ and hence $\dim(\alpha_1) = \dim(\beta) =2$.
Let $\alpha_1= xuv$ (where $x$ is the fixed vertex chosen in
$S^{\,3}_{16}$ to construct $B^{\,4}_{16}$). Then ${\rm
lk}_{S^{\,3}_{16}}(uv) = {\rm lk}_{B^{\,4}_{16}}(\alpha_1) =
\partial\beta$. This is not possible since $S^{\,3}_{16}$ does not
contain any edge of degree 3.

Thus $\alpha_1$ is a boundary face of $B^{\,4}_{16}$, i.e.,
$\alpha_1 \in S^{\,3}_{16}$. If $\alpha_2$ is the facet of $B^{\,e
}_{e+1}$ then ${\rm lk}_{S^{\,3}_{16}}(\alpha_1) = {\rm
lk}_{\widetilde{S}^{\,e+4}_{e+17}}(\alpha) = \partial \beta$.
Hence $\dim(\alpha_1) \geq 2$ and therefore $\dim(\beta) \leq 1$,
a contradiction. So, $\alpha_2$ is not the facet of
$B^{\,e}_{e+1}$ (and hence ${\rm lk}_{B^{\,e}_{e+1}}(\alpha_2)$
is a standard ball). Thus, the ball $B_1 := {\rm
lk}_{B^{\,4}_{16}}(\alpha_1) \ast {\rm lk}_{B^{\,e}_{e+
1}}(\alpha_2)$ is a non-trivial join of balls, so that all the
vertices of $B_1$ are in its boundary. But, $\partial B_1 = {\rm
lk}_{\widetilde{S}^{\,e+4}_{e+17}}(\alpha) =
\partial \beta$. Therefore, $B_1$ is the standard ball
$\overline{\beta}$ and hence ${\rm lk}_{B^{\,4}_{16}}(\alpha_1)$
is a standard ball. Therefore, ${\rm lk}_{S^{\,3}_{16}}(\alpha_1)$
is a standard sphere and hence $\dim(\alpha_1) \geq 2$. So, ${\rm
lk}_{B^{\,4}_{16}}(\alpha_1)$ is a standard ball of dimension
$\leq 1$, i.e., it is a vertex or an edge. Then the vertex set of
${\rm lk}_{B^{\,4}_{16}}(\alpha_1)$ is a face in $S^{\,3}_{16}$.
So, the vertex set $\beta$ of ${\rm lk}_{\widetilde{S}^{\,e
+4}_{e+17}}(\alpha)$ is a face of $S^{\,3}_{16} \ast B^{\,e}_{e+1}
\subseteq \widetilde{S}^{\,e+4}_{e+17}$. Therefore,
$\overline{\alpha} \ast \partial\beta$ is not an induced
subcomplex of $\widetilde{S}^{\,e+4}_{e+17}$, a contradiction.
Thus, for each $e\geq 0$, $\widetilde{S}^{\,e+4}_{e+17}$ is an
unflippable combinatorial $(e+4)$-sphere.

From this claim, it follows that $\partial( B^{\,4}_{16} \ast
B^{\hspace{.2mm}d-4}_{d-3})$ is a combinatorial $d$-sphere
which is not $d$-stellated. Since $B^{\,4}_{16}$ is a 2-stacked
ball, it follows that $B^{\,4}_{16} \ast B^{\hspace{.2mm}d-4}_{d
-3}$ is also 2-stacked. This implies that $\partial( B^{\,4}_{16}
\ast B^{\hspace{.2mm}d-4}_{d-3})$ is a 2-stacked combinatorial
$d$-sphere which is not $d$-stellated, for $d\geq 4$.
From this and the observation in (b), we find that for
each $d\geq 3$, there are $2$-stacked combinatorial $d$-spheres
which are not $d$-stellated.

Since the classes $\Sigma_k(d)$,
${\cal S}_k(d)$ are increasing in $k$, we get\,:}
\begin{itemize}
\item[$\bullet$] For $2 \leq k \leq l \leq d \geq 3$ there are
$k$-stacked combinatorial $d$-spheres which are not $l$-stellated.

\item[$\bullet$] For $2 \leq k \leq l \leq d+1 \geq 6$ there are
$k$-stacked triangulated $d$-spheres which are not $l$-stellated.
\end{itemize}
\item[{\bf (e)}] {\rm Let $S^{\hspace{.2mm}3}_{10}$ be the pure
simplicial complex of dimension three whose vertices are the
digits $0, 1, \dots, 9$ and whose facets are\,:
\begin{eqnarray*}
0123, 1234, 2345, 3456, 4567, 5678, 6789, 0128, 0139, 0189,
0238, 0356, 0358, 0369,  \\
0568, 0689, 1248, 1349, 1457, 1458,
1467, 1469, 1578, 1679, 1789, 2358, 2458, 3469.
\end{eqnarray*}
Let $S^{\hspace{.2mm}2}_{10}$ be the pure 2-dimensional
subcomplex of $S^{\hspace{.2mm}3}_{10}$ whose facets are\,:
$$012, 013, 023, 124, 134, 235, 245, 346, 356, 457, 467,
568, 578, 679, 689, 789.
$$
Then $S^{\hspace{.2mm}3}_{10}$ is a triangulated 3-sphere, and
$S^{\hspace{.2mm}2}_{10}$ is a triangulated 2-sphere embedded in
$S^{\hspace{.2mm}3}_{10}$. Being two-sided in
$S^{\hspace{.2mm}3}_{10}$, the ``equatorial"
$S^{\hspace{.2mm}2}_{10}$ divides $S^{\hspace{.2mm}3}_{10}$ into
two closed ``hemispheres", say $B_1$ and $B_2$. Of course, $B_1$
and $B_2$ are triangulated 3-balls. The facets of $B_1$ are the
first seven facets of $S^{\hspace{.2mm}3}_{10}$, while the facets
of $B_2$ are the remaining twentyone facets of
$S^{\hspace{.2mm}3}_{10}$.

The dual graph of the 3-ball $B_1$ is visibly a path. So, by
Proposition \ref{P7}, $B_1$ is 1-stacked. Since (from the above
discussion, or by direct verification) $\partial B_1 =
S^{\hspace{.2mm}2}_{10} = \partial B_2$, it follows that
$S^{\hspace{.2mm}2}_{10}$ is 1-stellated. But, it also bounds the
ball $B_2$ which is Ziegler's example \cite{z2} of a non-shellable
3-ball\,! (If $\alpha$ is a facet of a triangulated $d$-ball $B$,
then one says $\alpha$ is an {\em ear} of $B$ if $B \setminus
\{\alpha\}$ is also a triangulated $d$-ball. Clearly, if $B$ is
shellable, then the last facet, added while obtaining $B$ from
$B^{\hspace{.2mm}d}_{d+1}$ by a sequence of shelling moves, must
be an ear of $B$. Thus, if $B$ has no ears, then it must be
non-shellable. Such balls are ``strongly non-shellable" in the
terminology of Ziegler. A facet $\alpha$ of $B$ is an ear of $B$
if and only if the induced subcomplex of $\partial B$ on the
vertex set $\alpha$ is a $(d-1)$-ball. Using this criterion, it is
possible to verify that $B_2$ has no ears\,: it is strongly
non-shellable.) }

\item[{\bf (f)}] {\rm The following example of a shellable 3-ball
with a unique ear is due to Frank Lutz (personal communication).
Consider the pure 3-dimensional 2-neighbourly simplicial complex
$S^{\hspace{.2mm}3}_{8}$ with vertices $1, 2, \dots, 8$ and facets
\begin{eqnarray*}
1234, 2345, 3456, 4567, 5678, 1237, 1248, 1278, 1348, 1356, \\
1357, 1368, 1568, 1578, 2357, 2457, 2467, 2468, 2678, 3468.
\end{eqnarray*}
Let $S^{\hspace{.2mm}2}_{8}$ be the pure 2-dimensional subcomplex
of $S^{\hspace{.2mm}3}_{8}$ with facets
$$
123, 124, 134, 235, 245, 346, 356, 457, 467, 568, 578, 678.
$$
Again, $S^{\hspace{.2mm}2}_{8}$ is a triangulated 2-sphere
embedded in the triangulated 3-sphere $S^{\hspace{.2mm}3}_{8}$. As
in (e) above, $S^{\hspace{.2mm}2}_{8}$ divides
$S^{\hspace{.2mm}3}_{8}$ into two 3-balls $B_1$ and $B_2$. The
facets of $B_1$ are the first five facets of $S^{\hspace{.2mm}3}_{8}$,
while the facets of $B_2$ are the remaining fifteen facets of
$S^{\hspace{.2mm}3}_{8}$. Again, $B_1$ is an 1-stacked 3-ball
since its dual graph is a path. We have $\partial B_1 =
S^{\hspace{.2mm}2}_{8} = \partial B_2$. Thus,
$S^{\hspace{.2mm}2}_{8}$ is an 1-stellated sphere. The other ball
$B_2$ bounded by $S^{\hspace{.2mm}2}_{8}$ is shellable (indeed,
2-shelled). (A shelling of $B_2$\,: 1357, 1356, 1368, 1348, 1248,
3468, 1568, 1578, 1278, 2468, 2678, 1237, 2467, 2357, 2457.) But,
$B_2$ has only one ear, namely $2457$.

Clearly, the class ${\cal S}_k(d)$ of $k$-stacked $d$-spheres is
closed under connected sum. In consequence, the class
$\Sigma_1(d)$ of 1-stellated $d$-spheres is closed under connected
sums. However, consider the following construction. Take a
standard 2-ball $B^{\hspace{.2mm}2}_{3}$ with a vertex set $\{a,
b, c\}$ disjoint from $V(S^{\hspace{.2mm}3}_{8})$, and form the
join $B := B_2 \ast B^{\hspace{.2mm}2}_{3}$. Then $B$ is a
2-shelled 6-ball with a unique ear $2457abc$. Thus, $S := \partial
B$ is a 2-stellated 5-sphere. The facets $245abc$, $457abc$ are
two of the facets of $S$ in the unique ear of $B$. Take a vertex
disjoint copy $B^{\hspace{.2mm}\prime}$ of $B$, and let
$S^{\hspace{.2mm}\prime} = \partial B^{\hspace{.2mm}\prime}$, the
corresponding copy of $S$. Let $1^{\hspace{.1mm}\prime}, \dots,
8^{\hspace{.2mm}\prime}, a^{\hspace{.1mm}\prime},
b^{\hspace{.1mm}\prime}, c^{\hspace{.1mm}\prime}$ be the vertices
of $B^{\hspace{.2mm}\prime}$ corresponding to the vertices $1,
\dots, 8, a, b, c$ respectively. Form the connected sum
$\widetilde{B} = B \# B^{\hspace{.2mm}\prime}$ by doing the
identifications $2 \equiv 2^{\hspace{.2mm}\prime}, 4 \equiv
4^{\hspace{.2mm}\prime}, 5 \equiv 5^{\hspace{.2mm}\prime}, a
\equiv a^{\hspace{.2mm}\prime}, b \equiv b^{\hspace{.2mm}\prime},
c \equiv c^{\hspace{.2mm}\prime}$. Then $\widetilde{B}$ is a
16-vertex non-shellable 2-stacked 6-ball. Let $\widetilde{S} =
\partial \widetilde{B}$. Then $\widetilde{S}$ is a 16-vertex
2-stacked 5-sphere which is not 2-stellated (by Propositions
\ref{P6} and \ref{P4}). (It can be shown that $\widetilde{S}$ is
5-stellated.) But, $\widetilde{S} = S \# S^{\hspace{.2mm}\prime}$,
the connected sum of two 2-stellated 5-spheres. For $d\geq 5$, if
we take $B^{\,d- 3}_{d-2}$ in place of $B^{\,2}_3$ in the above
construction then, by the same argument,   we get a $d$-sphere which is not
2-stellated and is the connected sum of two 2-stellated
$d$-spheres. Thus
\begin{itemize}
\item[$\bullet$] For $d\geq 5$, the class $\Sigma_2(d)$ is not
closed under connected sum.
\end{itemize}

}
\end{itemize}
\end{eg}

By Proposition \ref{P4}, all the $k$-stellated spheres of
dimension $d \geq 2k-1$ are $k$-stacked. But, we are so far unable
to answer\,:

\begin{qn}\hspace{-1.8mm}{\bf .} \label{Q1}
{\rm Is there a $k$-stellated $d$-sphere which is not
$k$-stacked\,?}
\end{qn}
Note that, by Propositions \ref{P1} and \ref{P4}, for an
affirmative answer to Question \ref{Q1}, we must have $k+1 \leq d
\leq 2k-2$, and hence $k\geq 3$, $d\geq 4$.

Notice that any $(k+1)$-neighbourly $d$-sphere is (trivially)
$(d-k)$-stacked. A comparison of this observation with Proposition
\ref{P12} as well as a comparison between Proposition
\ref{P6} and Corollary \ref{C6} leads us to a strong suspicion\,:

\begin{conj}\hspace{-1.8mm}{\bf .} \label{Conj1}
For $d\geq 2k$, a polytopal $d$-sphere is $k$-stellated if
$($and only if\,$)$ it is $k$-stacked. Equivalently $($in view of
Corollary $\ref{C6}$ and Propositions $\ref{P2}$ and $\ref{P4})$, if
$S$ is a $k$-stacked polytopal sphere of dimension $d \geq 2k$,
then the $(d+1)$-ball $\overline{S}$ $($given by formula $(\ref{eq1}))$
is shellable.
\end{conj}

Let $S = S^{\hspace{.2mm}k-1}_{k+1} \ast S^{\hspace{.2mm}k-
1}_{k+1}$, $B_1 = S^{\hspace{.2mm}k-1}_{k+1} \ast
B^{\hspace{.2mm}k}_{k +1}$ and $B_2 = B^{\hspace{.2mm}k}_{k+1}
\ast S^{\hspace{.2mm}k-1}_{k+1}$. Then $B_1$, $B_2$ are
$k$-stacked polytopal $2k$-balls with $\partial B_1 = S = \partial
B_2$. Thus, $S$ is a $(2k-1)$-dimensional $k$-neighbourly
polytopal $k$-stacked sphere. Hence $S$ is $k$-stellated by
Proposition \ref{P12}. Thus, $S$ is an example of a $(2k-
1)$-dimensional $k$-stellated polytopal sphere which bounds two
distinct (though isomorphic) $k$-stacked balls. So, the bound
$d\geq 2k$ in Proposition \ref{P6} and Corollary \ref{C6} is
sharp.

A comparison of Propositions \ref{P5} and \ref{P6} above leads us
to the following query.

\begin{qn}\hspace{-1.8mm}{\bf .} \label{Q2}
{\rm (a) Are there $k$-stacked balls $B_1$, $B_2$ of dimension
$2k+1$ such that $B_1\neq B_2$ but $\partial B_1 = \partial
B_2$\,?

(b) If $S$ is a $k$-stacked sphere of dimension $\geq 2k+1$ then
by Proposition \ref{P5} there is a unique $k$-stacked ball
$\overline{S}$ such that $\partial \overline{S} = S$. Is
$\overline{S}$ always given by the formula (\ref{eq1})\,?}
\end{qn}

\begin{eg}[The Klee-Novik construction]\hspace{-1.5mm}{\bf .}
\label{E2} {\rm For $d\geq 1$, let $S^{\hspace{.2mm}d+1}_{2d+4}$
be the join of $d+2$ copies of $S^{\hspace{.2mm}0}_{2}$ with
disjoint vertex sets $\{x_i, y_i\}$, $1\leq i \leq d+2$. Then
$S^{\hspace{.2mm}d+1}_{2d+4}$ is a triangulated sphere with
missing edges $x_iy_i$, $1\leq i \leq d+2$ (cf. Example
\ref{E1}\,(a)). Each of the $2^{d+2}$ facets of $S^{\hspace{.2mm}d
+1}_{2d+4}$ may be encoded by a sequence of $d+2$ signs as
follows. If $\sigma$ is a facet, then for each index $i$ ($1\leq i
\leq d+2$) $\sigma$ contains either $x_i$ or $y_i$, but not both.
Put $\varepsilon_i = +$ if $x_i\in \sigma$ and $\varepsilon_i = -$
if $y_i\in \sigma$. Thus the sign sequence $(\varepsilon_1, \dots,
\varepsilon_{d+2})$ encodes the facet $\sigma$. For $0\leq k\leq
d$, let $\overline{M}(k, d)$ be the pure $(d+1)$-dimensional
subcomplex of $S^{\hspace{.2mm}d+1}_{2d+4}$ whose facets are those
facets $\sigma$ (of the latter complex) whose sign sequences have
at most $k$ sign changes. (A sign change in the sign sequence
$(\varepsilon_1, \dots, \varepsilon_{d+2})$ is an index $1\leq i
\leq d+1$ such that $\varepsilon_{i+1} \neq \varepsilon_{i}$.)
Then $\overline{M}(k, d)$ is a pseudomanifold with boundary.
Klee and Novik \cite{kn} proved that $M(k, d) := \partial
\overline{M}(k, d)$ is a triangulation of $S^{\hspace{.2mm}k}
\times S^{\hspace{.2mm}d-k}$ for $0\leq k \leq d$. (In their
paper, Klee and Novik use the notation $B(k, d+2)$ for
$\overline{M}(k, d)$.) The authors of \cite{kn} observed that
the permutation $D$, $E$ and $R$ are automorphisms of
$\overline{M}(k, d)$ (and hence of $M(k, d)$), where $D =
{\displaystyle \prod_{j=1}^{d+2}(x_j, y_j)}$, $E =
\hspace{-3mm}{\displaystyle \prod_{1\leq j< (d+3)/2}\hspace{-3mm}
(x_j, x_{d+3-j})(y_j, y_{d+3-j})}$ and $R = (x_1, \dots,
x_{d+2})(y_1, \dots, y_{d+2})$ when $k$ is even, $R = (x_1, \dots,
x_{d+2}, y_1, \dots, y_{d+2})$ when $k$ is odd. Clearly, these
three automorphisms generate a vertex-transitive automorphism
group of $\overline{M}(k, d)$. Therefore, the links in
$\overline{M}(k, d)$ (or in $M(k, d)$) of all the vertices are
isomorphic. The involution $A = {\displaystyle \prod_{j ~ {\rm
even}}(x_j, y_j)}$ is an isomorphism between $M(k, d)$ and $M(d-k,
d)$. Therefore, in discussing these constructions we may (and do)
assume $d\geq 2k$. (However, $A$ is not an isomorphism between
$\overline{M}(k, d)$ and $\overline{M}(d-k, d)$. Indeed, $A$
maps $\overline{M}(k, d)$ to the ``complement" of
$\overline{M}(d-k, d)$ in $S^{\hspace{.2mm}d+1}_{2d+4}$.)

Let $I = \{1, 2, \dots, d+1\}$. Define the linear order $\prec$ on
${I \choose \leq \,k\,}$ by\,: $A \prec B$ if either $\#(A) <
\#(B)$ or else $\#(A) = \#(B)$, $A <_{\rm lex} B$, where $<_{\rm
lex}$ is the usual lexicographic order. Let $L$ be the link of the
vertex $x_{d+2}$ in $\overline{M}(k, d)$. Clearly, for each $A \in
{I \choose \leq \,k\,}$, there is a unique facet $\tau_A$ of $L$
such that $A$ is precisely the set of sign-changes corresponding
to the facet $\tau_A\cup \{x_{d+2}\}$ of $\overline{M}(k, d)$. We
may transfer the linear order $\prec$ to the set of facets of $L$
via the bijection $A \mapsto \tau_A$. Then, Klee and Novik show in
\cite{kn} that $\prec$ is a shelling order for $L$. Thus, $L$ is a
shellable $d$-ball. What is more, if $\#(A) = j \leq k$ then the
facet $\tau_A$ of $L$ is obtained (from the $d$-ball with facets
$\tau_B$, $B\prec A$) by a shelling move of index $j-1$. In
consequence, $L$ is a $k$-shelled $d$-ball. Since the automorphism
group of $\overline{M}(k, d)$ is vertex transitive, it follows
that all vertex links of $\overline{M}(k, d)$ are $k$-shelled
$d$-balls. Thus, $\overline{M}(k, d)$ is a $(d+1)$-manifold with
boundary. Also, since the boundary of a $k$-shelled ball is a
$k$-stellated sphere (Corollary \ref{C1}), it follows that $M(k,
d) = \partial \overline{M}(k, d)$ has $k$-stellated vertex links.
Thus,
\begin{itemize}
\item[$\bullet$] $M(k, d)\in {\cal W}_k(d)$ for $d\geq 2k$.
\end{itemize}
Also note that, when $d\geq 2k+1$, the vertex links of
$\overline{M}(k, d)$ are the unique (Proposition \ref{P6})
$k$-stacked balls bounded by
the corresponding vertex links of $M(k, d)$. Therefore,
$\overline{M}(k, d)$ is the unique $(d+1)$-manifold $\overline{M}$
such that $\partial \overline{M} = M(k, d)$ and ${\rm
skel}_{d-k}(\overline{M}) = {\rm skel}_{d-k}(M(k, d))$. (In
consequence, when $d\geq 2k+2$, $\overline{M}(k, d)$ may be
recovered from $M(k, d)$ via the formula (\ref{eq2}) above; cf.
Proposition \ref{P9}.) Therefore, for $d \geq 2k+1$, every
automorphism of $M(k, d)$ extends to an automorphism of
$\overline{M}(k, d)$\,: they have the same automorphism group.
However, it is elementary to verify that the full automorphism
group of $\overline{M}(k, d)$ is of order $4d+8$. (Since this
group is transitive on the $2d+4$ vertices, it suffices to
show that the full stabilizer of the vertex $x_{d+2}$ is
of order 2. This is easy.) Thus,
\begin{itemize}
\item[$\bullet$] When $d\geq 2k+1$, the full automorphism
group of $M(k, d)$ is of order $4d+8$ (namely, the group
generated by $D$, $E$, $R$ above).
\end{itemize}
This leaves open the following tantalizing question.
}
\end{eg}

\begin{qn}\hspace{-1.8mm}{\bf .} \label{Q3}
{\rm What is the full automorphism group of $M(k, 2k)$\,?}
\end{qn}

Notice that the involution $A$ defined above is also an
automorphism of $M(k, 2k)$. However, since $A$ maps
$\overline{M}(k, 2k)$ to its complement in $S^{\hspace{.2mm}d
+1}_{2d+4}$, $A$ is not an automorphism of $\overline{M}(k, 2k)$.
Therefore, $A \not\in H := \langle D, E, R\rangle$. The
automorphism $A$ normalizes $H$, so that the group $G := \langle
D, E, R, A\rangle$ is of order $2\times \#(H) = 16(k+1)$. We
suspect that $G$ is the full automorphism group of $M(k, 2k)$. Is
it\,?

For any vertex $x$ of $M(k, d)$, let $L_x$ (respectively
$\overline{L}_x$) be the link of $x$ in $M(k, d)$ (respectively in
$\overline{M}(k, d)$). By the discussion above, $\overline{L}_x$
is obtained from $B^{\hspace{.2mm}d}_{d+1}$ by $\sum_{j=1}^k {d+1
\choose j}$ shelling moves, of which ${d+1 \choose j}$ are of
index $j-1$. Therefore, by Lemma \ref{L1}, $L_x = \partial
\overline{L}_x$ is obtained from $S^{\hspace{.2mm}d-1}_{d+1}$ by
$\sum_{j=1}^k {d+1 \choose j}$ bistellar moves, of which ${d+1
\choose j}$ are of index $j-1$. Since $L_x$ is a $k$-stellated
sphere of dimension $d-1 \geq 2k-1$, Proposition \ref{P14} implies
that $g_j(L_x) = {d+1 \choose j}$ for $0\leq j \leq k$. Since
$M(k, d)$ has $2d+4$ vertices, Lemma \ref{L4} implies that the
$g$-vector of $M(k, d)$ satisfies the recurrence $(d+2-j)g_j +
(j+1)g_{j+1} = {d+1 \choose j}(2d+4)$, $0\leq j \leq k$. Solving
this recurrence relation (with initial condition $g_0 =1$) we see
that the $g$-vector of $M(k, d)$ satisfies $g_j = {d+2 \choose j}$
for $0\leq j \leq k+1$. Since $d\geq 2k$ and $M(k, d)\in {\cal
W}_k(d)$, Corollary \ref{C10} determines the entire $g$-vector of
$M(k, d)$ from the above computation. By Proposition \ref{P15}, 
we find\,:
\begin{itemize}
\item[$\bullet$] $g_{l+1}(M(k, d)) = \left\{
\begin{array}{cl}
{d+2 \choose l+1} & \mbox{if} ~~ 0\leq l \leq k, \\[1.5mm]
(-1)^{l-k} {d+2 \choose l+1} & \mbox{if} ~~  k+1 \leq l \leq d-k.
\end{array}
\right.$
\end{itemize}
(The rest of the $g$-numbers may now be determined using Klee's
formula, namely (\ref{eq7}). These formulae are in agreement with
Theorem 5.2 in \cite{kn}, of course. )

The calculation above shows that, for $d\geq 2k$, $M(k,d)$
satisfies the inequalities in Proposition \ref{P21} (indeed, with
equality for $k\leq l \leq d-k-1$) even though $M(k, d)$ is not
2-neighbourly. (Actually, the $k$-skeleton of $M(k, d)$ agrees
with that of $S^{\hspace{.2mm}d +1}_{2d+4}$). We suspect that the
assumption of 2-neighbourliness in Proposition \ref{P21} should be
removable. These and other examples lead us to posit\,:

\begin{conj}[GLBC for triangulated manifolds]\hspace{-1.5mm}{\bf .}
\label{Conj2} Let $M$ be any triangulation of a connected closed
$d$-manifold. Then, for $1\leq l \leq \frac{d- 1}{2}$, the
$g$-numbers of $M$ satisfy $g_{l+1}(M) \geq \sum_{i=1}^l
(-1)^{l-i} \, \beta_i(M; \FF)$. Further, equality holds here for
some $l < \frac{d-1}{2}$ if and only if $M \in {\cal K}_l(d)$.
\end{conj}

Notice that Proposition \ref{P21} proves the inequality of
Conjecture \ref{Conj2} and the ``if\," part of the equality case
of the conjecture under the extra assumption that $M$ is
2-neighbourly and all the vertex links of $M$ are $\lfloor
\frac{d-1}{2}\rfloor$-stellated. The Klee-Novik manifolds $M(k, d)$,
$d \geq 2k$, satisfy all parts of Conjecture \ref{Conj2}.
The ``$l=1$" case of
Conjecture \ref{Conj2} (with $\FF = \QQ$) was a conjecture of
Kalai \cite{ka}; Novik and Swartz proved it (for any field $\FF$,
with the extra hypothesis that $M$ is $\FF$-orientable) in
\cite[Theorem 5.2]{ns}. Since Conjecture \ref{Conj2} includes the
GLBC for homology spheres, we do not expect it to be settled in a
hurry.

We expect that part (b) of Proposition \ref{P23} ($\equiv$
Theorem 5 in \cite{lss}) should generalize as follows (compare
K\"{u}hnel's conjecture \cite[Conjecture 18]{lu2})\,:

\begin{conj}\hspace{-1.8mm}{\bf .} \label{Conj3}
If $M$ is an $m$-vertex connected closed triangulated $d$-manifold
with Betti numbers $\beta_i$ $($with respect to some field\,$)$
then ${\,m+l-d-2\, \choose l+1} \geq {\,d+2\, \choose l+1}
{ \sum_{i=1}^l}(-1)^{l-i}\beta_i$ for $1 \leq l \leq (d-1)/2$.
Also, equality holds here for some $\,l< (d-1)/2$ if and only
if $M \in {\cal K}_l^{\ast}(d)$.
\end{conj}

In this connection, we may ask\,:

\begin{qn}\hspace{-1.8mm}{\bf .} \label{Q4}
{\rm Is it true that the $g$-vector of any triangulated closed
$d$-manifold on $m$ vertices satisfies $g_{l+1} \leq {m+l-d-2
\choose l+1}$, with equality (if and) only if the triangulation
is $(l+1)$-neighbourly\,? If this is true, then, of course,
Conjecture \ref{Conj2} implies Conjecture \ref{Conj3}.}
\end{qn}

By Proposition \ref{P24}, for $d\geq 2k+2\geq 6$, any member of
${\cal W}_k^{\ast}(d)$ has the same $\ZZ$-homology as the
connected sum of copies of $S^{\hspace{.2mm}k} \times
S^{\hspace{.2mm}d -k}$. This raises the question\,:
\begin{qn}\hspace{-1.8mm}{\bf .} \label{Q5}
{\rm Is it true that for $d \geq 2k+2$, any member of ${\cal
W}_k^{\ast}(d)$ triangulates a connected sum of (the total spaces of)
$S^{\hspace{.2mm}d -k}$-bundles over $S^{\hspace{.2mm}k}$\,?}
\end{qn}

\begin{eg}[Tight triangulations of closed manifolds]\hspace{-1.5mm}{\bf .}
\label{E3} {\rm We have noted that $S^{\hspace{.2mm}d}_{d +2}$ is
the only tight triangulation of $S^{\hspace{.2mm}d}$. This trivial
series apart, we know the following examples of tight
triangulations (cf. \cite{kl}).
\begin{itemize}
\item[{\bf (a)}] By Lemma \ref{L10}, all 2-neighbourly
2-dimensional closed triangulated manifolds are tight when
orientable and $\ZZ_2$-tight when non-orientable. For $n \geq 4$, there exist
$n$-vertex 2-neighbourly orientable (respectively, non-orientable)
triangulated 2-manifolds if and only if $n\equiv 0, 3, 4$ or 7
(mod 12) (respectively, $n \equiv 0$ or 1 (mod 3), except for $n =
4, 7$) (cf. \cite{r}).

\item[{\bf (b)}] For each $d\geq 2$, there is a $(2d+3)$-vertex
member $K^{d}_{2d+3}$ of ${\cal W}_1^{\ast}(d)$ found by
K\"{u}hnel \cite{k86}. For $d\geq 3$, it is the unique non-simply
connected $d$-manifold on $2d+3$ vertices (cf. \cite{bd8, css}).
It is orientable (triangulates $S^{\,d-1}\times S^{1}$) for $d$
even and non-orientable (triangulates $\TPSSDS$) for $d$ odd. By Proposition \ref{P22},
$K^{d}_{2d+3}$ is tight for $d$ even, and $\ZZ_2$-tight for $d$
odd.

\item[{\bf (c)}] \begin{enumerate} \item[(i)] The 15-vertex
triangulation of $(\TPSSF)^{\#3}$ obtained (in \cite{bd10}) is in
${\cal W}_1^{\ast}(4)$, hence $\ZZ_2$-tight by Proposition
\ref{P22}.
\item[(ii)] Recently, Nitin Singh, a student of the
second author, modified this construction to obtain (in \cite{si})
two 15-vertex triangulations of $(S^{\,3} \times S^{1})^{\#3}$ in
${\cal W}_1^{\ast}(4)$. Both are tight by Proposition \ref{P22}.
\end{enumerate}
\item[{\bf (d)}] Lutz constructed (in \cite{lu1}) two 12-vertex
triangulations of $S^{\,2} \times S^{\,3}$; they belongs to ${\cal
W}_2^{\ast}(5)$. By Proposition \ref{P25}, these are tight
triangulations.

\item[{\bf (e)}] Only finitely many $2k$-dimensional $(k+
1)$-neighbourly triangulated closed manifolds are known for $k
\geq 2$. By Lemma \ref{L10}, they are all tight. These examples
are\,: \begin{enumerate}
\item[(i)] The 9-vertex triangulation $\CC P^{\,2}_9$ of $\CC
P^{\,2}$ due to K\"{u}hnel (in \cite{kb}),
\item[(ii)] the
16-vertex triangulation of a $K3$-surface due to Casella and
K\"{u}hnel (in \cite{ck}), \item[(iii)] two 13-vertex triangulations of
$S^{\,3} \times S^{\,3}$ due to Lutz (in \cite{lu1}), and
\item[(iv)] six
15-vertex triangulations of homology $\HH P^{\,2}$ (three due to
Brehm and K\"{u}hnel in \cite{bk2} and three due to Lutz in
\cite{lu2}).
\end{enumerate}

\item[{\bf (f)}] Apart from the above list, we know only two tight
triangulated manifolds. These are\,: \begin{enumerate}
\item[(i)] A 15-vertex triangulation of $(\TPSSF)
\# (\CC P^{\,2})^{\#5}$ due to Lutz (in \cite{lu1}).
It is 2-neighbourly, non-orientable, $\ZZ_2$-tight and in ${\cal W}_2(4)$.
\item[(ii)] A 13-vertex triangulation of $SU(3)/SO(3)$ due to Lutz (in \cite{lu1}).
It is 3-neighbourly, orientable, $\ZZ_2$-tight and in ${\cal W}_3(5)$.
\end{enumerate}
The tightness of these two examples do not follow from the results presented here.
\end{itemize}
Corollary \ref{C11} implies that all the triangulations in Example
\ref{E3}\,(a), (b) and (c) are strongly minimal. By Theorem 4.4 of
\cite{ns}, all the
triangulations in Example \ref{E3}\,(d) and (e) are minimal.
As far as we know, the minimality of the
triangulations in Example \ref{E3}\,(f) is an open problem.
}
\end{eg}

Finally, we have the question of how to get new triangulations
meeting the hypothesis of Propositions \ref{P22} and \ref{P25}. In
particular, we may ask\,:

\begin{qn}\hspace{-1.8mm}{\bf .} \label{Q6}
{\rm Is there a 20-vertex triangulation of $(S^{\,2} \times
S^{1})^{\#12}$ or $(\TPSS)^{\#12}$ in ${\cal W}_1^{\ast}(3)$
or a 20-vertex triangulation of $(S^{\,3}\times S^{\,2})^{\#13}$
in ${\cal W}_2^{\ast}(5)$\,?}
\end{qn}

We might also ask if there are examples of triangulations from
${\cal W}_k(d) \setminus {\cal W}_{k-1}(d)$ satisfying the
hypothesis of Proposition \ref{P26} (cf. Remark \ref{R5}). We
suspect that there may not exist any such examples!

We do not know for a fact that, for $1<l< (d-1)/2$, the members
of ${\cal K}_l(d)$ (or even of ${\cal W}_l(d)$) actually attain
equality in Conjecture \ref{Conj2}. Thus, all parts of this conjecture
are wide open for $l>1$.
Notice that, as a consequence of Proposition \ref{P21}, the members
of ${\cal W}_l^{\ast}(d)$ do attain equality in Conjecture \ref{Conj3}
for $1\leq l < (d-1)/2$. However, we do not know if,
more generally, the members of ${\cal K}_l^{\ast}(d)$ attain these equalities for $1\leq l <
(d-1)/2$, as conjectured. The case $l=1$ is unrevealing since in
this case ${\cal W}_1(d) = {\cal K}_1(d)$ and ${\cal W}_1^{\ast}(d)
= {\cal K}_1^{\ast}(d)$ (Corollary \ref{C4}). Thus, the most important
question raised by this paper is whether (and to what extent)
the results can be extended from $k$-stellated spheres to $k$-stacked
spheres. A good place to begin this investigation is to address
the following\,:

\begin{qn}\hspace{-1.8mm}{\bf .} \label{Q7}
{\rm Is Proposition \ref{P19} (on the sigma-vector of $k$-stellated
spheres) valid for $k$-stacked spheres\,?}
\end{qn}

\addcontentsline{toc}{section}{~~~ Acknowledgement}

\noindent {\bf Acknowledgement\,:} The authors thank Frank H. Lutz
for drawing their attention to the reference \cite{dfm} and for
providing the unique ear 3-ball $B_2$ of Example \ref{E1}\,(f).




\addcontentsline{toc}{section}{~~~ References}

{\footnotesize

\begin{thebibliography}{99}
\bibitem{as}
A. Altshuler, L. Steinberg, Neighborly 4-polytopes with 9
vertices, {\em J. Combin. Theory} (A) {\bf 15} (1973), 270--287.

\bibitem{bd5}
B. Bagchi, B. Datta, Combinatorial triangulations of homology
spheres, {\em Discrete Math.} {\bf 305} (2005), 1--17.

\bibitem{bd8}
B. Bagchi, B. Datta, Minimal triangulations of sphere
bundles over the circle, {\em J. Combin. Theory} (A) {\bf
115} (2008), 737--752.

\bibitem{bd9}
B. Bagchi, B. Datta, Lower bound theorem for normal
pseudomanifolds, {\em  Expositiones Math.} {\bf 26} (2008),
327--351.

\bibitem{bd10}
B. Bagchi, B. Datta, On Walkup's class ${\cal K}(d)$ and a minimal
triangulation of $(S^{\hspace{.2mm}3} \mbox{$\times
\hspace{-2.5mm}_{-}$} \, S^{\hspace{.1mm}1})^{\#3}$, {\em Discrete
Math.} {\bf 311} (2011), 989--995.

\bibitem{bd13_v1}
B. Bagchi, B. Datta, On stellated spheres and a tightness
criterion for combinatorial manifolds, arXiv:\,1102.0856\,v1,
2011, 23 pages.

\bibitem{bl}
A. Bj\"{o}rner, F. H. Lutz, Simplicial manifolds, bistellar flips
and a 16-vertex triangulation of Poincar\'{e} homology 3-sphere,
{\em Experiment. Math.} {\bf 9} (2000), 275--289.

\bibitem{bk}
U. Brehm, W. K\"{u}hnel, Combinatorial manifolds with few
vertices, {\em Topology} {\bf 26} (1987), 465--473.

\bibitem{bk2}
U. Brehm, W. K\"{u}hnel, 15-vertex triangulations of an
8-manifold, {\em Math. Annalen} {\bf 294} (1992), 167--193.

\bibitem{ca}
J. W. Cannon, Shrinking cell-like decomposition of manifolds:
codimension three, {\em Ann. Math.} {\bf 110} (1979), 83--112.

\bibitem{ck}
M. Casella, W. K\"{u}hnel, A triangulated K3 surface with the
minimum number of vertices, {\em Topology} {\bf 40} (2001),
753--772.

\bibitem{css} J. Chestnut, J. Sapir, E. Swartz,
Enumerative properties of triangulations of spherical boundles
over $S^1$, {\em Euro. J. Combin.} {\bf 29} (2008), 662--671.

\bibitem{dk}
G. Danaraj, V. Klee, Shellings of spheres and polytopes, {\em Duke
Math. J.} {\bf 41} (1974), 443--451.

\bibitem{d}
J. Dancis, Triangulated $n$-manifolds are determined by their
$([n/2]+1)$-skeletons, {\em Topo. Appl.} {\bf 18} (1984), 17--26.

\bibitem{dfm} R. Dougherty, V. Faber, M. Murphy, Unflippable
tetrahedral complexes, {\em  Discrete Comput. Geom.} {\bf 32}
(2004), 309-–315.

\bibitem{ed} R. D. Edwards, The topology of manifolds and cell-like
maps, {\em Proc. I. C. M.} (Helsinki, 1978), pp. 111-–127, {\em
Acad. Sci. Fennica}, Helsinki, 1980.

\bibitem{ef} F. Effenberger, Stacked polytopes and tight
triangulations of manifolds, {\em J. Combin. Theory} (A) {\bf 118}
(2011), 1843--1862.

\bibitem{fo}
R. Forman, Morse theory for cell complexes, {\em Adv. in Math.}
{\bf 134} (1998), 90--145.

\bibitem{f} M. Freedman, The topology of four dimensional
manifolds, {\em J. Diff. Geom.} {\bf 17} (1982), 357--454.

\bibitem{g}
B. Gr\"{u}nbaum, {\em Convex Polytopes} - 2nd ed. (GTM 221),
Springer-Verlag, New York, 2003.

\bibitem{ka}
G. Kalai, Rigidity and the lower bound theorem 1, {\em Invent.
math.} {\bf 88} (1987), 125--151.

\bibitem{kn}
S. Klee, I. Novik, Centrally symmetric manifolds with few
vertices, {\em Adv. in Math.} {\bf 229} (2012), 487--500.

\bibitem{k}
V. Klee, A combinatorial analogue of Poincar\'e's duality theorem,
{\em Can. J. Math.} {\bf 16} (1964), 517--531.

\bibitem{k86}
W. K\"{u}hnel, Higher dimensional analogues of Cs\'{a}sz\'{a}r's
torus, {\em Results in Mathematics} {\bf 9} (1986) 95--106.

\bibitem{k95}
W. K\"{u}hnel, {\em Tight Polyhedral Submanifolds and Tight
Triangulations}, Lecture Notes in Mathematics {\bf 1612},
Springer-Verlag, Berlin, 1995.

\bibitem{kb}
W. K\"{u}hnel, T. F. Banchoff, The 9-vertex complex projective
plane, {\em Math. Intelligencer} {\bf 5} (3) (1983), 11--22.

\bibitem{kl}
W. K\"{u}hnel, F. Lutz, A census of tight triangulations, {\em
Period. Math. Hunger.} {\bf 39} (1999), 161--183.

\bibitem{li}
W. B. R. Lickorish, Simplicial moves on complexes and manifolds,
{\em Geometry \& Topology Monographs}, {\bf 2} (1999), 299--320.
 maths.warwick.ac.uk/gt/GTMon2/paper16.abs.html

\bibitem{lu1}
F. H. Lutz, {\em Triangulated Manifolds with Few Vertices and
Vertex-Transitive Group Actions}, Thesis (D 83, TU Berlin),
Shaker Verlag, Aachen, 1999.

\bibitem{lu2}
F. H. Lutz, Triangulated manifolds with few vertices\,:
Combinatorial manifolds, arXiv:math/0506372v1, 2005, 37 pages.


\bibitem{lss}
F. H. Lutz, T. Sulanke, E. Swartz, $f$-vector of
3-manifolds, {\em Electron. J. Comb.} {\bf 16} (2009),
\#R\,13, 1--33.

\bibitem{mc1}
P. McMullen, On simple polytopes, {\em Invent. Math.} {\bf 113}
(1993), 419--444.

\bibitem{mc2}
P. McMullen, Triangulations of simplicial polytopes, {\em
Beitr\"{a}ge Algebra Geom.} {\bf 45} (2004), 37--46.

\bibitem{mw}
P. McMullen, D. W. Walkup, A generalized lower bound conjecture
for simplicial polytopes, {\em Mathematika} {\bf 18} (1971),
264--273.

\bibitem{ns} I. Novik, E. Swartz, Socles of Buchsbaum modules,
complexes and posets, {\em Adv. in Math.} {\bf 222} (2009),
2059--2084.

\bibitem{p}
U. Pachner, Konstruktionsmethoden und das kombinatorische
Hom\"{o}omorphieproblem f\"{u}r Triangulationen kompakter
semilinearer Mannigfaltigkeiten, {\em Abh. Math. Sem. Univ.
Hamburg} {\bf 57} (1987), 69--86.

\bibitem{r} G. Ringel, {\em Map color theorem}, Springer-Verlag,
New York - Heidelberg, 1974.

\bibitem{si}
N. Singh, Minimal triangulation of $(S^{\hspace{.2mm}3} \times
S^{\hspace{.2mm}1})^{\#3}$ (preprint).

\bibitem{st}
R. P. Stanley, The number of faces of simplicial polytopes and
spheres, {\em Annals New York Academy of Sciences} {\bf 440}
(1985), 212--223.

\bibitem{s} E. Swartz, Face enumeration - from spheres to manifolds,
{\em J. Eur. Math. Soc.} {\bf 11} (2009), 449--485.

\bibitem{wa}
D. W. Walkup, The lower bound conjecture for 3- and 4-manifolds,
{\em Acta Math.} {\bf 125} (1970) 75--107.

\bibitem{z} G. M. Ziegler, {\em Lectures on Polytopes},
Springer-Verlag, New York, 1995.

\bibitem{z2} G. M. Ziegler, Shelling polyhedral 3-balls and
4-polytopes, {\em Discrete Comput. Geom.} {\bf 19} (1998), 159--174.

\end{thebibliography}
}
\end{document}